\documentclass[lettersize,journal]{IEEEtran}
\usepackage{cite}
\usepackage{graphicx}
\usepackage{amsmath,amsfonts}
\usepackage[linesnumbered,ruled,vlined]{algorithm2e}
\usepackage{hyperref}
\usepackage{amssymb}
\usepackage{stfloats}
\usepackage{subfig}
\usepackage{tabularx}
\usepackage{multirow}    
\usepackage{booktabs}
\usepackage{diagbox}
\usepackage{makecell}
\usepackage{caption}
\usepackage{color} 
\usepackage[font=footnotesize]{caption}

\newtheorem{remark}{Remark}
\newtheorem{assumption}{Assumption}

\newtheorem{corollary}{Corollary}

\newtheorem{lemma}{Lemma}
\newtheorem{theorem}{Theorem}
\newtheorem{proof}{Proof}
\newtheorem{example}{Example}

\begin{document}

\title{Resilient Decentralized Preconditioned Stochastic Proximal Optimization with Gradient Sketching over Time-Varying Networks}
\author{Jinhui Hu, Guo Chen,~\IEEEmembership{Member,~IEEE}, Huaqing Li,~\IEEEmembership{Senior Member,~IEEE}, Tao Han,~\IEEEmembership{Member,~IEEE},\\ and Tingwen Huang,~\IEEEmembership{Fellow,~IEEE}
\thanks{This work is supported in part by National Natural Science Foundation of China under Grants 62603244, 62473136, 62673401, and 62573364, in part by Australian Research Council under Grants FT190100156 and DP230100801, in part by the Scientific Research Project of the Hubei Provincial Department of Education under Grant Q20252504, in part by Hubei Normal University Research Funds for Talent Cultivation under Grant YC2026124, and in part by the Chongqing Science Fund for Distinguished Young Scholars under Grant 2024NSCQ-JQX0103.
}
\thanks{J. Hu and T. Han are with the School of Artificial Intelligence and Computer Science, Hubei Normal University, Huangshi 435002, China, and J. Hu is also with the Department of Mechanical Engineering, City University of Hong Kong, Kowloon Tong, Kowloon, Hong Kong SAR, China (e-mail: jinhuihu3-c@my.cityu.edu.hk; taohan@hbnu.edu.cn).}
\thanks{G. Chen is with the School of Electrical Engineering and Telecommunications, University of New South Wales, Sydney, NSW 2052, Australia (e-mail: guo.chen@unsw.edu.au).}
\thanks{H. Li is with the Chongqing Key Laboratory of Nonlinear Circuits and Intelligent Information Processing, College of Electronic and Information Engineering, Southwest University, Chongqing 400715, China (e-mail: huaqingli@swu.edu.cn).}
\thanks{T. Huang is with the School of Computer Science and Control Engineering, Shenzhen University of Advanced Technology, Shenzhen 518055, China (email: huangtingwen@suat-sz.edu.cn).}
}

\markboth{}%
{Shell \MakeLowercase{\textit{et al.}}: A Sample Article Using IEEEtran.cls for IEEE Journals}


\maketitle

\begin{abstract}\label{sec1}
Variance reduction is indispensable in Byzantine-resilient decentralized stochastic optimization over multi-agent systems (MASs), as it effectively mitigates gradient noise to facilitate the resilient aggregation process. However, most existing Byzantine-resilient decentralized variance-reduced (VR) stochastic gradient algorithms rely on random data sampling, which proves inefficient in data-scarce yet high-dimensional tasks, for instance, image deblurring. This paper aims to pursue an alternative approach that achieves variance reduction via gradient sketching. We first formulate a class of structural risk minimization (SRM) problems over time-varying networks, where the local objectives are not necessarily decomposable and their gradients may be inaccessible. To solve the SRM problems in a decentralized manner, we incorporate a localized gradient-sketching technique (i.e., SEGA with random dimension sampling) into a preconditioned stochastic proximal gradient framework, proposing a novel gossip-based decentralized algorithm with variance reduction, dubbed \textit{Gossip-SEGA}. To further enhance its resilience, we substitute the weighted-average aggregation in \textit{Gossip-SEGA} with a norm-penalized approximation, thereby developing a Byzantine-resilient decentralized algorithm, namely \textit{RED-SEGA}. Theoretically, we derive sufficient conditions for both consensus among reliable agents and (sub)linear convergence rate of \textit{RED-SEGA} over time-varying networks in the presence of Byzantine agents with unknown identities.
This is challenging, as the convergence and consensus errors must be bounded within a transformed variable-metric geometry, where the inherent complexities are further compounded by adversarial Byzantine agents and time-varying networks. The effectiveness and resilience of \textit{RED-SEGA} are validated by solving a constrained least-squares problem and an image deblurring task under various Byzantine attacks.
\end{abstract}


\section{Introduction}\label{sec2}
\subsection{Literature Review}\label{sec2-1}
\IEEEPARstart{D}{ecentralized} optimization presents significant advantages for both cooperative and noncooperative tasks in MASs, as it enhances privacy protection during information exchange \cite{Wang2024}, addresses the communication bottlenecks inherent in centralized parameter-server architectures \cite{Lian2017a}, and boosts scalability for large-scale MAS deployments \cite{Shorinwa2024}. Benefiting from these advantages, decentralized optimization algorithms (DOAs) have garnered considerable attention in collaborative machine/deep learning \cite{Ye2025,Huang2025}, noncooperative games \cite{Li2025a}, and energy coordination \cite{Li2021}.

While many elegant DOAs \cite{Li2021c,Makridis2024,Doostmohammadian2025,Lu2026} exhibit strong efficacy and robustness for large-scale optimization tasks under the assumption of universal reliability, they remain vulnerable to the disruptive influence of malfunctioning or malicious agents. Such agents are known as the Byzantine agents, arising in the course of multi-agent optimization and learning. They may either unintentionally disseminate untrue or misleading information, or deliberately execute various Byzantine (backdoor) attacks, often operating with an omniscient viewpoint \cite{Chang2023,Hu2025,Peng2025}. For instance, Byzantine agents can perform label-flipping attacks via changing the label of its local training data to disturb the aggregation process in federated learning (FL) \cite{Karimireddy2021}. With the omniscient viewpoint of both node- and network-level information, Byzantine agents are also able to launch dissensus attacks via delivering misleadingly untrue models to different reliable neighbors, thereby isolating the local gradient descent step in decentralized learning \cite{He2022}. Prior studies \cite{Karimireddy2021,Chang2023,Hu2025,Peng2025} have demonstrated that Byzantine agents can severely undermine the utility of final trained models via diverse attack modes.

Data heterogeneity and anonymity of individual agents within decentralized systems present formidable obstacles to accurately detect Byzantine agents \cite{Karimireddy2021,Hu2025,Peng2025}, which leads to severe performance degradation of the whole MASs, for instance, substantial disagreement among reliable agents \cite{Sun2026,Duan2026}. Therefore, significant efforts have been dedicated to developing Byzantine-resilient optimization algorithms. For instance, within a centralized parameter-server architecture, standard stochastic gradient descent (SGD) is frequently integrated with resilient aggregation strategies, including geometric median \cite{Chen2017a}, coordinate-wise median \cite{Yin2018}, trimmed mean \cite{Xie2018}, Krum \cite{Blanchard2017}, FABA \cite{Xia2019}, norm-penalized approximation \cite{Li2019n}, centered clipping \cite{Karimireddy2021}, and best-effort voting \cite{Fan2022}. Nevertheless, an inherent limitation of the centralized parameter-server architecture is their dependence on an infallible central server for model or gradient aggregation, posing risks such as single-point failures and communication bottlenecks. To address these issues, extensive efforts have sought to extend Byzantine-resilient centralized methods \cite{Chen2017a, Yin2018, Xie2018, Blanchard2017, Xia2019, Li2019n, Karimireddy2021} to decentralized settings by exploring the connection between resilient aggregation strategies and gossip-based communication protocols, leading to many elegant Byzantine-resilient decentralized stochastic optimization algorithms \cite{Peng2021,He2022,Wu2022,Han2025}. However, the noise incurred by evaluating local stochastic gradients obscures malicious inputs from Byzantine agents and complicates the identification, filtering, screening, or penalization in the resilient aggregation, particularly exacerbated under non-independent and identically distributed (non-i.i.d.) data settings \cite{Peng2021,Liu2023c}.

To asymptotically eliminate stochastic gradient noise, an earlier work \cite{Wu2020} combines \textit{SAGA} with geometric median aggregation \cite{Chen2017a} to propose \textit{Byrd-SAGA} for centralized FL. A subsequent study \cite{Gorbunov2023} simultaneously investigates variance reduction and communication compression in Byzantine-resilient centralized optimization to devise \textit{Byz-VR-MARINA}, where tighter theoretical results are provided. In decentralized settings, a recent study \cite{Hu2025a}, aiming to solve a class of \text{SRM} problems in the presence of Byzantine agents, develops a proximal-gradient mapping algorithmic framework, termed \textit{Prox-DBRO-VR}, by integrating two VR techniques via random data sampling with a resilient aggregation inspired by \cite{Ben-ameur2016,Li2019n,Peng2021}. Both theoretically and empirically, \cite{Wu2020, Gorbunov2023, Hu2025a} demonstrate that VR techniques are an effective countermeasure against Byzantine agents.

While existing Byzantine-resilient stochastic gradient algorithms attain VR via random data sampling \cite{Wu2020, Gorbunov2023, Hu2025a}, this approach demonstrates high efficacy in large-scale data regimes. However, in data-scarce yet high-dimensional scenarios, for instance, image deblurring \cite{Yang2025}, these methods are often hampered by the curse of dimensionality, leading to low efficiency. More importantly, the application of these methods \cite{Wu2020, Gorbunov2023, Hu2025a} is premised on the condition that local objectives must be decomposable. This may not always hold in domain-specific problems, such as reinforcement learning \cite{Wu2022}, energy coordination \cite{Li2021}, and model predictive control \cite{Li2021b}. Moreover, in deep learning tasks, for instance, large language model (LLM) fine-tuning, the computation of full backpropagation gradients often becomes a major bottleneck due to its prohibitive cost \cite{Chen2023a,Chen2024a}. A promising remedy for these bottlenecks is offered by \textit{SEGA} \cite{Hanzely2018}, which leverages a gradient-sketching technique to flexibly enable both first-order (FO) stochastic and zeroth-order (ZO) derivative-free optimization by tailoring its sketching strategy. Notably, the FO version of \textit{SEGA} achieves variance reduction via random dimension sampling, thereby circumventing the need to compute full-batch gradients with possibly prohibitive cost, while the ZO version of \textit{SEGA} approximates the directional derivative by using the difference of function values and thus free of gradient computation. To further cut down the computational cost, follow-up work \cite{Gorbunov2020} introduces random data sampling into \textit{SEGA} to propose \textit{N-SEGA}.

\subsection{Motivations and Challenges}\label{sec2-2}
Both \textit{SEGA} \cite{Hanzely2018} and its extension \textit{N-SEGA} \cite{Gorbunov2020} are confined to single-agent optimization and lack scalability to large-scale MASs. However, directly extending them via popular \textit{gossip} communication to a decentralized setting renders the resulting algorithms highly vulnerable to Byzantine attacks \cite{Peng2021,He2022,Wu2022,Han2025,Hu2025a}. These gaps motivate us to design a Byzantine-resilient decentralized algorithm that achieves variance reduction while escaping from the curse of dimensionality in data-scarce yet high-dimensional scenarios, obviating the decomposability prerequisite for local objectives, and preserving its applicability to derivative-free optimization. Nevertheless, developing such an algorithm entails non-trivial theoretical hurdles, specifically in constructing an unbiased stochastic gradient estimator and in rigorously bounding the convergence and consensus behaviors within a transformed variable-metric geometrical complexities that are further compounded by adversarial Byzantine agents and time-varying networks.

\subsection{Contributions}\label{sec2-3}
This paper addresses the aforementioned challenges by extending a \textit{sketch-and-project} process to the decentralized domain, integrating the localized gradient-sketching technique into a preconditioned stochastic proximal gradient framework, and calibrating the algorithmic parameters. This rigorous calibration secures statistical unbiasedness while balancing the gradient estimation error, gradient-learning error, and resilient consensus conditions. The main contributions of this work are threefold:
\begin{enumerate}
  \item This paper investigates a class of SRM problems over time-varying networks in the presence of Byzantine agents. To address the challenges posed by non-decomposable local objectives and potentially inaccessible exact gradients, we first extend \textit{SEGA} with random dimension sampling \cite{Hanzely2018} to a gossip-based decentralized variant, termed \textit{Gossip-SEGA}. Recognizing its inherent vulnerability to Byzantine attacks, we further develop \textit{RED-SEGA}, a Byzantine-resilient decentralized preconditioned stochastic proximal gradient algorithm with variance reduction. This resilient extension is achieved by replacing the naive weighted-average consensus in \textit{Gossip-SEGA} with a resilient aggregation mechanism based on a norm-penalized approximation.
  \item In contrast to existing Byzantine-resilient methods \cite{Wu2020, Gorbunov2023, Hu2025a} that achieve variance reduction via random data sampling, \textit{RED-SEGA} attains variance reduction through a different technical line, i.e., gradient sketching. This strategy is particularly advantageous for data-scarce yet high-dimensional tasks, as only few partial derivatives of local objectives are evaluated at each iteration. The improved computational efficiency of \textit{RED-SEGA} compared to the baseline methods is validated by extensive experiments in Section \ref{sec6-1}.
  \item Both theoretically and empirically, \textit{RED-SEGA} is proven to work effectively over time-varying networks, even with the existence of Byzantine agents. This overcomes a major limitation of most existing Byzantine-resilient DOAs \cite{Han2025,He2022,Wu2022,Yang2019c,Fang2022,Hu2025,Hu2025a,Wang2024}, which are only confined to static networks.
\end{enumerate}

\subsection{Organization}\label{sec2-4}
The remainder of this paper is organized as follows. Section \ref{sec3} introduces the necessary preliminaries. Section \ref{sec4} details the algorithmic development of the proposed algorithms. The theoretical consensus and convergence analyses are established in Section \ref{sec5}. Section \ref{sec6} verifies the effectiveness and resilience of the proposed algorithms through extensive numerical experiments. The detailed mathematical derivations have been deferred to the Supplementary Material due to space constraints.

\section{Preliminaries}\label{sec3}
\subsection{Basic Notation}\label{BasNo}
All vectors in this paper are column vectors unless stated otherwise.
\begin{table}[!h]
\centering
\caption{Basic notations.}
\begin{tabularx}{8.8cm}{lX}  
\toprule                
\bf{Symbols}  & \bf{Definitions}  \\
\midrule
${\mathbb R}$, ${{\mathbb R}^n}$, ${{\mathbb R}^{m \times n}}$ & the sets of real numbers, $n$-dimensional column real vectors, $m \times n$ real matrices, respectively\\
:= & the definition symbol\\
${\cdot ^ \top }$ & the transpose of any matrices or vectors \\
$\left| \cdot \right|$ & an operator to represent the absolute value of a constant or the cardinality of a set\\
$\left\|  \cdot   \right\|$ & the standard Euclidean norm for vectors or spectral norm for matrices\\
${\left\| x \right\|_p}$ & the $\ell_p$-norm of a vector $x \in \mathbb{R}^n$, defined as ${\left( \sum\nolimits_{i=1}^n |x_i|^p \right)^{1/p}}$ for a scalar $p \ge 1$. \\
${\langle x, y \rangle}_W$ & the $W$-weighted inner product of vectors $x, y \in \mathbb{R}^n$, defined as $x^\top W y$ for a positive-definite matrix $W \succ 0$. \\
${\left\| x \right\|_W}$ & the $W$-weighted Euclidean norm of a vector $x \in \mathbb{R}^n$, defined as $\sqrt{x^\top W x}$ for a positive-definite matrix $W \succ 0$, or its induced matrix norm. \\
$0_n$ & an n-dimensional all-zero vector \\
$1_n$ & an n-dimensional all-one vector \\
${\mathbf{I}}$ & an identity matrix of appropriate dimension\\
\bottomrule
\end{tabularx}
\label{Tab1}
\end{table}
For arbitrary three vectors $\tilde x, \tilde y, \tilde z\in \mathbb{R}^{ n}$, we consider a positive scalar $a$ and a closed, proper, convex function $r:{\mathbb{R}^n} \to \mathbb{R} $, the proximal operator of $g$ is defined by
${\mathbf{prox}}_{{a, r}}\left \{ \tilde x \right \} = \arg {\min _{\tilde y}} \left \{ {r\left( \tilde y \right) + \frac{1}{{2a}}{{\left\| {\tilde y - \tilde x} \right\|}^2}} \right \}$; With a slight abuse of notation,
let $\partial r\left( \tilde x \right)$ denote the subdifferential or subgradient of the proper, closed, and convex function $r:{\mathbb{R}^{ n}} \to \mathbb{R} $ at $\tilde x$, such that
\begin{equation*}
\partial r\left( {\tilde x} \right) = \left\{ {\tilde y|\forall \tilde z \in {\mathbb{R}^n}, r\left( {\tilde z} \right) \ge r\left( {\tilde x} \right) + \left\langle {\tilde y,\tilde z - \tilde x} \right\rangle } \right\}.
\end{equation*}
${\lambda _{\min }}\left( X \right) $ and ${\lambda _{\max }}\left( X \right)  $ denote the maximum and minimum eigenvalue of any matrix $X$, respectively. Some basic notaions are defined in Table \ref{Tab1}.

\subsection{Problem Formulation}\label{sec3-1}
 This paper aims to solve a class of SRM problems over a networked MAS comprising finite numbers of reliable and Byzantine agents, which is formulated as follows:
\begin{equation}\label{E3-1-1}
\text{SRM}: \quad  \mathop {\min }\limits_{\tilde x} \sum\limits_{i \in {\mathcal R}} {{f_i}\left( {\tilde x} \right) + r_i\left( {\tilde x} \right)},
\end{equation}
where $\tilde x \in {\mathbb{R}^n}$ is the decision variable, ${f_i}:{\mathbb{R}^n} \to \mathbb{R}$ and $r_i:{\mathbb{R}^n} \to \mathbb{R} \cup \left\{ { + \infty } \right\}$ are the smooth and (possibly) nonsmooth local objective function private to agent $i$, $i \in {\mathcal R}$. Let ${\mathcal R}$ and $\mathcal{B}$ denote the sets of reliable agents and Byzantine agents in the network, respectively. In this paper, we consider a severe Byzantine threat model where Byzantine agents are omniscient and completely unconstrained in their behavior. Specifically, they possess complete a \textit{priori} knowledge of the entire MASs, including network topologies, reliable agents' instantaneous states, as well as all algorithmic parameters. They can collude and seamlessly adapt their attacks over time without any bounds. Notably, despite this omniscience, their disruptive capabilities are limited to falsifying their own broadcasted states, as they cannot tamper with the local update rules of other agents.
We denote the optimal solution of (\ref{E3-1-1}) by $\tilde x ^*$. Unlike existing decentralized stochastic gradient algorithms, such as \cite{Xin2022,Ye2025,Hu2025a}, positing a decomposable structure for the local objectives, this paper does not require ${f_i}$, ${i \in {\mathcal R}}$, to admit any further decomposition, aligning with numerous practical applications, such as reinforcement learning \cite{Wu2022}, energy coordination \cite{Li2021}, and model predictive control \cite{Li2021b}. In fact, the true gradients of ${f_i}$, ${i \in {\mathcal R}}$, are not necessarily to be available during iterations (see Section \ref{sec4-2}), satisfying with the principles of ZO optimization that finds extensive applications in black-box training \cite{Chen2023a} and LLM fine-tuning \cite{Chen2024a} recently. To facilitate the sequel convergence analysis, we need to make the following assumptions. We consider positive-definite weight matrices $W_i$, $\forall i \in \mathcal{R}$, whose exact characterization will be provided in Section \ref{sec4-2}.
\begin{assumption}(Strong convexity)\label{A1}
For $i \in \mathcal{R}$, the local objective function $f_i$ is $\mu_i$-strongly convex in $W_i$-norm, meaning that for all $\tilde x$ and $\tilde y$, we have
\begin{equation}\label{E3-1-2}
  {f_i}\left( {\tilde x} \right) \ge {f_i}\left( {\tilde y} \right) + {\left\langle {\nabla {f_i}\left( {\tilde y} \right),\tilde x - \tilde y} \right\rangle} + \frac{{{\mu _i}}}{2}\left\| {\tilde x - \tilde y} \right\|_{W_i}^2,
\end{equation}
where $\mu_i \ge 0$ is the strongly-convex parameter.
\end{assumption}

\begin{assumption}(General smoothness)\label{A2}
For $i \in \mathcal{R}$, the local objective function $f_i$ is $Q_i$-smooth with respect to $W_i$, meaning that for all $\tilde x$ and $\tilde y$, we have
\begin{equation}\label{E3-1-3}
    {f_i}\left( {\tilde x} \right) - {f_i}\left( {\tilde y} \right) - {\left\langle {\nabla {f_i}\left( {\tilde y} \right),\tilde x - \tilde y} \right\rangle} \ge \frac{\left\| {\nabla {f_i}\left( {\tilde x} \right) - \nabla {f_i}\left( {\tilde y} \right)} \right\|_{{Q_i}}^2}{2},
\end{equation}
where ${Q_i} \succ 0$ is a positive-definite matrix.
\end{assumption}
\begin{remark}\label{R1}
Assumption \ref{A1} ensures that the optimal solution to problem (\ref{E3-1-1}) exists uniquely. If we set $W_i = {\mathbf{I}}$ for all $i \in \mathcal{R}$, then Assumption \ref{A1} reduces to the standard strongly-convex condition assumed in most existing literature, such as \cite{Yang2019c,Li2021b,Peng2021,Fang2022,Chang2023,Hu2025a}. It follows from \cite{Hanzely2018} that Assumption \ref{A2} reduces to a standard Lipschitz continuous condition under certain conditions. Therefore, Assumptions \ref{A1}-\ref{A2} are two flexible variants of assumptions commonly adopted in \cite{Yang2019c,Li2021b,Peng2021,Fang2022,Chang2023,Hu2025a}.
\end{remark}

\subsection{Time-Varying Network Model}\label{sec3-2}
We consider a time-varying network ${{\mathcal{G}}_k} = \left( {{\mathcal{V}},{\mathcal{E}}\left( {{\zeta _k}} \right)} \right)$, where ${\mathcal{V}}$ denotes the set of all agents including both reliable and Byzantine agents in the network, and ${\mathcal{E}}\left( {{\zeta _k}} \right)$ represents the set of connected (communication) edges at iteration (time) $k$, where the random vector ${\zeta _k}$ is defined by ${\zeta _k}: = {\mathbf{col}}{\left\{ {{{\tilde \zeta }_{\left( {i,j} \right),k}}} \right\}_{\left( {i,j} \right) \in {\mathcal E}}}$ with ${\bar \zeta _{\left( {i,j} \right)}} = \left\{ \begin{array}{l}
 1,\left( {i,j} \right) \in {\mathcal{E}_k},\\
 0,{\rm{otherwise.}}
\end{array} \right.$, being a random vector associated with the connectivity of the edge ${\left( {i,j} \right)}$ at iteration $k \ge 0$. Specifically, if agent $i$ is connected with agent $j$ at iteration $k$, i.e., $\left( {i,j} \right) \in {{\mathcal{E}}_k}$, then ${{{\bar \zeta }_{\left( {i,j} \right),k}}} = 1$; otherwise, ${{{\bar \zeta }_{\left( {i,j} \right),k}}} = 0$.
We denote by ${\mathcal{E}}\left( \zeta  \right)$ the set of candidate connected edges collected throughout iterations, i.e., ${\mathcal{E}}: = \bigcup\nolimits_{k = 0}^\infty  {{\mathcal{E}_k}} $. For any two different edges $\left\{ {\left( {i,j} \right), \left( {i',j'} \right)} \right\} \in {\mathcal{E}}$,
the entry ${\bar \zeta _{\left( {i,j} \right),k}}$ of ${\zeta _k}$ is independent of ${\bar \zeta _{\left( {i',j'} \right),k}}$ at any iteration $k$. However, ${\zeta _k}$ is not necessarily independent of each other across iterations. Building on this, there exists a countably infinite random sequence ${\left\{ {{\zeta _k}} \right\}_{k \ge 0}}$ characterizing the time-varying network during the optimization process. To simplify the notations, we use ${{\mathcal{E}}_k}$ to briefly denote ${\mathcal{E}}\left( {{\zeta _k}} \right)$. At each iteration $k$, let ${{\mathcal{N}}_{i,k}}$, ${{\mathcal{R}}_{i,k}}$, and ${{\mathcal{B}}_{i,k}}$ denote the instantaneous sets of total, reliable, and Byzantine neighbors of agent $i$, respectively, such that ${{\mathcal{N}}_{i,k}} = {{\mathcal{R}}_{i,k}} \cup {{\mathcal{B}}_{i,k}}$. Note that the sets ${{\mathcal{N}}_{i,k}}$, ${{\mathcal{R}}_{i,k}}$, and ${{\mathcal{B}}_{i,k}}$ exclude agent $i$. To guarantee sufficiently frequent information aggregation among all reliable agents, we will introduce an \textit{average} network based on ergodicity \cite{Walters1982a}. Specifically, we denote ${\omega _{\left( {i,j} \right)}}$ by the limiting connected frequency associated with the reliable edge $\left( {i,j} \right)$ such that for $i,j \in {\mathcal{R}}$ and $i < j$, we have
\begin{equation}\label{E3-2-1}
{\omega _{\left( {i,j} \right)}} = \mathop {\lim \inf }\limits_{T \to \infty } \frac{1}{{T}}\sum\limits_{l = 0}^{T-1} {{{\bar \zeta }_{\left( {i,j} \right),l}}}.
\end{equation}
Therefore, this ergodic property allows us to define an average edge set as $\bar{\mathcal{E}}_{\mathcal{R}} := \{(i,j) \mid \omega_{(i,j)} > 0\}$, which forms the expected average network among all reliable agents, denoted by $\bar{\mathcal{G}}_{\mathcal{R}} := (\mathcal{R}, \bar{\mathcal{E}}_{\mathcal{R}})$. We next introduce an assumption regarding the time-varying network topology to rigorously guarantee sufficient connectivity and support the subsequent non-asymptotic convergence analysis.
\begin{assumption}(Ergodic connectivity)\label{A3}
The sequence of random network topologies ${\left\{ {{\zeta _k}} \right\}_{k \ge 0}}$ is independent and identically distributed (i.i.d.) over time. For any reliable edge $(i, j)$, its connection probability is strictly positive, given by $\mathbb{E}\left[ {{\zeta_{\left( {i,j} \right),k}}} \right] = {\omega _{\left( {i,j} \right)}} > 0$. Furthermore, the expected network graph $\bar {\mathcal{G}}_{\mathcal{R}}$ is bidirectionally connected.
\end{assumption}
\begin{remark}\label{R2}
Assumption \ref{A3} guarantees the ergodic connectivity for time-varying networks, which relaxes the persistent connectivity condition mandated by many existing algorithms over static networks \cite{Fang2022,Yang2019c,He2022,Hu2025a}. Specifically, Assumption \ref{A3} accommodates temporary disconnections or disruptions of reliable communication edges, provided that the expected graph ensures sufficient information aggregation among the reliable agents. The adopted i.i.d. random graph model practically encompasses a variety of stochastic network scenarios, such as dynamic Erdős-Rényi topologies, probabilistic link failures, and networks with randomly activated edges \cite{Nedic2017a,Saadatniaki2020,Chen2025,Nguyen2025,Han2025}. It is worth mentioning that while this paper primarily models intermittent connectivity through time-varying networks, practical MASs inevitably suffer from communication delays. As demonstrated in a previous work \cite{Hu2024}, the adopted norm-penalized aggregation inherently provides a certain degree of tolerance against communication delays by filtering them in a manner akin to malicious Byzantine perturbations, ensuring that the network structure and system stability are not disrupted. Nevertheless, simply penalizing delayed information as malicious faults undermines the convergence performance of the proposed algorithm. Albeit beyond the scope of this paper, quantifying the theoretical bounds of communication delays and exploring delay-tolerant techniques \cite{Makridis2024,Doostmohammadian2025} to restore optimal convergence performance constitute an interesting direction for future research.
\end{remark}

\subsection{Resilient Consensus Problem Setup}\label{sec3-3}
We introduce a global decision variable $x: = {\mathbf{col}}{\left\{ {{x_i}} \right\}_{i \in \mathcal{R}}} \in {\mathbb{R}^{n\left| \mathcal{R} \right|}}$, which collects sequentially local decision variables of all reliable agents. Building upon the time-varying network model defined in Section \ref{sec3-2}, \text{SRM} (\ref{E3-1-1}) can be equivalently formulated as a decentralized constrained consensus (\text{DCC}) problem subject to (s.t.) an equality constraint ${x_i} = {x_j}$, $\left( {i,j} \right) \in {{ \mathcal{E}}_{\mathcal{R},k}}$, which is described as follows:
\begin{equation}\label{E3-3-1}
\begin{aligned}
\text{DCC}: \quad &\mathop {\min }\limits_x \sum\limits_{i \in \mathcal{R}} {{f_i}\left( {{x_i}} \right) + r_i\left( {{x_i}} \right)},\\
&{\text{s.t.}} \; {x_i} = {x_j}, \left( {i,j} \right) \in {{ \mathcal{E}}_{\mathcal{R},k}}.
\end{aligned}
\end{equation}
In Byzantine-free cases, one popular line of existing decentralized algorithms addresses the consensus constraint ${x_i} = {x_j}$, $\left( {i,j} \right) \in {{ \mathcal{E}}_{\mathcal{R},k}}$ in (\ref{E3-3-1}) via a family of contractive stochastic matrices \cite{Nedic2017a,Saadatniaki2020,Chen2025,Nguyen2025} to ensure consensual information aggregation. However, it has been shown that this kind of consensus protocols will fail when there are Byzantine agents in the network \cite{He2022,Hu2025a}, which will be discussed in Section \ref{sec4-2}. For the subsequent analysis, we let $\bar{\mathcal{R}}_i = \{j \mid (i,j) \in \bar{\mathcal{E}}_{\mathcal{R}}\}$ denote the set of reliable neighbors of agent $i$ over the \textit{average} network $\bar {\mathcal{G}}_{\mathcal{R}}$.

Building on existing remedies \cite{Ben-ameur2016,Hu2025a}, a decentralized resilient consensus (\text{DRC}) problem is formulated as follows:
\begin{equation}\label{E3-3-2}
\text{DRC}:\mathop {\min }\limits_x \sum\limits_{i \in \mathcal{R}} { {{f_i}\left( {{x_i}} \right) + r_i\left( {{x_i}} \right) + \phi \sum\limits_{i \in \mathcal{R}} {\sum\limits_{j \in {{\bar {\mathcal{R}}}_i}} {{\omega _{\left( {i,j} \right)}}{{\left\| {{x_i} - {x_j}} \right\|}_p}} }   } },
\end{equation}
where $\phi$ is a penalty parameter. \text{DRC} (\ref{E3-3-2}) is a soft-approximation variant of \text{DCC} (\ref{E3-3-1}), typically allowing some pairs of neighboring states to differ throughout the optimization process. This is an essential property for MASs when Byzantine agents are present \cite{Peng2021,Hu2025a}, as enforcing the hard consensus constraint in (\ref{E3-3-1}) would allow Byzantine agents to compromise the system. To formally quantify the systemic disagreement among reliable agents, we define the deterministic expected consensus penalty function over the \textit{average} network $\bar{\mathcal{G}}_{\mathcal{R}}$ as $\bar{P}(x) := \phi \sum_{i\in\mathcal{R}} \sum_{j\in\bar{\mathcal{R}}_i} \omega_{(i,j)} \|x_i - x_j\|_p$. Let $x_i^*$ be the local optimal solution for agent $i$, and define the global optimal solution vector to the DRC problem (\ref{E3-3-2}) as ${x^*}: = {\mathbf{col}}{\left\{ {x_i^*} \right\}_{i \in \mathcal{R}}}$.

\section{Algorithm Development}\label{sec4}
\subsection{Connection to Existing Work}\label{sec4-1}
Most existing Byzantine-resilient decentralized vanilla SGD algorithms \cite{Yang2019c,Peng2021,He2022,Ye2024b,Hu2025} reduce the per-iteration gradient computation complexity for each agent compared to decentralized batch gradient methods \cite{Sundaram2019,Yang2019c,Fang2022}, but they also introduce an accumulated convergence error when evaluating local stochastic gradients. Therefore, inspired by VR techniques \cite{Gorbunov2023}, a recent work \cite{Hu2025a} proposes an algorithmic framework, dubbed \textit{Prox-DBRO-VR}, via incorporating localized VR techniques, \textit{SAGA} \cite{Defazio2014c} and \textit{LSVRG} \cite{Kovalev2019}, into a Byzantine-resilient decentralized algorithm \cite{Sundaram2019}, to achieve linear convergence without incurring the accumulated convergence error. However, this approach assumes the local objectives to be decomposable. This assumption is not applicable to some engineering optimization problems, such as decentralized model predictive control \cite{Li2021b}, decentralized energy management \cite{Chang2023}, and decentralized reinforcement learning \cite{Wu2022}, exhibiting a non-decomposable structure of local objectives.
Moreover, this paper seeks Byzantine resilience over time-varying networks, which is significantly weaker than the persistent connectivity assumption mandated by many Byzantine-resilient works operating over static networks (e.g., \cite{Yang2019c,Li2021b,Fang2022,He2022,Chang2023,Hu2025a,Hu2025}). Specifically, it accommodates temporary disconnections or disruptions of reliable communication edges, provided that information aggregation remains sufficiently frequent over time. This relaxation is beneficial to the MASs with the existence of Byzantine agents as it eliminates the need for continuously fixed connectivity among all reliable agents, required by \cite{Wu2022,Hu2025a,Fang2022,Ye2024b,Han2025}, given that the identities of these reliable agents are not known a \textit{priori}. Furthermore, conventional Byzantine-resilient decentralized methods (e.g., \cite{Yang2019c,Peng2021,Fang2022,He2022,Wu2022,Ye2024b,Han2025,Hu2025,Hu2025a}) predominantly depend on standard FO gradients. In contrast, \textit{RED-SEGA} admits a versatile gradient sketching oracle that inherently supports both FO and ZO optimization. This ZO capability proves advantageous for complex tasks where exact gradient evaluation is computationally prohibitive or even inaccessible, notably in black-box training \cite{Chen2023a} and LLM fine-tuning \cite{Chen2024a}.

\subsection{Localized Gradient Sketching}\label{sec4-2}
Each reliable agent $i$ is assumed to be able to query an $\mathbb{S} $-oracle that returns a random linear transformation of the gradient with respect to the local objective function $f_i$ rather than the true gradient $\nabla {f_i}$, $i \in \mathcal{R}$. We define the sketched matrices by $S_i \in {\mathbb{R}^{n \times m}}$ ($m \ge 1$) and a fixed local distribution $\mathcal{D}_i$ over sketching matrices $S_i \in {\mathbb{R}^{n \times p}}$, $\forall i \in \mathcal{R}$. Considering a query point $z \in {\mathbb{R}^n}$, the $\mathbb{S} $-oracle returns a random linear transformation of the gradient, which is defined by
\begin{equation}\label{E4-2-1}
\mathbb{S} \left( {S_i,z} \right): = {S_i^ \top }\nabla {f_i}\left( z \right).
\end{equation}
where $S_i \sim \mathcal{D}_i$. There are two specific examples of the $\mathbb{S} $-oracle defined in (\ref{E4-2-1}): i) If the local distribution $\mathcal{D}_i$ is set uniformly over standard unit basis vectors $e_i$, $i \in \mathcal{R}$, such that $\mathbb{S}\left( {{e_i},z} \right) = e_i^ \top \nabla {f_i}\left( z \right)$ represents the $i$-th partial derivative of the local objective function $f_i$ at $z$, then $\mathbb{S} $-oracle provides the gradient estimate defined by the decentralized randomized coordinate descent algorithm \cite{Richtarik2016};
ii) If the local sketched matrix reduces to a vector, i.e., $S_i = s_i \in {\mathbb{R}^n}$, then one can approximate the directional derivatives by applying the difference of function values, i.e., $\mathbb{S}\left( {{s_i},z} \right) \approx \left( {{f_i}\left( {z + \xi {s_i}} \right) - {f_i}\left( z \right)} \right)/\xi $ with $\xi > 0$, which implies that $\mathbb{S} $-oracle can be applied to decentralized ZO optimization \cite{Akhavan2023,Yuan2024}. However, it is worth noting that while the subsequent theoretical analysis assume exact and unbiased directional derivatives, the finite difference approximation used in the ZO implementation inherently introduces a bias. Therefore, we explicitly present the ZO method as a heuristic extension of our \textit{RED-SEGA} framework, whose empirical effectiveness will be demonstrated in the experimental section.

It is worth emphasizing that the gradient-sketching oracle $\mathbb{S} $ defined in (\ref{E4-2-1}) is not required to return an unbiased estimator of the local gradient $\nabla {f_i}$, in contrast to unbiased stochastic oracles commonly assumed in decentralized FO algorithms \cite{Yang2019c,Fang2022,He2022,Wu2022,Ye2024b,Hu2025}. Nevertheless, it remains necessary to construct an unbiased stochastic gradient estimator using certain techniques, as such an estimator is critical for convergence analysis. We next demonstrate how this can be accomplished in two steps. \newline
\textbf{Step I:} For any $k \ge 0$, we consider the current status $x_{i,k}$ and define the current estimate of the local gradient $\nabla {f_i}\left( x_{i,k} \right)$ by ${h _{i,k}}$, $\forall i \in \mathcal{R}$. In view of the gradient-sketching step (\ref{E4-2-1}), we extend a \textit{sketch-and-project} process \cite{Arik2020,Hanzely2018} to a decentralized domain such that the next estimate ${h_{i,k+1} }$ can be calculated as follows:
\begin{equation}\label{E4-2-2}
{h _{i,k + 1}} = {h _{i,k}} + W_i^{ - 1}{M_{i,k}}\left( {\nabla {f_i}\left( {{x_{i,k}}} \right) - {h _{i,k}}} \right),
\end{equation}
where ${M_{i,k}}: = {S_{i,k}}{\left( {{{\left( {{S_{i,k}}} \right)}^ \top }W_i^{ - 1}{S_{i,k}}} \right)^\dag }{\left( {{S_{i,k}}} \right)^ \top }$ is a symmetric matrix with ${S_{i,k}}$ being a time-varying sketched matrix over $k \ge 0$, $\forall i \in \mathcal{R}$. According to pseudo-inverse properties, one can verify that
\begin{equation}\label{E4-2-2+}
{\left( {{M_{i,k}}} \right)^ \top }W_i^{ - 1}{M_{i,k}} = {M_{i,k}}.
\end{equation}
The estimate ${h _{i,k + 1}}$ is produced by seeking the close-form solution of the sequel constrained convex optimization problem
\begin{equation}\label{E4-2-3}
\begin{aligned}
  &\mathop {\min }\limits_{h} \left\| {h  - {h _{i,k}}} \right\|_{{W_i}}^2, \\
  &{\text{s}}{\text{.t}}{\text{. }}S_{i,k}^ \top h  = S_{i,k}^ \top \nabla {f_i}\left( {{x_{i,k}}} \right). \\
\end{aligned}
\end{equation}
\textbf{Step II:} For any $k \ge 0$, ${h _{i,k + 1}}$ is not guaranteed to be an unbiased estimator of $\nabla {f_i}\left( {{x_{i,k}}} \right)$. To secure the unbiasedness, we introduce a positive (time-varying) local variable ${\theta _{i,k}} = \theta \left( {{S_{i,k}}} \right)$, such that
\begin{equation}\label{E4-2-4}
{\mathbb{E}_{\mathcal{D}_i}}\left[ {{\theta _{i,k}}{M_{i,k}}} \right] = {W_i},
\end{equation}
where ${\mathbb{E}}_{\mathcal{D}_i}\left[ \cdot \right]$ implies the expectation taken over the randomness in $\mathcal{D}_i$. Let $\mathcal{D}: = \prod\nolimits_{i \in \mathcal{R}} {{\mathcal{D}_i}} $, where ${\prod _{\cdot}}$ implies the Cartesian product.
Based on this relationship, we define the stochastic gradient by ${g_{i,k}} = \left( {1 - {\theta _{i,k}}} \right){h _{i,k}} + {\theta _{i,k}}{h _{i,k + 1}}$, such that applying (\ref{E4-2-2}) yields
\begin{equation}\label{E4-2-5}
{g_{i,k}} = {h _{i,k}} + {\theta _{i,k}}W_i^{ - 1}{M_{i,k}}\left( {\nabla {f_i}\left( {{x_{i,k}}} \right) - {h _{i,k}}} \right),
\end{equation}
which produces an unbiased estimator of $\nabla {f_i}\left( {{x_{i,k}}} \right)$, i.e.,
\begin{equation}\label{E4-2-6}
{\mathbb{E}}_{\mathcal{D}_i}\left[ {{g_{i,k}}} \right] = \nabla {f_i}\left( {{x_{i,k}}} \right).
\end{equation}

\begin{assumption}(Mutual independence)\label{A4}
The local distributions $ {{\mathcal{D}_i}}$, ${i \in \mathcal{R}}$, are mutually independent.
\end{assumption}

We exemplify localized gradient sketching in the sequel.
\begin{example}(Multi-coordinate gradient sketching)\label{Exa1}
For any agent $i \in \mathcal{R}$, let the local weighted matrix ${W_i} = {\mathbf{Diag}}\left\{ {w_i^1, w_i^2, \ldots, w_i^n} \right\} \succ 0$. Suppose the local distribution $\mathcal{D}_i$ corresponds to selecting a random subset of coordinate indices ${\mathcal{J}_{i,k}} \subseteq \left\{ {1, 2, \ldots, n} \right\}$ of size $b_i$ uniformly without replacement. Thus, the marginal probability of any coordinate $j$ being selected is $b_i/n$. We define the sketched matrix as ${S_{i,k}} = \left[ {e_i^{j_1}, e_i^{j_2}, \ldots, e_i^{j_{b_i}}} \right]$ for $j \in \mathcal{J}_{i,k}$, where ${\left\{ {e_i^j} \right\}_{j = 1}^n}$ are the standard unit basis vectors in $\mathbb{R}^n$. According to (\ref{E4-2-2}), the gradient estimate is updated by:
\begin{equation}\label{E4-2-8}
{h_{i,k + 1}} = {h_{i,k}} + \sum\limits_{j \in {\mathcal{J}_{i,k}}} {e_i^j{{\left( {e_i^j} \right)}^\top}\left( {\nabla {f_i}\left( {{x_{i,k}}} \right) - {h_{i,k}}} \right)},
\end{equation}
which is equivalent to ${\left[ {{h_{i,k + 1}}} \right]_j} = {\left( {e_i^j} \right)^\top}\nabla {f_i}\left( {{x_{i,k}}} \right)$ for $j \in {\mathcal{J}_{i,k}}$, and ${\left[ {{h_{i,k + 1}}} \right]_j} = {\left[ {{h_{i,k}}} \right]_j}$ otherwise. From (\ref{E4-2-8}), we observe that ${h_{i,k + 1}}$ is independent of the local weighted matrix ${W_i}$.

Recalling the definition of $M_{i,k}$, we have $M_{i,k} = \sum_{j \in \mathcal{J}_{i,k}} w_i^j e_i^j (e_i^j)^\top$. By setting the scaling parameter to ${\theta_{i,k}} = n/b_i$, the expectation yields:
\begin{equation}\label{E4-2-9}
\begin{aligned}
\mathbb{E}_{\mathcal{D}_i}\left[ {{\theta_{i,k}}{M_{i,k}}} \right] &= \frac{n}{{{b_i}}}\sum\limits_{j = 1}^n {\mathbb{P}\left( {j \in {\mathcal{J}_{i,k}}} \right)w_i^je_i^j{{\left( {e_i^j} \right)}^\top}} \\
&= \frac{n}{{{b_i}}}\sum\limits_{j = 1}^n {\left( {\frac{{{b_i}}}{n}} \right)w_i^je_i^j{{\left( {e_i^j} \right)}^\top}} \\
&= {W_i},
\end{aligned}
\end{equation}
which satisfies the unbiasedness condition and generates the stochastic gradient estimator:
\begin{equation}\label{E4-2-10}
{g_{i,k}} = {h_{i,k}} + \frac{n}{{{b_i}}}\sum\limits_{j \in {\mathcal{J}_{i,k}}} {e_i^j{{\left( {e_i^j} \right)}^\top}\left( {\nabla {f_i}\left( {{x_{i,k}}} \right) - {h_{i,k}}} \right)}.
\end{equation}
\end{example}

\subsection{Decentralized Stochastic (Sub)Gradient Descent}\label{sec4-3}
To resolve \text{DCC} (\ref{E3-3-1}), each reliable agent $i$, $i \in \mathcal{R}$, performs a decentralized SGD step as follows:
\begin{equation}\label{E4-3-0}
\begin{aligned}
{{\bar x}_{i,k}} = & {\left[ {{A_k}} \right]_{ii}}{x_{i,k}} + \sum\limits_{j \in {\mathcal{N}_{i,k}}} {{{\left[ {{A_k}} \right]}_{ij}}{x_{j,k}}}  - \alpha {g_{i,k}}\\
= & \sum\limits_{j \in {\mathcal{R}_{i,k}} \cup \left\{ i \right\}} {{{\left[ {{A_k}} \right]}_{ij}}{x_{j,k}}}   - \alpha {g_{i,k}} + \underbrace {\sum\limits_{j \in {\mathcal{B}_{i,k}}} {{{\left[ {{A_k}} \right]_{ij}}}{x_{j,k}^i}} }_{\text{Distraction}},\\
\end{aligned}
\end{equation}
where $A_k$ denotes a row-stochastic time-varying weight matrix, ${\left[ {{A_k}} \right]_{ij}}$ denotes the $(i,j)$-th entry of $A_k$, and ${x_{j,k}^i} := \left\{ \begin{gathered}
  {x_{j,k}}, {\text{if}} \; j \in {\mathcal{R}_{i,k}} \hfill \\
  {z_{j,k}^i}, {\text{if}} \; j \in {\mathcal{B}_{i,k}} \hfill \\
\end{gathered}  \right.$ with ${z_{j,k}^i}$ being an arbitrary vector in $ {\mathbb{R}^n}$. Therein, we use $x_{j,k}^i$ to imply that Byzantine agent $j$ can send different values to its different reliable neighbors $i$, $i \in {\mathcal{N}_{j,k}}$. The distraction term of reliable agent $i$ is induced by its indistinguishable Byzantine neighbors, which could deviate its state arbitrarily. For instance, when Byzantine agent $j$ sends ${x_{j,k}^i} = {x_{i,k}} - \sum\nolimits_{j \in {\mathcal{R}_{i,k}}} {{{\left[ {{A_k}} \right]}_{ij}}\left( {{x_{j,k}} - {x_{i,k}}} \right)} /\left( {\sum\nolimits_{j \in {\mathcal{B}_{i,k}}} {{{\left[ {{A_k}} \right]}_{ij}}} } \right)$ to its reliable neighbor $i$, then the update (\ref{E4-3-0}) becomes an isolated gradient descent step without any information exchange, i.e., ${{\bar x}_{i,k}} = {x_{i,k}} - \alpha {g_{i,k}}$, which is known as the dissensus attack \cite{He2022,Hu2025}.

To mitigate the influence of Byzantine agents, each reliable agent $i$, $i \in \mathcal{R}$, performs a resilient decentralized preconditioned stochastic subgradient descent step according to \text{DRC} (\ref{E3-3-2}) as follows:
\begin{equation}\label{E4-3-1}
{\bar x_{i,k}} = {x_{i,k}} - \alpha W_i^{ - 1}\left( {{g_{i,k}} + \phi \sum\limits_{j \in {\mathcal{N}_{i,k}}} {\partial {{\left\| {{x_{i,k}} - {x_{j,k}^i}} \right\|}_p}} } \right),
\end{equation}
It should be noted that the identity of each agent $j$, $j \in {\mathcal{N}_{i,k}}$ (whether it is reliable or Byzantine) is not necessarily known to its neighboring reliable agent $i$, $i \in \mathcal{R} $. It is also worth noting that the instantaneous resilient aggregation in the \textit{RED-SEGA} update (\ref{E4-3-1}) intrinsically relates to the deterministic frequency-weighted objective (\ref{E3-3-2}) via stochastic subgradient estimation. For any reliable agent $i$, the local penalty subgradient evaluated at iteration $k$ is driven by the time-varying network topology $\mathcal{E}_k$. If we isolate the information from reliable neighbors, taking the expectation over the network randomness $\zeta_k$ yields
\begin{equation}\label{E4-3-2}
\mathbb{E}_{\zeta_k}\left[ \sum_{j\in\mathcal{R}_{i,k}} \partial\vert{}\vert{}x_{i,k} - x_{j,k}\vert{}\vert{}_p\right] =  \sum_{j\in\bar{\mathcal{R}}_i} \omega _{\left( {i,j} \right)} \partial\vert{}\vert{}x_{i,k} - x_{j,k}\vert{}\vert{}_p.
\end{equation}
This reveals that the instantaneous consensus step over the time-varying reliable network is exactly an unbiased stochastic subgradient of the frequency-weighted penalty term in the deterministic DRC problem (\ref{E3-3-2}). The malicious messages injected by Byzantine agents $j \in \mathcal{B}_{i,k}$ in (\ref{E4-3-1}) can thus be mathematically interpreted as bounded adversarial perturbations added to this unbiased stochastic subgradient, which are subsequently controlled by the algorithm's resilience mechanisms.

\subsection{Decentralized Preconditioned Proximal Mapping}\label{sec4-4}
We then define a $W_i$-weighted proximal operator to minimize the local nonsmooth objective function $r_i$ as follows:
\begin{equation}\label{E4-3-2}
{\mathbf{prox}}_{\alpha ,{r_i}}^{{W_i}}\left\{ {\tilde x} \right\}: = \arg \mathop {\min }\limits_{\tilde y} \left\{ {{r_i}\left( {\tilde y} \right) + \frac{1}{{2\alpha }}\left\| {\tilde y - \tilde x} \right\|_{{W_i}}^2} \right\},
\end{equation}
where ${W_i} \in {\mathbb{R}^{n \times n}}$ is an arbitrary but positive-definite matrix, such that the next update $x_{i, k+1}$ takes the form of
\begin{equation}\label{E4-3-3}
x_{i,k+1} = {\mathbf{prox}}_{\alpha ,{r_i}}^{{W_i}}\left\{ {\bar x_{i,k}} \right\},
\end{equation}
which is a preconditioned proximal mapping step. Note that the local nonsmooth objective functions correspond to either nonsmooth penalty, for instance the Manhattan norm \cite{Hu2025a,Ye2025}, or the indicator function with respect to convex set constraints \cite{Li2021,Li2021b}.

\subsection{Development of \textit{Gossip-SEGA} and \textit{RED-SEGA}}\label{sec4-5}
To address \text{DRC} (\ref{E3-3-2}) in a decentralized manner, we summarize updates (\ref{E4-2-2}), (\ref{E4-2-5}), (\ref{E4-3-0}), and (\ref{E4-3-3}) to develop a non-resilient decentralized stochastic proximal gradient algorithm presented in Algorithm \ref{Algo1}. To further enable a resilient aggregation of Algorithm \ref{Algo1}, we replace the naive \textit{Gossip} aggregation in (\ref{E4-3-0}) with a resilient aggregation in (\ref{E4-3-1}) to devise a Byzantine-resilient decentralized preconditioned stochastic proximal gradient algorithm with variance reduction, which is elaborated in Algorithm \ref{Algo2}.
\begin{algorithm}[!htp]
	\small
    \SetKwInput{KwInit}{Initialize}
	\SetKwBlock{Repeat}{Repeat for $ k = 0,1,2, \dots$}{}
	\SetKwBlock{End}{End for a required criterion is satisfied}{}
	\KwIn{a proper constant $\alpha>0$.}
	\KwInit{arbitrary starting points $x_{i,0} \in {\mathbb{R}^n}$, $\forall i \in \mathcal{V}$, and a proper penalty parameter ${\phi}>0$.}
	\For{$k=0,1,\ldots,$}{
		\For{\rm{each reliable agent} $i \in \mathcal{R}$}{
			{\textbf{Transmit} its current local model ${x_{i,k}}$ to its neighbors $j \in {\mathcal{N}_i}$ and receive the true information ${x_{j,k}}$ or untrue information ${z_{j,k}^i}$ from its neighbors;}\\
			{\textbf{Sample} a sketched operator $S_{i,k}$ locally according to a fixed distribution ${\mathcal{D}_i}$, i.e., ${S_{i,k}} \sim {\mathcal{D}_i}$;}\\
            {\textbf{Reconstruct} the local true gradient incrementally according to the sketch-and-project step (\ref{E4-2-2});}\\
			{\textbf{Evaluate} the local stochastic gradient $g_{i,k}$ according to (\ref{E4-2-5});}\\
			{\textbf{Estimate} an intermediate variable according to the decentralized stochastic gradient descent step (\ref{E4-3-0});}\\
			{\textbf{Update} the local model according to the preconditioned proximal gradient descent step (\ref{E4-3-3}).}
		}
		\For{\rm{each Byzantine agent} $i \in \mathcal{B}$}{\textbf{Send} an arbitrary vector $ {z_{j,k}^i} \in {\mathbb{R}^n}$ to its neighbor $i$, $i \in {\mathcal{N}_j}$.}
	}
	\KwOut{the local decision variable $x_{i,K}$ of each reliable agent $i \in \mathcal{R}$ until a prescribed criterion is satisfied.}
	\caption{\textit{Gossip-SEGA}.}
	\label{Algo1}
\end{algorithm}
\begin{algorithm}[!htp]
	\small
    \SetKwInput{KwInit}{Initialize}
	\SetKwBlock{Repeat}{Repeat for $ k = 0,1,2, \dots$}{}
	\SetKwBlock{End}{End for a required criterion is satisfied}{}
	\KwIn{a proper constant $\alpha>0$.}
	\KwInit{arbitrary starting points $x_{i,0} \in {\mathbb{R}^n}$, $\forall i \in \mathcal{V}$, and a proper penalty parameter ${\phi}>0$.}
	\For{$k=0,1,\ldots,$}{
		\For{\rm{each reliable agent} $i \in \mathcal{R}$}{
            {Perform the same steps 3-8 as in Algorithm \ref{Algo1}, except that (\ref{E4-3-0}) in Step 7 is replaced by (\ref{E4-3-1}).}
		}
		\For{\rm{each Byzantine agent}}{ Perform the same as in Algorithm \ref{Algo1}.}
	}
	\KwOut{the local decision variable $x_{i,K}$ of each reliable agent $i \in \mathcal{R}$ until a prescribed criterion is satisfied.}
	\caption{\textit{RED-SEGA}.}
	\label{Algo2}
\end{algorithm}
\begin{remark}\label{R4}
To solve SRM problems over time-varying networks in the presence of Byzantine agents, \textit{RED-SEGA} integrates gradient sketching with a norm-penalized approximation. This represents a different mechanism for achieving variance reduction compared to the baseline \textit{Prox-DBRO-VR} \cite{Hu2025a}. By evaluating only a sketched subset of $b$ coordinates across the local dataset, the gradient evaluation complexity is reduced from $\mathcal{O}(n D_{\max})$ to $\mathcal{O}(b D_{\max})$, where $D_{\max} := \max_i D_i$ denotes the maximum local data size $D_i$. This implies that \textit{RED-SEGA} not only demonstrates strong adaptability to time-varying networks but also exhibits a significant computational edge over \textit{Prox-DBRO-VR} in data-scarce yet high-dimensional real-world applications, i.e., $b D_{\max} \ll n$, for instance high-resolution image deblurring. As revealed by the experimental results in Section \ref{sec6-1}, \textit{RED-SEGA} achieves consensus and convergence performance comparable to \textit{Prox-DBRO-VR} \cite{Hu2025a} in terms of both iterations and oracle calls (scalar derivative evaluations) under fair comparative settings.
\end{remark}

\section{Theoretical Analysis}\label{sec5}
This paper aims to provide a theoretical analysis for \textit{RED-SEGA} while the convergence performance of \textit{Gossip-SEGA} is empirically studied. To this end, we first divide the theoretical analysis of \textit{RED-SEGA} into two parts: consensus and convergence analysis.
\subsection{Consensus Analysis}\label{sec5-1}
We seek a condition on the penalty parameter $\phi$ to ensure the equivalence between \text{SRM} (\ref{E3-1-1}) and \text{DRC} (\ref{E3-3-2}) in the sequel theorem. We first define a weighted node-edge incidence matrix $\Omega$ associated with the \textit{average} network $\bar {\mathcal{G}}_{\mathcal{R}}$. The $\left( {i,e} \right)$th entry of $\Omega$ is $\omega_e $ while the $\left( {j,e} \right)$th entry of $\Omega$ is $-\omega_e $. The minimum non-zero singular value of $\Omega$ is denoted as ${{\sigma _{\min }}\left( \Omega  \right)}$. Different to \cite{Hu2025a}, we provide a more rigorous analysis on the resilient consensus condition in Theorem \ref{T1}.
\begin{theorem}(Resilient consensus condition)\label{T1}
Suppose that Assumptions \ref{A1} and \ref{A3} hold. Let $\Omega$ be the weighted node-edge incidence matrix associated with the \textit{average} network $\bar {\mathcal{G}}_{\mathcal{R}}$, and let $C_{d,2}$ be a positive constant characterizing the equivalence between the $\ell_d$-norm and the $\ell_2$-norm. There exists a set of subgradients $\{ \partial r_i^* \in \partial r_i(\tilde{x}^*) \}_{i \in \mathcal{R}}$ satisfying the global optimality condition $\sum_{i \in \mathcal{R}} ( \nabla f_i(\tilde{x}^*) + \partial r_i^* ) = 0_n$, such that if the penalty parameter $\phi$ satisfies
\begin{equation}\label{E5-1-1}
\phi  \ge {\phi _{\min }}: = \frac{{{C_{d,2}}\sqrt {|\mathcal{R}|} }}{{{\sigma _{\min }}\left( \Omega  \right)}}{\max _{i \in \mathcal{R}}}{\left\| {\nabla {f_i}\left( {{{\tilde x}^*}} \right) + \partial r_i^*} \right\|_2},
\end{equation}
where $C_{d,2} = {n^{\max \left\{ {0,\frac{1}{d} - \frac{1}{2}} \right\}}}$, then the optimal solution to \text{DRC} (\ref{E3-3-2}) is strictly equivalent to that of \text{DCC} (\ref{E3-3-1}) and \text{SRM} (\ref{E3-1-1}), i.e., ${x^*} = {1_{\left| \mathcal{R} \right|}} \otimes {{\tilde x}^*}$.
\begin{proof}
The detailed proof has been deferred to Section A of the Supplementary Material due to space constraints.
\end{proof}
\end{theorem}

\begin{remark}\label{R5}
To ensure consensus among all reliable agents, \textit{RED-SEGA} requires a lower bound on the penalty parameter $\phi$, as proven in Theorem \ref{T1}. Notably, this threshold is fundamentally independent of both the states and the number of Byzantine agents. In practice, the tolerable number of Byzantine agents is constrained only by the topological necessity of ensuring that reliable agents avoid isolation to maintain sufficient information aggregation under Assumption \ref{A3}. However, while Theorem \ref{T1} establishes that a sufficiently large penalty parameter $\phi$ can guarantee exact equivalence between the original problem \text{SRM} (\ref{E3-1-1}) and its resilient decentralized surrogate \text{DRC} (\ref{E3-3-2}), the subsequent analysis will reveal an inherent trade-off: a larger $\phi$ inadvertently inflates the convergence error bound. Consequently, the theoretical mandate for a sufficiently large penalty acts primarily as an analytical baseline. Therefore, $\phi$ should be treated as a tunable hyperparameter for practical deployments, optimized empirically to strike a delicate balance between Byzantine-resilient consensus and minimal steady-state error.
\end{remark}

\subsection{Convergence Analysis}\label{sec5-2}
To facilitate the subsequent analysis, we define the following vectors and matrices: $\Delta _{i,k}^\psi : = {h _{i,k}} - \nabla {f_i}\left( {{{\tilde x}^*}} \right)$, $\Delta _{i,k}^f: = \nabla {f_i}\left( {{x_{i,k}}} \right) - \nabla {f_i}\left( {{{\tilde x}^*}} \right)$, ${\chi _{i,k}}: = \phi \sum\nolimits_{j \in {\mathcal{R}_{i,k}}} {{{\left\| {{x_{i,k}} - {x_{j,k}}} \right\|}_p}} $, and ${\delta _{i,k}}: = \phi \sum\nolimits_{j \in {\mathcal{B}_{i,k}}} {{{\left\| {{x_{i,k}} - {z_{j,k}^i}} \right\|}_p}} $. Given that the random sketched matrices $S_{i,k}$ are identically and independently sampled from a fixed local distribution $\mathcal{D}_i$ at each iteration $k \ge 0$, both $\theta_{i,k}$ and $M_{i,k}$ exhibit stationary statistical properties over time. By virtue of this statistical stationarity, their expectations yield the time-independent constant matrices: ${M_i}: = {\mathbb{E}_{{\mathcal{D}_i}}}\left[ {{M_{i,k}}} \right] = {\mathbb{E}_{{\mathcal{D}_i}}}\left[ {{M_{i,k + 1}}} \right]$ and ${\Phi _i}: = {\mathbb{E}_{{\mathcal{D}_i}}}\left[ {\theta _{i,k}^2{M_{i,k}}} \right] = {\mathbb{E}_{{\mathcal{D}_i}}}\left[ {\theta _{i,k + 1}^2{M_{i,k + 1}}} \right]$. Let $\mathcal{F}_k$ denote the filtration generated by all sources of randomness including both the gradient sketching distributions and the sequence of network topologies up to iteration $k$. For brevity, we define ${\mathbb{E}_k}\left[ {\cdot} \right]: = \mathbb{E}\left[ {\cdot |{\mathcal{F}_k}} \right]$ as the conditional expectation given the historical information up to step $k$ as appropriate. Under the conditional expectation $\mathbb{E}_k[\cdot]$, the historical variables such as $x_{i,k}$ and $h_{i,k}$ become deterministic constants, while the expected values of the instantaneous random variables, for instance, $M_{i,k}$ and $\zeta_k$, can be evaluated directly based on their independent distributions. The total (unconditional) expectation taken over all randomness from step $0$ to $k$ is denoted simply by $\mathbb{E}[\cdot]$, which naturally satisfies the law of total expectation $\mathbb{E}\left[ {\cdot} \right] = \mathbb{E}\left[ {{\mathbb{E}_k}\left[ {\cdot} \right]} \right]$.
\begin{lemma}(Gradient estimation error)\label{L1}
Suppose that Assumption \ref{A4} holds. For all $i \in \mathcal{R}$ and $k \ge 0$, we have
\begin{equation}\label{E5-2-1}
{\mathbb{E}_k}\left[\left\| {\Delta _{i,k + 1}^\psi } \right\|_{{W_i}}^2 \right] = {\left\| {\Delta _{i,k}^\psi } \right\|_{{W_i} - {M_i}}^2 + \left\| {\Delta _{i,k}^f} \right\|_{{M_i}}^2}.
\end{equation}
\begin{proof}
The detailed proof has been deferred to Section B of the Supplementary Material due to space constraints.
\end{proof}
\end{lemma}

\begin{lemma}(Gradient-learning error)\label{L2}
Suppose that Assumption \ref{A4} holds. For all $i \in \mathcal{R}$ and $k \ge 0$, we have
\begin{equation}\label{E5-2-2}
\begin{aligned}
{\mathbb{E}_k}\left[{\left\| {{g_{i,k}} - \nabla f_i\left( {{{\tilde x}^*}} \right)} \right\|_{{W_i}}^2}\right]  \le 2 {\left\| {\Delta _{i,k}^\psi } \right\|_{{\Phi _i} - {W_i}}^2 + \left\| {\Delta _{i,k}^f} \right\|_{{\Phi _i}}^2},
\end{aligned}
\end{equation}
where ${\Phi _i} - {W_i} \succcurlyeq 0$.
\begin{proof}
The detailed proof has been deferred to Section C of the Supplementary Material due to space constraints.
\end{proof}
\end{lemma}

Before presenting the main theorem, we define $\mu : = {\min _{i \in \mathcal{R}}}{\mu _i}$, ${x_k}: = {\mathbf{col}}{\left\{ {{x_{i,k}}} \right\}_{i \in \mathcal{R}}}$, $W: = {\mathbf{diag}}{\left\{ {{W_i}} \right\}_{i \in \mathcal{R}}}$, $\Delta _k^\psi : = {\mathbf{col}}{\left\{ {\Delta _{i,k}^\psi } \right\}_{i \in \mathcal{R}}}$, and ${U_k^W}: = \left\| {{x_k} - {x^*}} \right\|_W^2 + c\alpha \left\| {\Delta _k^\psi} \right\|_W^2$ with $c$ being a positive constant. For brevity, we denote $\Delta _k^F : = {\mathbf{col}}{\left\{ {\Delta _{i,k}^f } \right\}_{i \in \mathcal{R}}}$, ${\Gamma _1}: = 4n{\phi ^2}\sum\nolimits_{i \in \mathcal{R}} {\left( {4{{\left| {{\mathcal{R}_i}} \right|}^2} + {{\left| {{\mathcal{B}_i}} \right|}^2}} \right)/{\lambda _{\min }}\left( {{W_i}} \right)}  $, ${\Gamma _2}: = 16n{\phi ^2}\sum\nolimits_{i \in \mathcal{R}} {\left( {4{{\left| {{\mathcal{R}_i}} \right|}^2} + {{\left| {{\mathcal{B}_i}} \right|}^2}} \right)} {\lambda _{\max }}\left( {{W_i}} \right)/\left( {\mu {\lambda _{\min }}\left( {{W_i}} \right)} \right)$, $Q: = {\mathbf{diag}}{\left\{ {{Q_i}} \right\}_{i \in \mathcal{R}}}$, $M := \mathbf{diag}\{M_i\}_{i \in \mathcal{R}} $, $\Phi := \mathbf{diag}\{\Phi_i\}_{i \in \mathcal{R}}$, $\tilde{\Phi} = \text{diag}\{\lambda_{\min}(W_i)^{-2}\Phi_i\}_{i\in\mathcal{R}}$, and $\tilde{W} := \text{diag}\{\lambda_{\min}(W_i)^{-2}W_i\}_{i\in\mathcal{R}}$, and $\gamma   := \mu /2$.
\begin{theorem}(Linear convergence)\label{T2}
Suppose that Assumptions \ref{A1}-\ref{A4} hold. Under the conditions of Theorem \ref{T1} and a constant $c$ satisfying $0 < c < {{\lambda _{\min }}\left( Q \right)}/\left( {2{\lambda _{\max }}\left( {M} \right)} \right)$, if a constant step-size is employed and satisfies
$0 < \alpha  \le \min \left\{ {\frac{{c{\lambda _{\min }}\left( {M} \right)}}{{{\lambda _{\max }}\left( {\gamma  cW + 4\left( {\tilde \Phi  - \tilde W} \right)} \right)}},\frac{{{\lambda _{\min }}\left( {Q - cM} \right)}}{{4{\lambda _{\max }}\left(\tilde \Phi  \right)}}} \right\}$, then the sequence ${\left\{ {{x_k}} \right\}_{k \ge 0}}$ generated by \textit{RED-SEGA} takes the following convergence form
\begin{equation}\label{E5-2-3}
\mathbb{E}\left[ {{U_k^W}} \right] \le {\left( {1 - \gamma  \alpha } \right)^{k}}{U_0^W} + \frac{{ {{\Gamma  _1}\alpha  + {\Gamma  _2}} }}{\gamma }\left( {1 - {{\left( {1 - \gamma \alpha } \right)}^{k }}} \right),
\end{equation}
where ${U_0^W} = \left\| {{x_0} - {x^*}} \right\|_W^2 + c\alpha \left\| {\Delta _0^\psi} \right\|_W^2$. This implies that the $\{x_k\}_{k \ge 0}$ generated by \textit{RED-SEGA} converges linearly to an error ball centered at the global optimal solution $x^*$ solution $x^*$ with a rate of $(1 - \mathcal{O}(\gamma \alpha))^k$, where the radius of this error ball is proportional to $(\Gamma_1 \alpha + \Gamma_2) / \gamma$.
\begin{proof}
The detailed proof has been deferred to Section D of the Supplementary Material due to space constraints.
\end{proof}
\end{theorem}

\begin{remark}\label{R8}
To provide concrete guidance for practitioners, we illustrate the abstract step-size condition in Theorem \ref{T2} based on the multi-coordinate gradient sketching scenario from Example \ref{Exa1}. Let us consider a Euclidean metric where $W_i = {\mathbf{I}}$ for all $i \in \mathcal{R}$. Suppose the local objectives are $L$-smooth and $\mu$-strongly convex, yielding $Q = L^{-1} {\mathbf{I}}$. Furthermore, assume each agent uniformly samples $b$ out of $n$ coordinates at each iteration without replacement. As derived in Example \ref{Exa1}, under this uniform sampling regime, the time-independent sketching matrix strictly reduces to $M = \left( {b/n} \right){\mathbf{I}}$. To satisfy the unbiasedness requirement, the scaling parameter is determined as $\theta_{i,k} = n/b$, which consequently yields the time-independent variance scaling matrix $\Phi = \left( {n/b} \right){\mathbf{I}}$. Substituting these specific scalar matrices into the analytical bounds established in Theorem 2, the abstract matrix-based step-size condition elegantly simplifies into the following tractable scalar constraint $0 < \alpha \le \min\left\{ \frac{c \left( {n/b} \right)}{ c \gamma + 4\left(\left( {n/b} \right) - 1\right)}, \frac{L^{-1} - c\left( {b/n} \right)}{4\left( {n/b} \right)} \right\}$. This explicit formulation mathematically reveals the direct impact of the sketching ratio $b/n$ on the permissible step-size. Notably, when $b=n$, the variance-induced penalty term $4(n/b - 1)$ strictly vanishes, allowing the proposed algorithm to recover the standard step-size bounds.
\end{remark}

\begin{remark}\label{R9}
A key advantage of \textit{RED-SEGA} lies in its scalability for high-dimensional problems. For each agent, the gradient sketching step, e.g., multi-coordinate sketching with $b \ge 1$, introduces an additional projection step. However, as demonstrated in Example \ref{Exa1}, this structured sketching algebraically reduces the projection matrix into a sparse masking operation. Consequently, its computational overhead is negligible, requiring only $\mathcal{O}(b D_{\max})$ operations for local gradient evaluation. Furthermore, unlike traditional filtration- or screening-based resilient aggregation methods, such as coordinate-wise median \cite{Yin2018} and trimmed mean \cite{Xie2018}, which incur a heavy local computational overhead of $\mathcal{O}(n |\mathcal{N}_{\max}| \log |\mathcal{N}_{\max}|)$ (where $|\mathcal{N}_{\max}| := \max_i |\mathcal{N}_i|$) due to sorting operations across all dimensions, \textit{RED-SEGA} utilizes a norm-penalized approximation for resilient consensus. Evaluating the subgradient of this penalty involves only element-wise vector operations, reducing the aggregation overhead per agent to $\mathcal{O}(n |\mathcal{N}_{\max}|)$. Consequently, the overall per-iteration computational complexity for any agent is bounded by $\mathcal{O}(n |\mathcal{N}_{\max}| + b D_{\max})$. This linear complexity mathematically guarantees both efficiency and robust scalability in large-scale MASs, as the per-agent computational burden scales linearly with respect to both the problem dimension $n$ and the maximum local data size $D_{\max}$.
\end{remark}

To demonstrate that the convergence error of \textit{RED-SEGA} vanishes in the absence of Byzantine agents, the following corollary establishes its sublinear convergence under a decaying step-size when $W = \mathbf{I}$. Prior to presenting the corollary, we define ${U_k}: = {\left\| {{x_k} - {x^*}} \right\|^2} +{\tilde c}{\alpha _k}{\left\| {\Delta _k^\psi } \right\|^2}$, where $\tilde c$ is a positive constant satisfying $0 < \tilde c < {\lambda _{\min }}\left( Q \right)/\left( {2{\lambda _{\max }}\left( M \right)} \right)$, ${{\Gamma _3}} := {16n{\phi ^2}\sum\nolimits_{i \in \mathcal{R}} {\left( {{{\left| {{\mathcal{R}_i}} \right|}^2} + {{\left| {{\mathcal{B}_i}} \right|}^2}} \right)} }$, ${{\Gamma _4}}:= {n\gamma {\phi ^2}\sum\nolimits_{i \in \mathcal{R}} {{{\left| {{\mathcal{B}_i}} \right|}^2}} }$, $\varpi = \min \left\{ {\frac{{{\lambda _{\min }}\left( {Q - \tilde c  M} \right)}}{{4{\lambda _{\max }}\left( \Phi  \right)}},\frac{{\tilde c{\lambda _{\min }}\left( M \right)}}{{{\lambda _{\max }}\left( {\gamma \tilde c + {\text{4}}\left( {\Phi  - {\mathbf{I}}} \right)} \right)}}} \right\}$, $\xi  > \beta  \max \left\{ {\gamma ,\varpi } \right\} $ with $\beta   > 1/\gamma $,  and $\Xi : = \max \left\{ {\frac{{{\beta  ^2}{\Gamma _3}}}{{\gamma \beta   - 1}},\left( {\xi  - \gamma \beta } \right){U_0} + \frac{{{\beta  ^2}}}{\xi }{\Gamma _3} + \beta  {\Gamma _4} - \frac{{{\xi \Gamma _4}}}{\gamma } } \right\}$.
\begin{corollary}(Sublinear convergence)\label{C1}
Suppose that Assumptions \ref{A1}-\ref{A4} hold. Under the conditions of Theorem \ref{T1} and for any constant $c \ge 0$, if a decaying step-size ${\alpha _k} = \beta  /\left( {\xi  + k} \right)$, $\forall k \ge 0$, is employed and the weight matrix satisfies $W = \mathbf{I}$, then the sequence ${\left\{ {{x_k}} \right\}_{k \ge 0}}$ generated by \textit{RED-SEGA} takes a convergence form of
\begin{equation}\label{E5-2-4}
\mathbb{E}\left[ U_k \right] \le \frac{\Xi }{{\xi  + k}} + \frac{{{\Gamma _4}}}{\gamma },
\end{equation}
where the error term $\Gamma _4$ vanishes in the absence of Byzantine agents, i.e., $\left| {{\mathcal{B}_i}} \right| = 0$. This implies that the sequence $\{x_k\}_{k \ge 0}$ generated by \textit{RED-SEGA} converges sublinearly to an error ball centered at the global optimal solution $x^*$ with a rate of $\mathcal{O}\left( {1/k } \right)$, where the radius of this error ball is proportional to $\Gamma _4 / \gamma$.
\begin{proof}
The detailed proof has been deferred to Section E of the Supplementary Material due to space constraints.
\end{proof}
\end{corollary}
\begin{remark}\label{R10}
Given that Theorem \ref{T2} establishes a non-asymptotic linear convergence rate to a bounded steady-state error ball with a constant step-size, a strictly vanishing-error result, i.e., exact convergence, is unachievable due to the fundamental impossibility of accurately detecting and isolating Byzantine agents without any \textit{prior} knowledge, thereby rendering their disruptive impact on the convergence of reliable agents a persistent factor \cite{Wu2022,He2022,Han2025,Liu2023c,Hu2025a,Hu2025}. Nevertheless, Corollary \ref{C1} provides a theoretical resolution to the fundamental tension between the consensus penalty parameter $\phi$ and the convergence accuracy (characterized by the radius of the convergence error ball). Specifically, while Theorem \ref{T1} indicates that a sufficiently large $\phi $ is necessitated to resist Byzantine attacks, which, however, inevitably enlarges the convergence error ball as verified by Theorem \ref{T2}, Corollary \ref{C1} demonstrates that employing a decaying step-size asymptotically reduces the convergence error and the radius of the error ball is proportional to the number of Byzantine agents since $\Gamma_4 \propto \sum\nolimits_{i \in \mathcal{R}} |\mathcal{B}_i|^2$. This implies that the error term $\Gamma_4 /\gamma$ in (\ref{E5-2-4}) vanishes in the absence of Byzantine agents, i.e., $\left| \mathcal{B} \right| = 0$, which enables \textit{RED-SEGA} to achieve exact convergence.
\end{remark}

\section{Numerical Experiments}\label{sec6}
This section provides two case studies to illustrate the effectiveness and practicality of \textit{RED-SEGA} by solving a constrained least-squares problem and an image deblurring task, respectively.
\subsection{Decentralized Least Square Under Byzantine Attacks}\label{sec6-1}
\begin{figure}[!htp]
  \centering
  \includegraphics[scale=0.22]{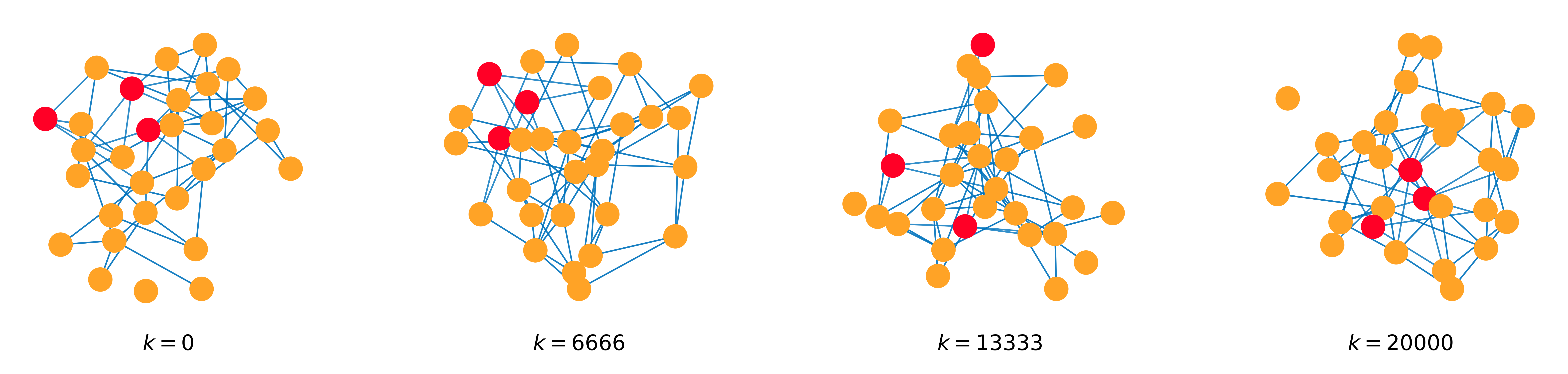}
  \caption{Example of dynamic evolution of Byzantine time-varying networks.}\label{Fig0-3}
\end{figure}

In this section, we compare the performance of \textit{Gossip-SEGA} and \textit{RED-SEGA} under various Byzantine attacks, including 1) dropout attacks: Byzantine agents occasionally drops out of communication with their reliable neighbors and the dropout probability is randomly selected from $\left[ {0.5,1.0} \right]$ over iterations; 2) Gaussian attacks: Byzantine agents inject Gaussian noise to the true values (see \cite[Section VII]{Wu2022} for specific settings of the injected noise), introducing random perturbations that aim to degrade the quality of model aggregation; 3) A-Little-Is-Enough attacks: Byzantine attacks generate malicious values that are slightly below the statistical mean of reliable nodes' values (see \cite{Baruch2019} for specific settings of the malicious values), rendering the attack subtle and hard to detect while still being effective in disrupting convergence. We consider a multi-agent system, including both reliable and Byzantine agents, to minimize a constrained least-square problem as follows:
\begin{equation}\label{E6-1}
\mathop {\min }\limits_{\left\| \tilde{x} \right\| \le 1} \sum\limits_{i \in \mathcal{R}} {\left\| {{\mathcal{A}_i}{\tilde x} - {\mathcal{B}_i}} \right\|^2} ,
\end{equation}
where the entries of ${\mathcal{A}_i} \in {\mathbb{R}^{m \times n}}$ are i.i.d. synthetic samples from a standard normal distribution and ${{\mathcal{B}}_i} \in {\mathbb{R}^m}$ is the all-ones vector. The number of data samples is $m$ and the variable dimension is set to $n=100*m$ (i.e., $m \ll n$).  Similar to \cite{Li2021b,Li2021}, the inequality constraint in (\ref{E6-1}) is tackled by introducing a nonsmooth indicator function. We assume that there is an oracle that can return the partial derivatives of the function $f_i$. All time-varying communication networks are generated based on an Erdős-Rényi random graph model. The communication link between any two distinct agents is established independently with a predefined probability, which is specified later. The variation probability of the time-varying networks is set to 1 such that the instantaneous network topology $\mathcal{G}_k$ changes at each step and may frequently be disconnected. However, the expected network forms a connected graph. This property asymptotically guarantees ergodic connectivity, which simulate intermittent real-world communication environments (e.g., wireless sensor networks) while satisfying the sufficient connectivity conditions outlined in Assumption \ref{A3}. Fig. \ref{Fig0-3} provides an illustrative example of the dynamic evolution of time-varying network topologies across four distinct time steps, where red and yellow nodes denote Byzantine and reliable agents, respectively.

\begin{figure*}[!htp]
\begin{center}
\subfloat[Consensus error over iterations.]{\includegraphics[width=1.7in,height=1.2in]{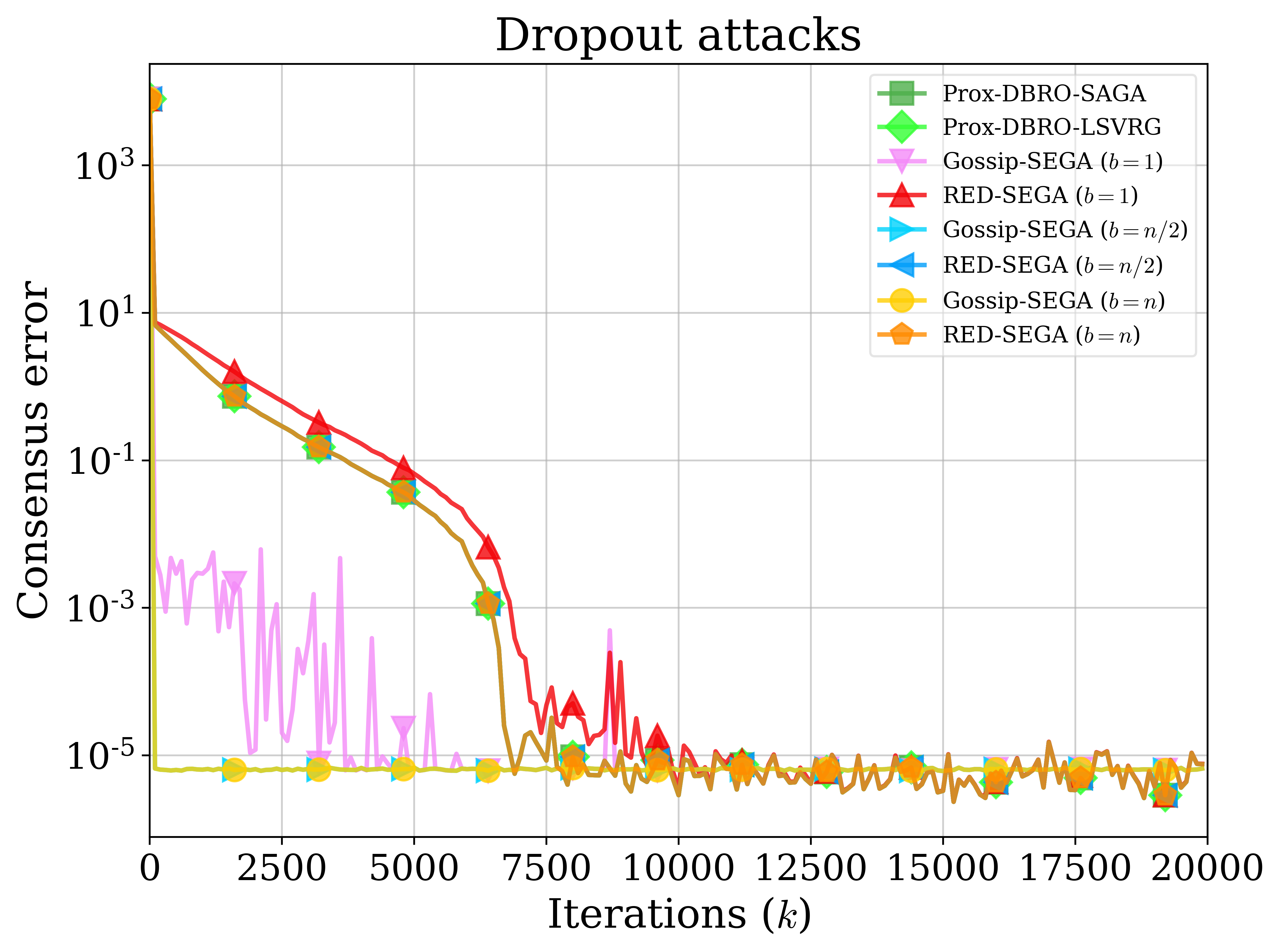}\label{Fig1-1}} \hfill
\subfloat[Residual over iterations.]{\includegraphics[width=1.7in,height=1.2in]{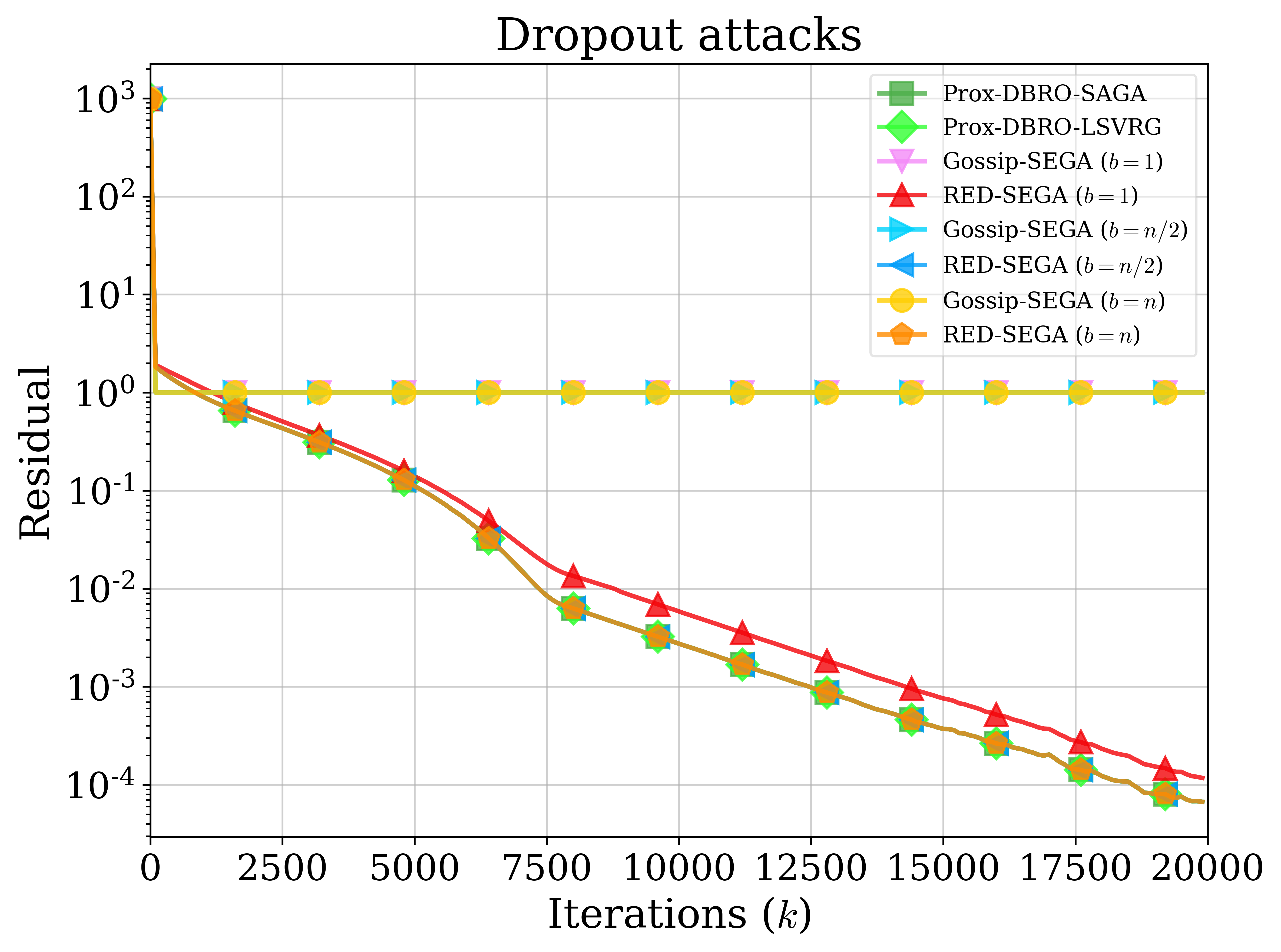}\label{Fig1-2}} \hfill
\subfloat[Consensus error over oracle calls.]{\includegraphics[width=1.7in,height=1.2in]{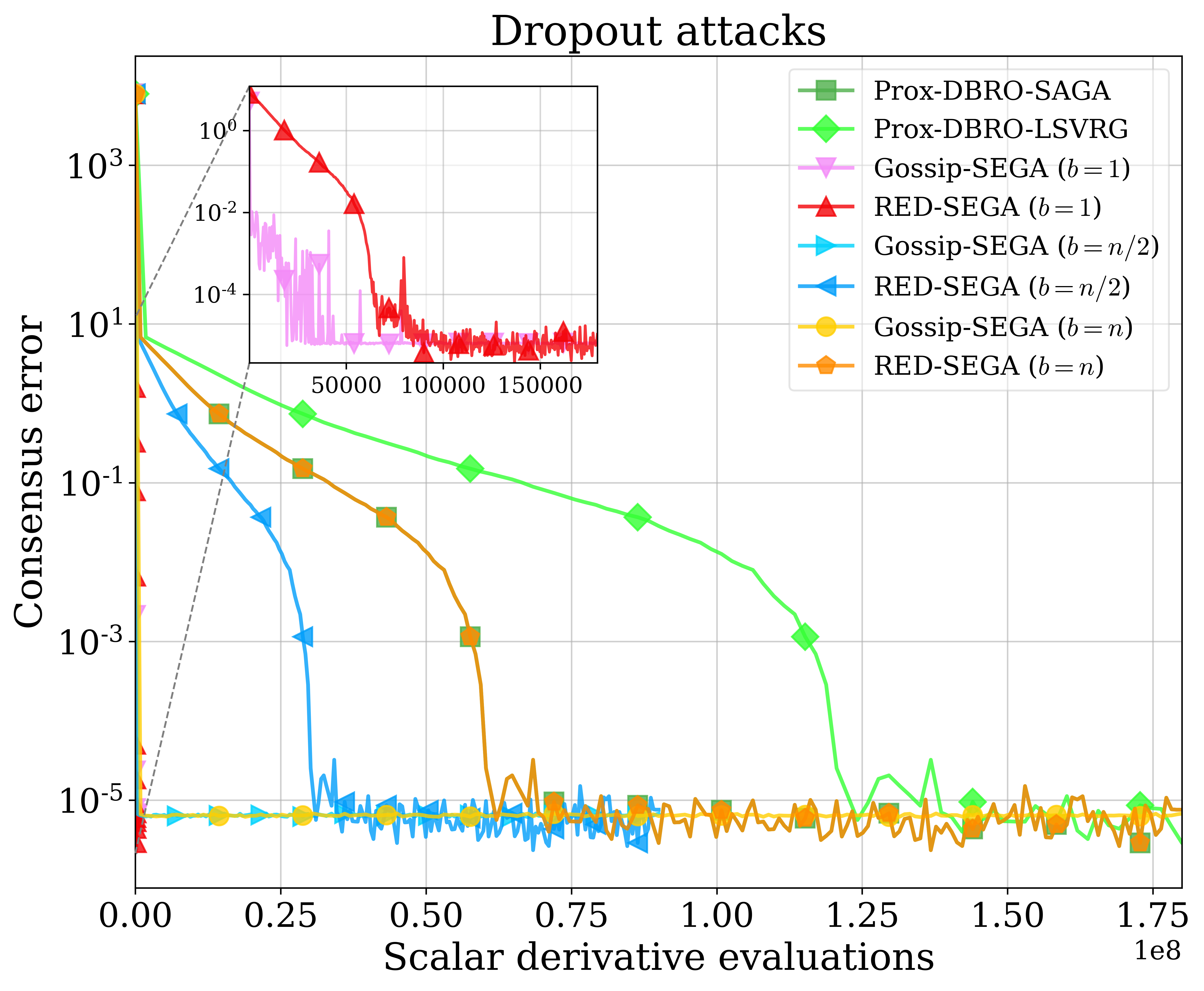}\label{Fig1-3}} \hfill
\subfloat[Residual over oracle calls.]{\includegraphics[width=1.7in,height=1.2in]{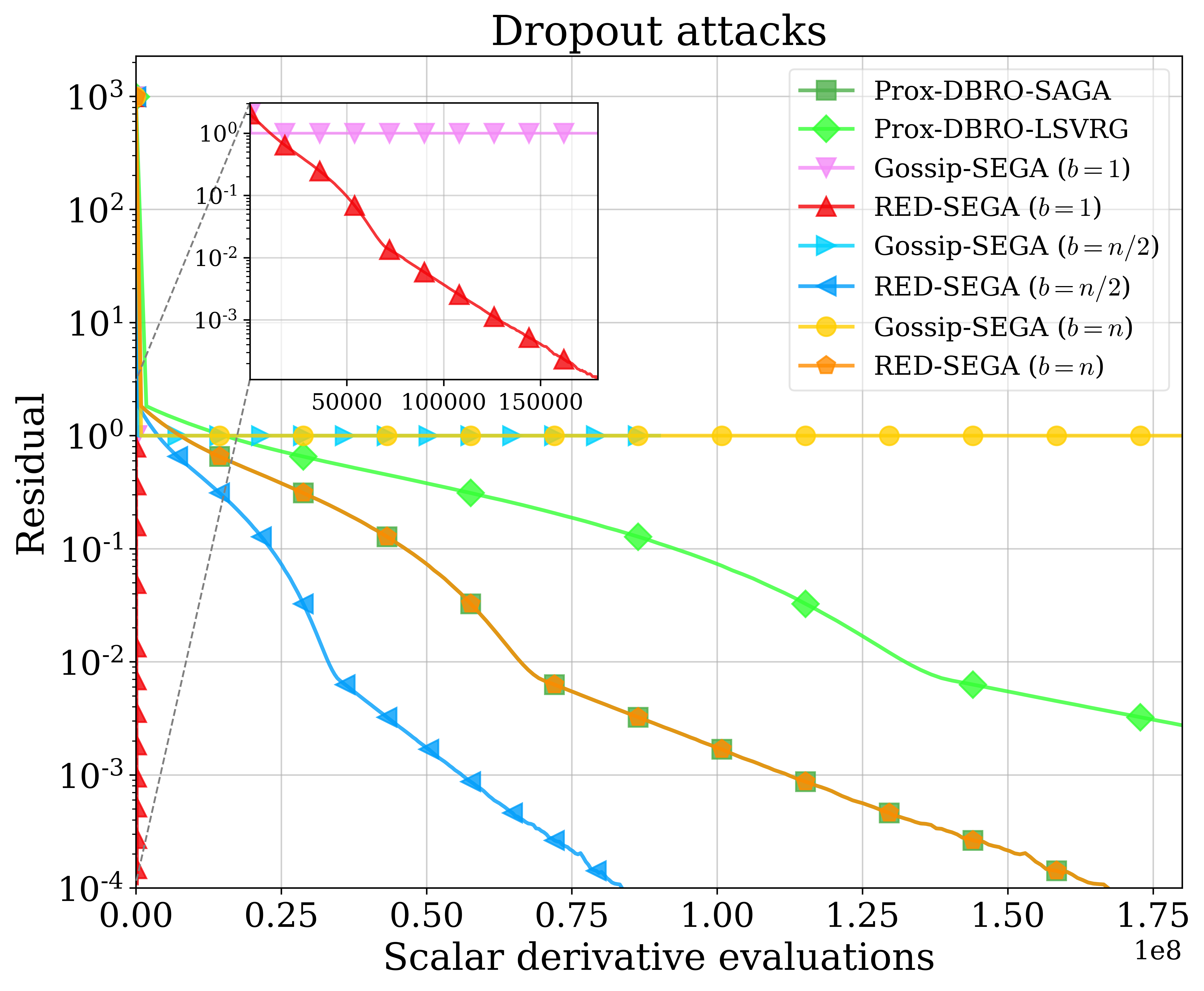}\label{Fig1-4}} \hfill
\end{center}
\caption{Convergence and consensus performance comparison under dropout attacks over time-varying networks comprising 1 Byzantine agents.}
\label{Fig1}
\end{figure*}

\begin{figure*}[!htp]
\begin{center}
\subfloat[Consensus error over iterations.]{\includegraphics[width=1.7in,height=1.2in]{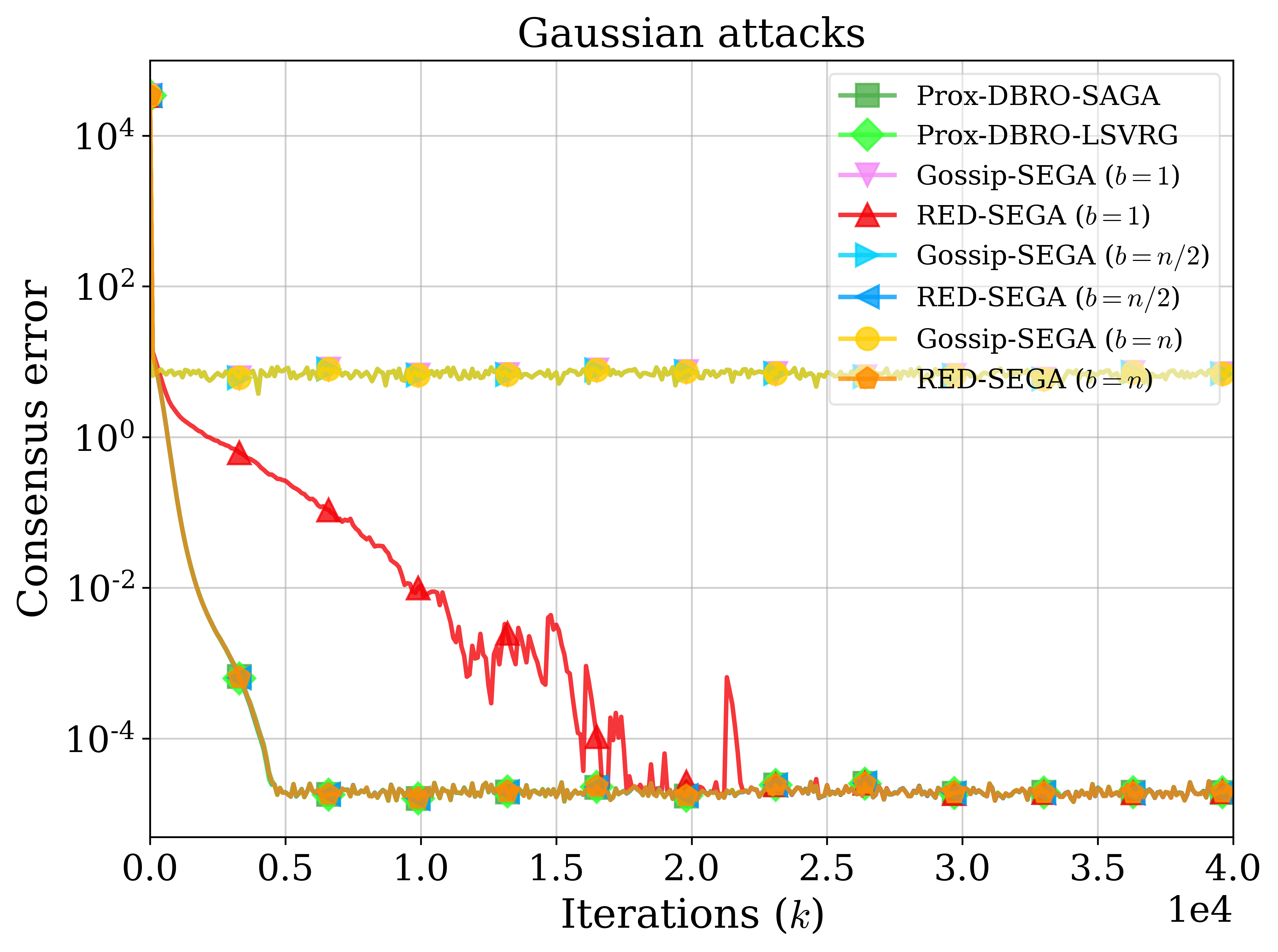}\label{Fig2-1}} \hfill
\subfloat[Residual over iterations.]{\includegraphics[width=1.7in,height=1.2in]{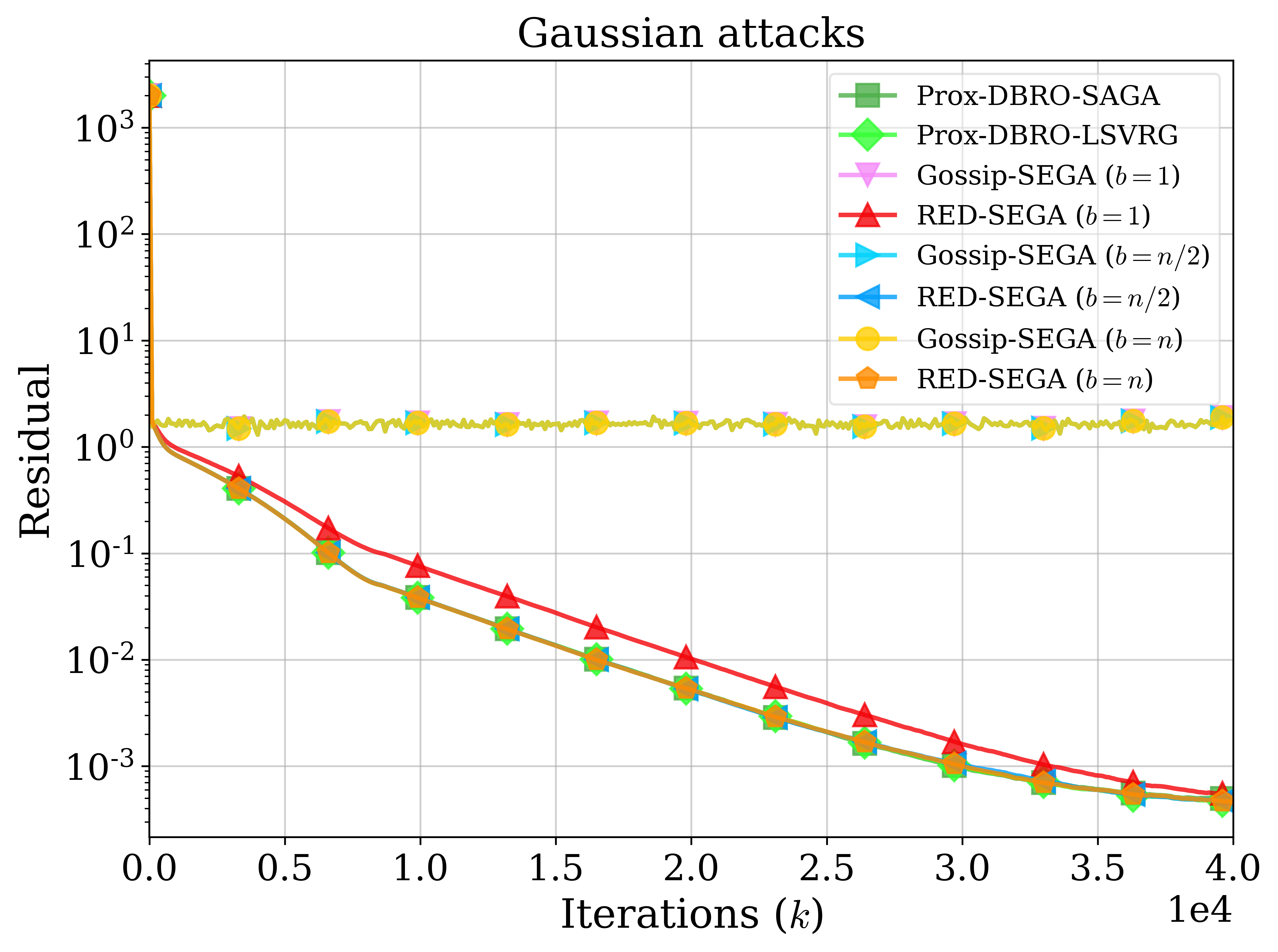}\label{Fig2-2}} \hfill
\subfloat[Consensus error over oracle calls.]{\includegraphics[width=1.7in,height=1.2in]{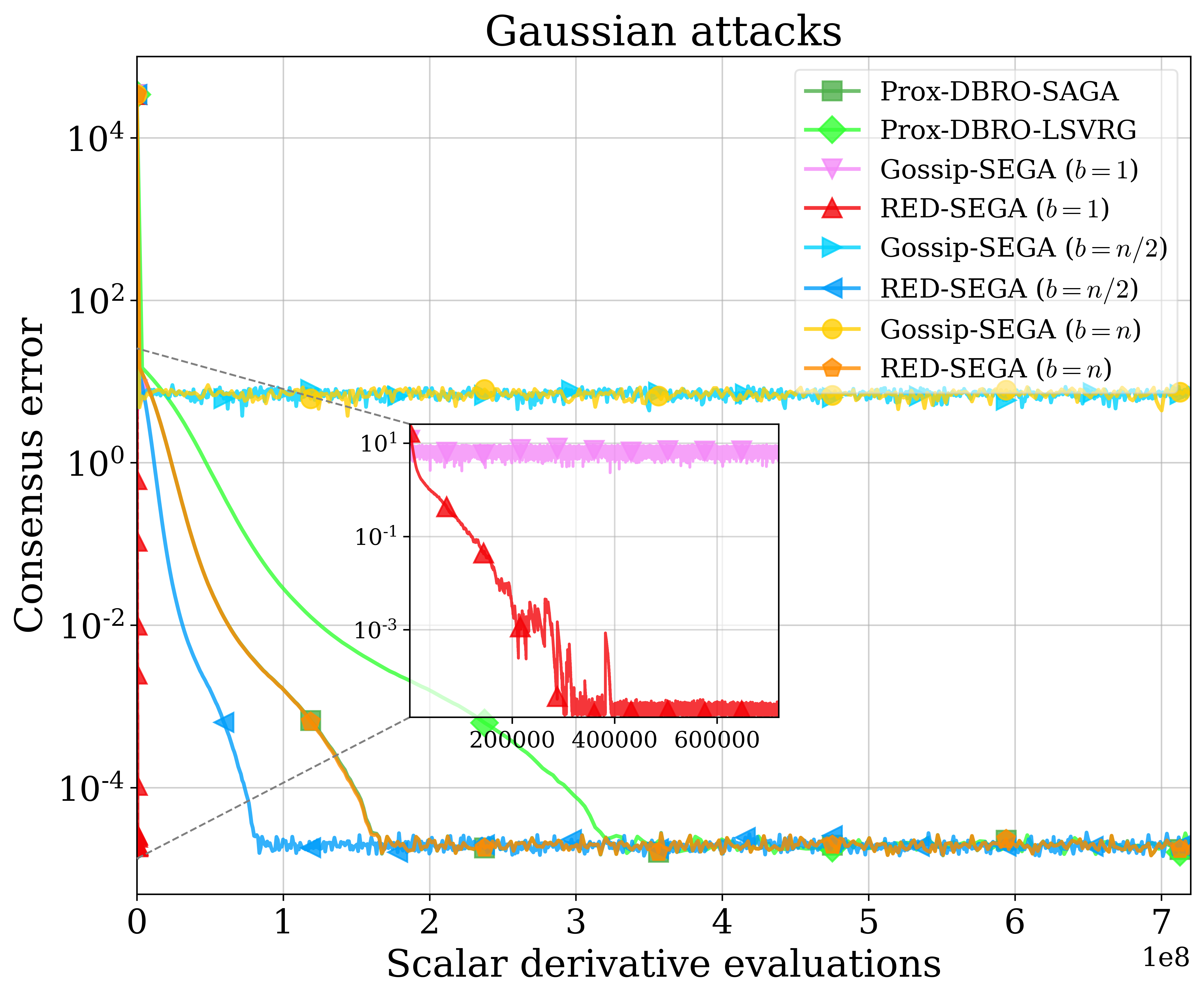}\label{Fig2-3}} \hfill
\subfloat[Residual over oracle calls.]{\includegraphics[width=1.7in,height=1.2in]{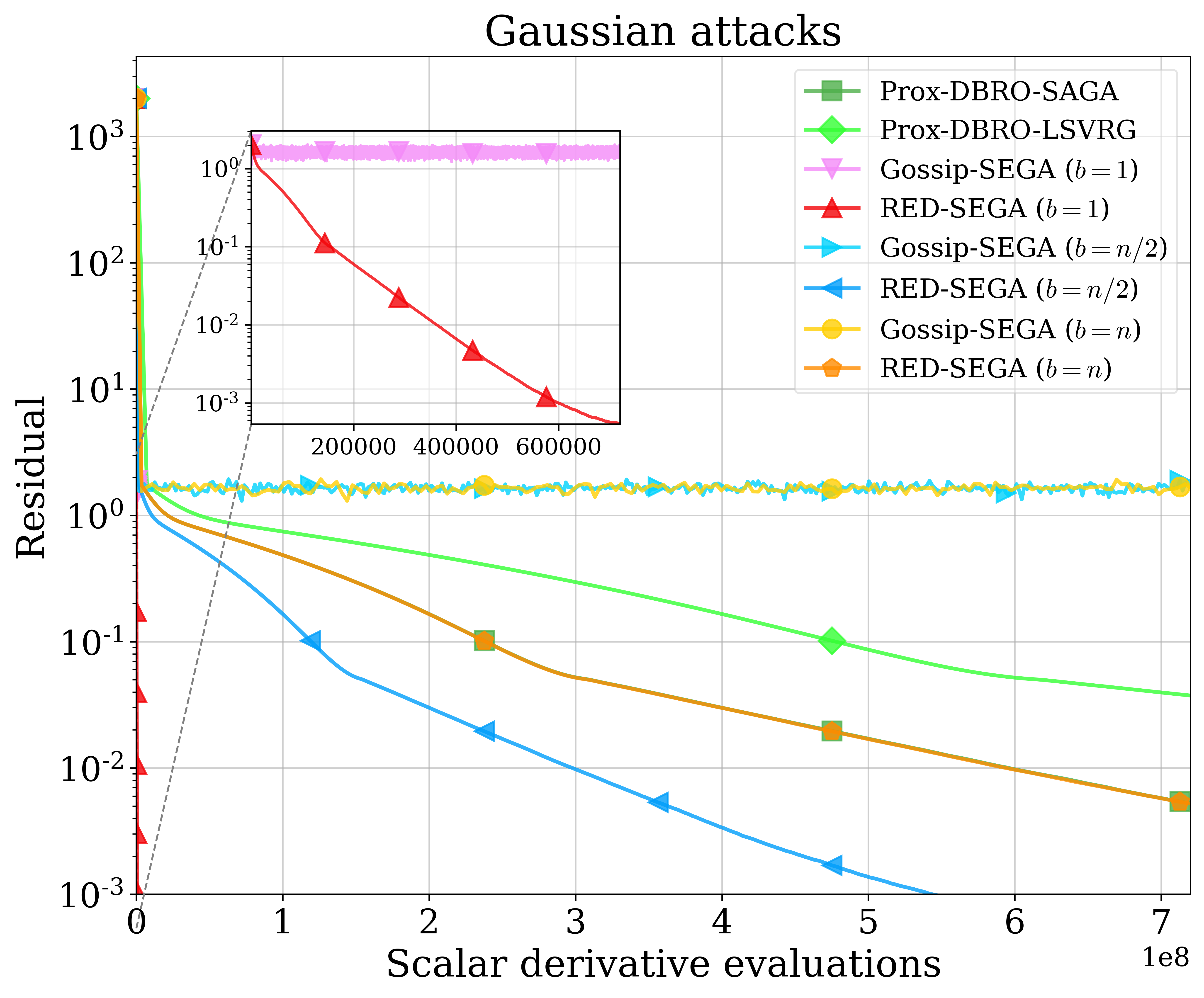}\label{Fig2-4}} \hfill
\end{center}
\caption{Convergence and consensus performance comparison under Gaussian attacks over time-varying networks comprising 2 Byzantine agents.}
\label{Fig2}
\end{figure*}

\begin{figure*}[!htp]
\begin{center}
\subfloat[Consensus error over iterations.]{\includegraphics[width=1.7in,height=1.2in]{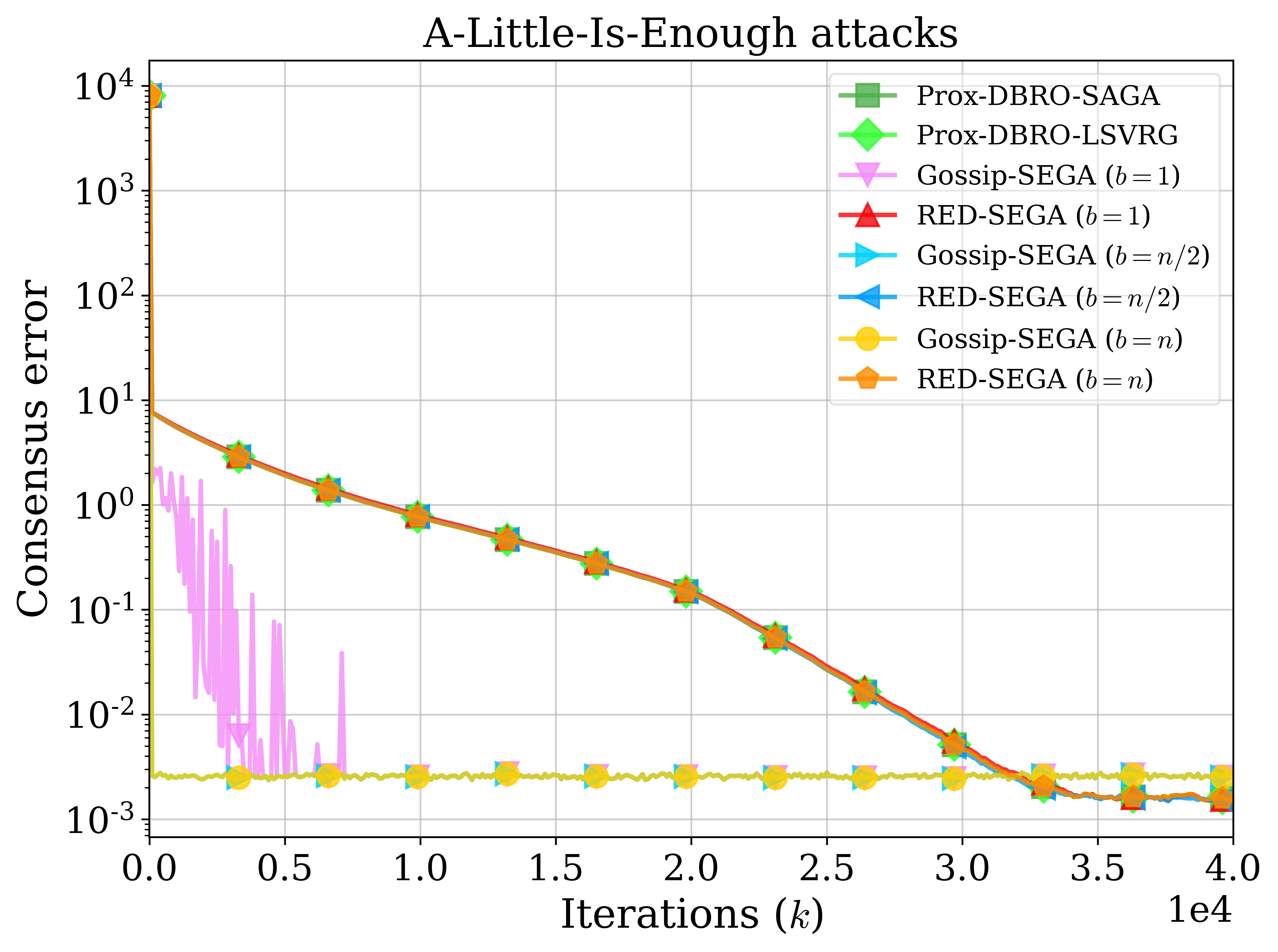}\label{Fig3-1}} \hfill
\subfloat[Residual over iterations.]{\includegraphics[width=1.7in,height=1.2in]{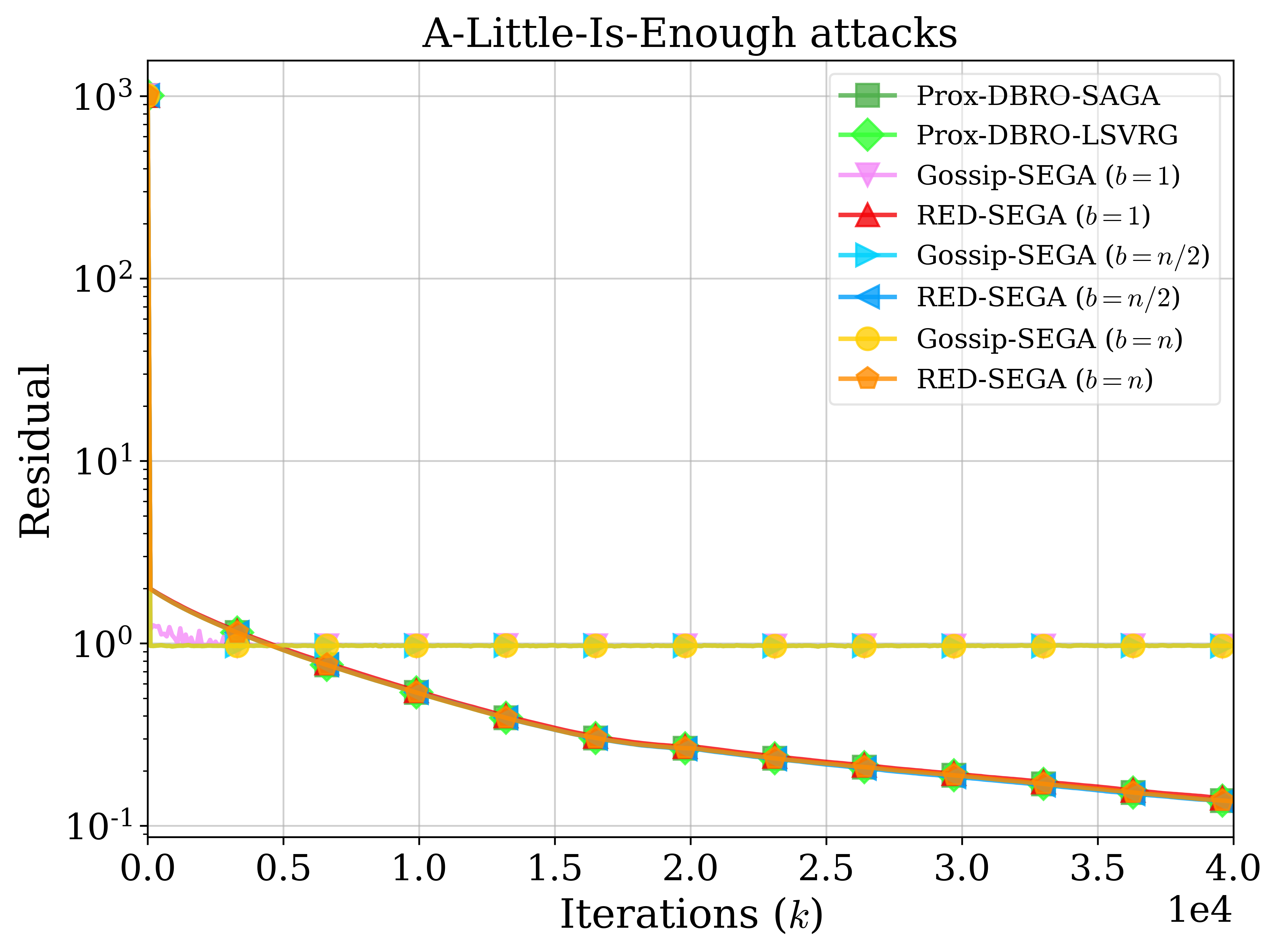}\label{Fig3-2}} \hfill
\subfloat[Consensus error over oracle calls.]{\includegraphics[width=1.7in,height=1.2in]{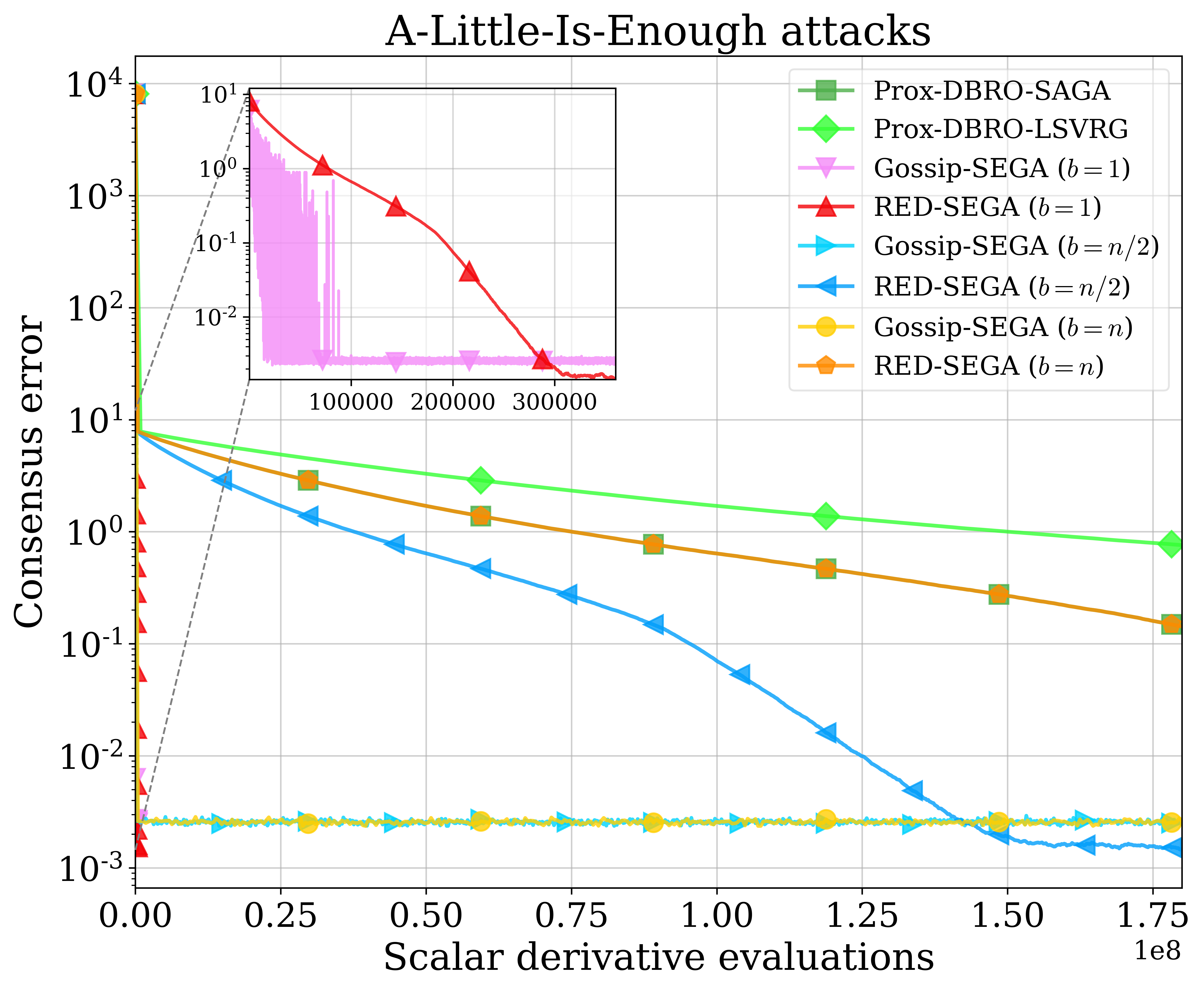}\label{Fig3-3}} \hfill
\subfloat[Residual over oracle calls.]{\includegraphics[width=1.7in,height=1.2in]{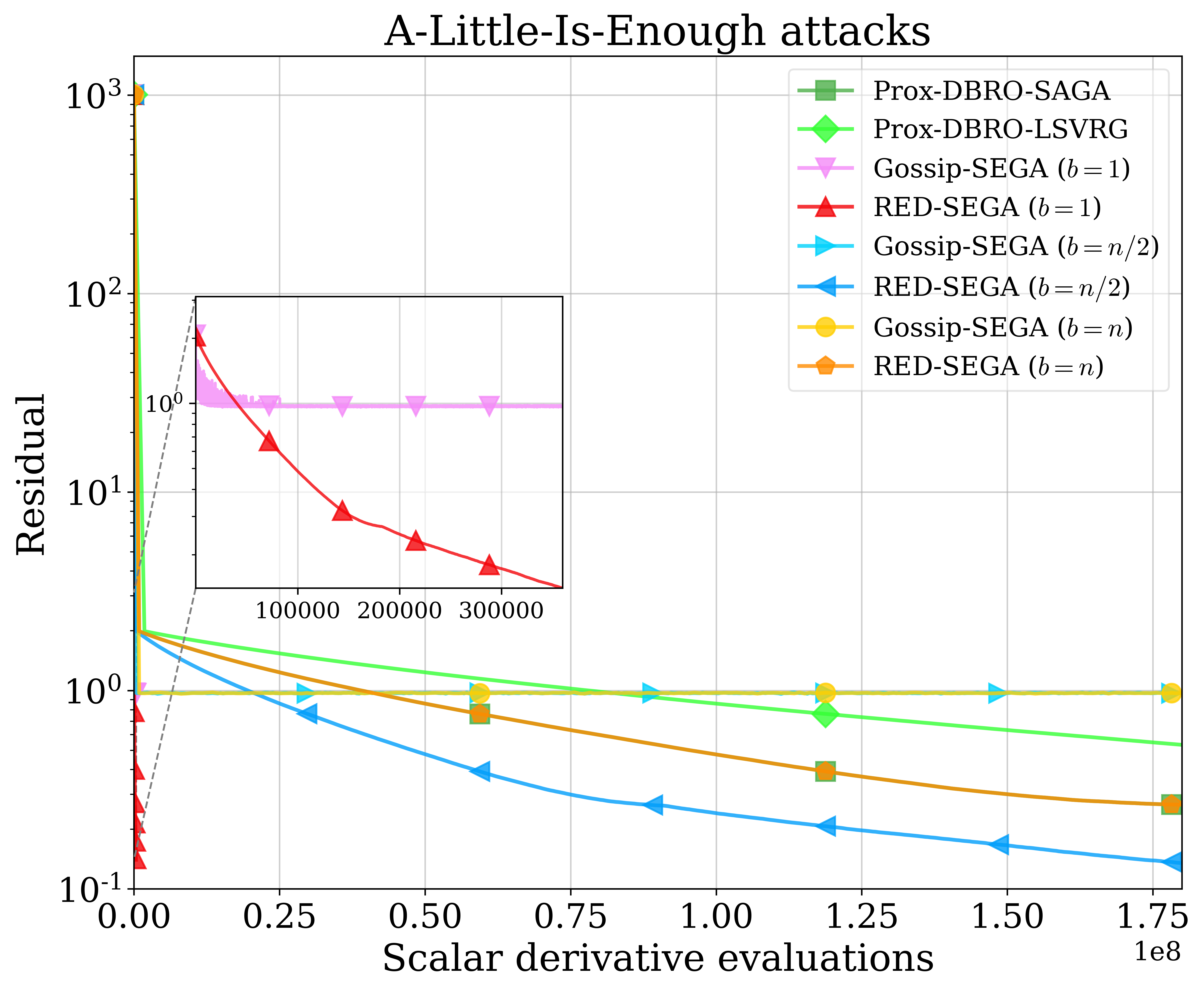}\label{Fig3-4}} \hfill
\end{center}
\caption{Convergence and consensus performance comparison under A-Little-Is-Enough attacks over time-varying networks comprising 3 Byzantine agents.}
\label{Fig3}
\end{figure*}

\begin{figure*}[!h]
\begin{center}
\subfloat[Consensus error over iterations.]{\includegraphics[width=1.7in,height=1.2in]{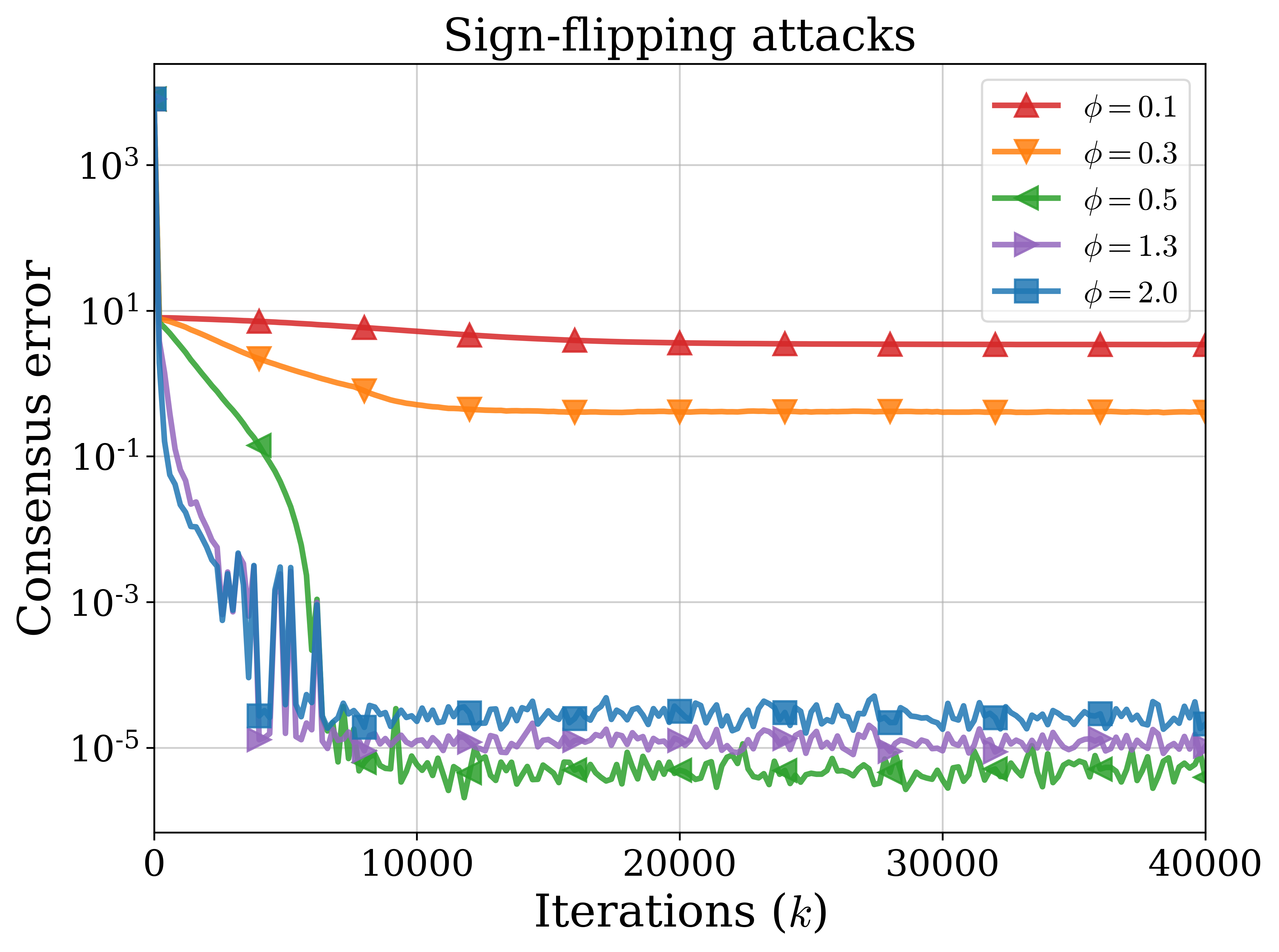}\label{Fig4-1}} \hfill
\subfloat[Residual over iterations.]{\includegraphics[width=1.7in,height=1.2in]{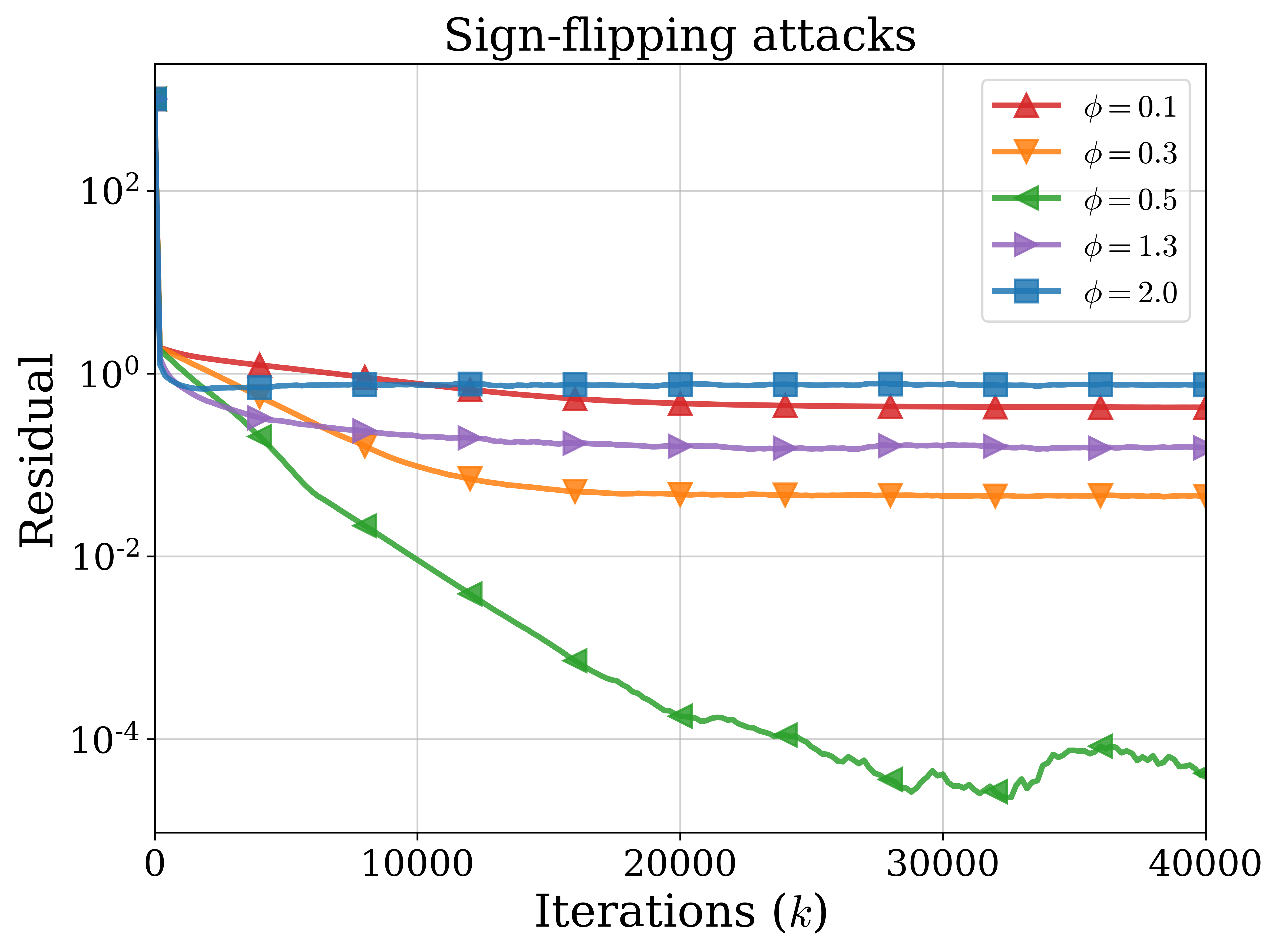}\label{Fig4-2}} \hfill
\subfloat[Consensus error over oracle calls.]{\includegraphics[width=1.7in,height=1.2in]{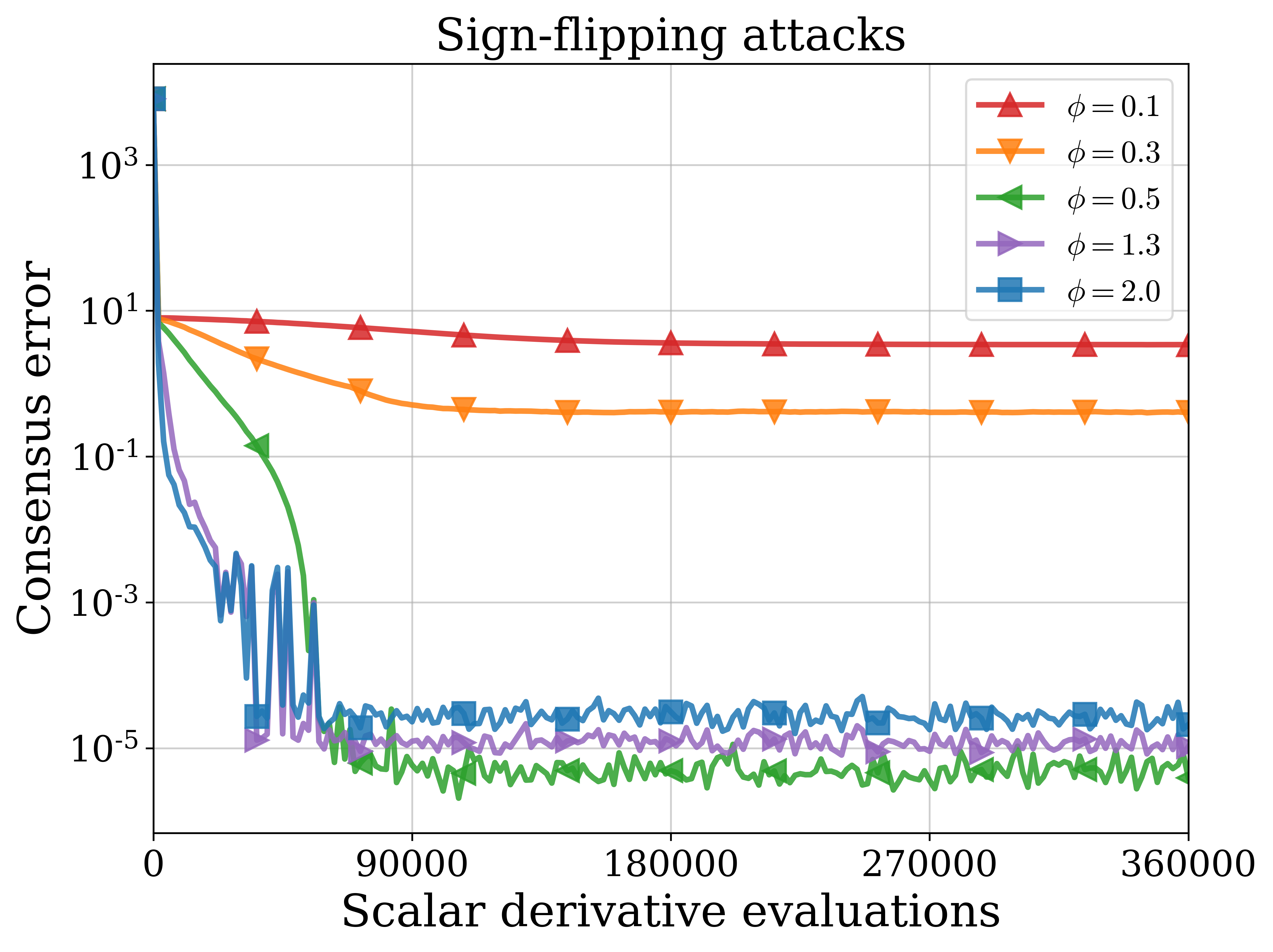}\label{Fig4-3}} \hfill
\subfloat[Residual over oracle calls.]{\includegraphics[width=1.7in,height=1.2in]{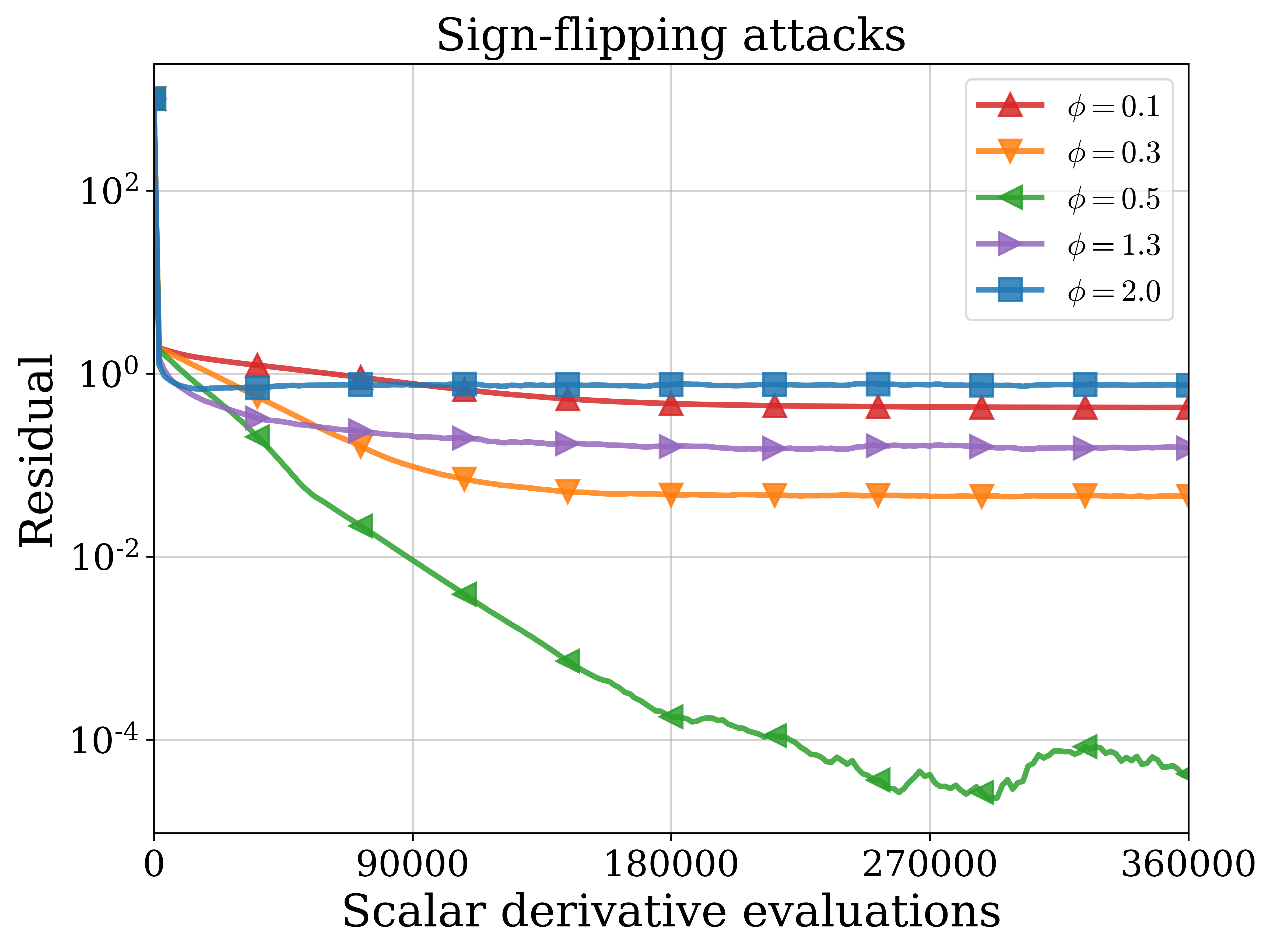}\label{Fig4-4}} \hfill
\end{center}
\caption{Ablation study of the penalty parameter $\phi$ for \textit{RED-SEGA} under sign-flipping attacks over time-varying networks comprising 2 Byzantine agents..}
\label{Fig4}
\end{figure*}
In this section, we compare the performance of the proposed algorithm \textit{RED-SEGA} with two baseline methods \textit{Prox-DBRO-SAGA} and \textit{Prox-DBRO-LSVRG} from recent literature \cite{Hu2025a} under four types of Byzantine attacks. In these evaluations, we scale the total number of agents from 10 to 30, with a corresponding increase in the number of Byzantine agents, which will be reported later. We evaluate the convergence performance of all tested algorithms using the residual, defined as $\sum\nolimits_{i \in \mathcal{R}} {\left\Vert{} {{x_{i,k}} - {{\tilde x}^*}} \right\Vert{}_{{W_i}}^2}$, to quantify the convergence error, where ${W_i} = {\mathbf{I}}$ in this case. This metric is plotted against both the number of iterations and oracle calls (scalar derivative evaluations) to illustrate the iteration complexity and overall computational efficiency, respectively. To ensure a fair comparison, the starting points of all tested algorithms are initialized identically. and the step-size of each tested algorithm was manually tuned to achieve its best performance in each case study.\\
\textbf{Dropout attacks:} Any Byzantine agent $j$, $j \in \mathcal{B}_i$, $i \in \mathcal{R}$, drops out of the communication by masking its actual state. The transmitted message is modeled as $\tilde{x}_{j,k}^i = \vartheta_{j,k} x_{j,k}$, where the masking variable $\xi_{j,k} \in \{0, 1\}$ follows a Bernoulli distribution. Specifically, the dropout probability $\mathbb{P}(\xi_{j,k} = 0)$ is randomly drawn from the interval $[0.5, 1.0]$ throughout the iterations. This case study is conducted over a time-varying network of $m = 10$ agents, including one Byzantine agent, where the connectivity probability of each communication link is set to 0.35. We adopt ${\ell _1}$-norm penalty for \textit{Prox-DBRO-SAGA}, \textit{Prox-DBRO-LSVRG}, and \textit{RED-SEGA}, where the penalty parameter for  these algorithms is set as 0.5. As shown in Figs. \ref{Fig1}: (a) and (b), \textit{RED-SEGA} and the baseline methods (\textit{Prox-DBRO-SAGA} and \textit{Prox-DBRO-LSVRG}) achieve almost identical convergence rates in terms of iterations. However, Figs. \ref{Fig1}: (c) and (d) demonstrate that \textit{RED-SEGA} ($b = 1$) and \textit{RED-SEGA} ($b = n/2$) exhibit significantly higher computational efficiency. In contrast, the baseline methods can only match the convergence speed of the full-coordinate version, \textit{RED-SEGA} ($b = n$), when evaluated against scalar derivative evaluations (oracle calls). Furthermore, it is worth noting that while \textit{Gossip-SEGA} attains the smaller consensus error in a fastest manner, it reaches an agreement at a state that severely deviates from the global optimal solution as evidenced by Figs. \ref{Fig1}: (b) and (d). This is undesirable, indicating that \textit{Gossip-SEGA} is misled by the Byzantine agent. These empirical observations align with our theoretical analysis, confirming that \textit{Gossip-SEGA} inherently lacks Byzantine resilience exhibited by \textit{RED-SEGA}.\\
\textbf{Gaussian attacks \cite{Hu2025}:} Any Byzantine agent $j$, $j \in \mathcal{B}_i$, $i \in \mathcal{R}$, corrupts its true model before transmitting it to its reliable neighbor $i$. The falsified vector is constructed as $\tilde{x}_{j,k}^i = x_{j,k} + \nu_{j,k}$, wherein $\nu_{j,k}$ represents the random noise matrix identically sampled from a predefined Gaussian distribution $\mathcal{N}(0, \sigma^2 \textbf{I})$. This case study is evaluated on a time-varying network comprising a total of $m = 20$ agents, including $2$ Byzantine agents, where the connectivity probability of each communication link is set to 0.25. We adopt ${\ell _2}$-norm penalty for \textit{Prox-DBRO-SAGA}, \textit{Prox-DBRO-LSVRG}, and \textit{RED-SEGA}, where the penalty parameter for  these algorithms is set as 0.03. As shown in Figs. \ref{Fig2}: (a) and (b), although the baseline methods achieve network consensus slightly faster, \textit{RED-SEGA} and the baseline algorithms (\textit{Prox-DBRO-SAGA} and \textit{Prox-DBRO-LSVRG}) exhibit comparable linear convergence rates in terms of the residual reduction over iterations. Consistent with the previous scenario, Figs. \ref{Fig2}: (c) and (d) demonstrate that \textit{RED-SEGA} ($b = 1$) and \textit{RED-SEGA} ($b = n/2$) significantly outperform the baseline methods in computational efficiency. In contrast, the baseline methods merely match the convergence efficiency of the full-coordinate version, \textit{RED-SEGA} ($b = n$), when evaluated against scalar derivative evaluations (oracle calls). Furthermore, a striking difference from the dropout attack scenario is observed for \textit{Gossip-SEGA}. As evidenced by the flat horizontal trajectories across all subfigures, \textit{Gossip-SEGA} fails to achieve both network consensus and optimality, remaining stagnant at a high error plateau. This catastrophic failure further corroborates that \textit{Gossip-SEGA} is vulnerable to continuous Gaussian perturbations, thereby validating the essential role of the resilient aggregation mechanism equipped in \textit{RED-SEGA}. \\
\textbf{A-Little-Is-Enough attacks \cite{Baruch2019}:} A Byzantine neighbor $j$, $j \in \mathcal{B}_i$, $i \in \mathcal{R}$, crafts its message to the reliable neighbor $i$ as $\tilde{x}_{j,k}^i = \hat{\mu}_{\mathcal{R}_{i,k}} - \pi_j \cdot \hat{\sigma}_{\mathcal{R}_{i,k}}$. In this formulation, $\hat{\mu}_{\mathcal{R}_{i,k}}$ and $\hat{\sigma}_{\mathcal{R}_{i,k}}$ denote the empirical mean and standard deviation statistically computed from the local models of reliable agents. To bypass detection defenses, the tuning parameter $\pi_j$ is tightly bounded by the cumulative distribution function of the standard normal distribution. This case study is evaluated on a time-varying network comprising a total of $m = 30$ agents, including $3$ Byzantine agents, where the connectivity probability of each communication link is set to 0.15. We adopt ${\ell _\infty}$-norm penalty for \textit{Prox-DBRO-SAGA}, \textit{Prox-DBRO-LSVRG}, and \textit{RED-SEGA}, where the penalty parameter for these algorithms is set as 10. As shown in Figs. \ref{Fig3}: (a) and (b), \textit{RED-SEGA} and the baseline algorithms (\textit{Prox-DBRO-SAGA} and \textit{Prox-DBRO-LSVRG}) effectively withstand this statistically crafted attack, maintaining comparable linear convergence rates in terms of residual reduction over iterations. Consistent with the previous scenarios, Figs. \ref{Fig3}: (c) and (d) clearly illustrate that the gradient-sketching variants, \textit{RED-SEGA} ($b = 1$) and \textit{RED-SEGA} ($b = n/2$), significantly surpass the baseline methods in overall computational efficiency. The baseline methods, once again, merely match the convergence speed of the full-coordinate version, \textit{RED-SEGA} ($b = n$), when measured against scalar derivative evaluations (oracle calls). Notably, the deceptive nature of the attack completely incapacitates \textit{Gossip-SEGA}. As evidenced by the flat horizontal residual trajectories in Figs. \ref{Fig3}: (b) and (d), \textit{Gossip-SEGA} fails entirely to optimize, remaining stagnant at a high error plateau. This failure confirms that while subtle, statistically constrained malicious messages easily bypass naive gossip averaging, the norm-penalized resilient aggregation mechanism equipped in \textit{RED-SEGA} maintains robust and effective defensive capabilities.

The numerical results in Figs. \ref{Fig1}-\ref{Fig3} demonstrate that \textit{RED-SEGA} achieves consensus and convergence performance comparable to the adapted \textit{Prox-DBRO-VR} baselines, namely \textit{Prox-DBRO-SAGA} and \textit{Prox-DBRO-LSVRG}. Next, we conduct an ablation study under sign-flipping attacks to empirically validate the theoretical trade-off regarding the penalty parameter $\phi$, as established in Theorem \ref{T2}.\\
\textbf{Sign-flipping attacks \cite{Li2019n}:} A Byzantine agent $j$, $j \in \mathcal{B}_i$, $i \in \mathcal{R}$, transmits a maliciously scaled and inverted vector to its reliable neighbor $i$. The fake message is generated via $\tilde{x}_{j,k}^i = - \rho_j \cdot \bar{x}_i$, where $\bar{x}_i$ denotes the average state calculated among the reliable neighbors (including agent $i$ itself), and $\rho_j > 0$ serves as a scaling factor that magnifies the adversarial deviation. In this ablation study, we consider the multi-coordinate version of \textit{RED-SEGA} ($b=n/2$), which is evaluated on a time-varying network comprising a total of $m = 20$ agents, including $2$ Byzantine agents, where the connectivity probability of each communication link is set to 0.45. We adopt ${\ell _\infty}$-norm penalty for \textit{Prox-DBRO-SAGA}, \textit{Prox-DBRO-LSVRG}, and \textit{RED-SEGA}. As illustrated in Fig. \ref{Fig4}, an inadequate penalty, i.e., $\phi \le 0.3$, fails to bound the malicious perturbations, leading to a persistent consensus error (Figs. \ref{Fig4}: (a) and (c)) and a large convergence error (Figs. \ref{Fig4}: (b) and (d)). Conversely, an excessively large penalty, i.e., $\phi \ge 1.3$, forces a rapid consensus but dominates the objective function's true gradient, thereby amplifying the steady-state convergence error. Moderate values of $\phi$, i.e., $\phi = 0.5$, strike a balance. Under this setting, \textit{RED-SEGA} achieves a tight consensus among reliable agents while successfully minimizing the convergence error to the lowest achievable convergence error. These empirical observations corroborate the theoretical insights and provide practical guidance for hyperparameter tuning.

In summary, the experimental results demonstrate that \textit{RED-SEGA} not only exhibits enhanced resilience compared to \textit{Gossip-SEGA}, which corroborates the theoretical analysis in Section \ref{sec4-3}, but also achieves superior computational efficiency over baseline methods \textit{Prox-DBRO-SAGA} and \textit{Prox-DBRO-LSVRG} through single- or multi-coordinate sampling.

\subsection{Decentralized Image Deblurring Under Byzantine Attacks}\label{sec6-2}
\begin{figure}[!htp]
\begin{center}
\begin{tabular}{c}
\includegraphics[width=0.65in,height=0.65in]{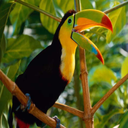} \\
(a)
\end{tabular}
\begin{tabular}{c}
\includegraphics[width=0.65in,height=0.65in]{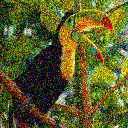} \\
(b)
\end{tabular}
\begin{tabular}{c}
\includegraphics[width=0.65in,height=0.65in]{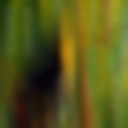} \\
(c)
\end{tabular}
\begin{tabular}{c}
\includegraphics[width=0.65in,height=0.65in]{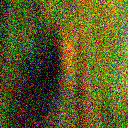} \\
(d)
\end{tabular}
\caption{A clean image from Set5 dataset \cite{Bevilacqua2012} and its blurred versions: (a) clean; (b) Gaussian noise; (c)
motion blurring; (d) Gaussian noise and motion blurring.}
\label{Fig8}
\end{center}
\end{figure}

In the second case study, we aim to show the performance of the multi-coordinate FO and ZO versions of \textit{RED-SEGA} via an image deblurring task. In this case study, a time-varying network containing 4 reliable agents and 1 Byzantine agent is generated, and the Byzantine attacks studied in Section \ref{sec6-1} are considered. In practice, a clean image can be simultaneously degraded by motion deblurring (arising from relative sensor-scene displacement during exposure) and Gaussian noise. Specifically, according to \cite{Ren2019,Yang2025,Lu2026}, the linear degradation process is expressed as
\begin{equation}\label{E6-2-1}
{y_i} = {H_i}{\tilde x} + {\tilde n_i},
\end{equation}
where $\tilde x$, ${\tilde n_i}$, $H_i$, and $y_i$ are the clean image, Gaussian noise, a degradation operator, and a degraded image, respectively. As an example, Figs. \ref{Fig8}: (a)-(d) illustrate four states of an image from the Set5 dataset \cite{Bevilacqua2012}: (a) clean; (b) with Gaussian noise only; (c) with motion blurring only; and (d) with both Gaussian noise and motion blurring. We vectorize the pixel matrix of RGB images and convert the ﬁlter into Toeplitz matrices, which allows us to perform convolution through matrix multiplication. Owing to this, the image deblurring task can be formulated as a least-squared problem with a nonsmooth sparse term
\begin{equation}\label{E6-2-2}
\mathop {\min }\limits_{{\tilde x}} \sum\limits_{i \in \mathcal{R}} {\left\| {{H_i}{\tilde x} - {y_i}} \right\|^2} + \beta {\left\| {\tilde x} \right\|_1},
\end{equation}
where $\beta > 0 $ is the regularized constant and the regularized term is used to eliminate small and random fluctuations caused by noise.

\begin{figure}[!htp]
\begin{center}
\begin{tabular}{c}
\includegraphics[width=0.65in,height=0.65in]{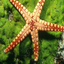} \\
\includegraphics[width=0.65in,height=0.65in]{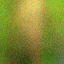} \\
\includegraphics[width=0.65in,height=0.65in]{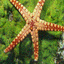} \\
\includegraphics[width=0.65in,height=0.65in]{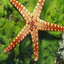} \\
(a)
\end{tabular}
\begin{tabular}{c}
\includegraphics[width=0.65in,height=0.65in]{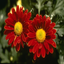} \\
\includegraphics[width=0.65in,height=0.65in]{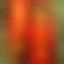} \\
\includegraphics[width=0.65in,height=0.65in]{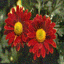} \\
\includegraphics[width=0.65in,height=0.65in]{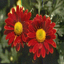} \\
(b)
\end{tabular}
\begin{tabular}{c}
\includegraphics[width=0.65in,height=0.65in]{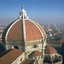} \\
\includegraphics[width=0.65in,height=0.65in]{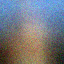} \\
\includegraphics[width=0.65in,height=0.65in]{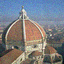} \\
\includegraphics[width=0.65in,height=0.65in]{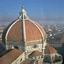} \\
(c)
\end{tabular}
\begin{tabular}{c}
\includegraphics[width=0.65in,height=0.65in]{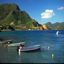} \\
\includegraphics[width=0.65in,height=0.65in]{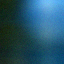} \\
\includegraphics[width=0.65in,height=0.65in]{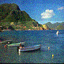} \\
\includegraphics[width=0.65in,height=0.65in]{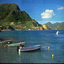} \\
(d)
\end{tabular}
\caption{Deblurring performance of ZO and FO versions \textit{RED-SEGA} under various Byzantine attacks: (a) dropout attacks; (b) Gaussian attacks; (c) A-Little-Is-Enough attacks; (d) sign-flipping attacks. The first to the fourth rows show the clean images, degraded images, deblurring results from the ZO version of \textit{RED-SEGA}, and deblurring results from the FO version of \textit{RED-SEGA}, respectively.}
\label{Fig9}
\end{center}
\end{figure}

\begin{table}[!htp]
\centering
\caption{Performance on image deblurring under various Byzantine attacks: (a) dropout attacks; (b) Gaussian attacks; (c) A-Little-Is-Enough attacks; (d) sign-flipping attacks.}
\label{Tab2}
\scalebox{0.95}{
\begin{tabular}{c|cccc}
\toprule
\hline
Clean images & (a) & (b) & (c) & (d)  \\
\hline 
PSNR   & inf  & inf  & inf  & inf  \\
SSIM   & 1.00 & 1.00 & 1.00 & 1.00 \\
Loss   & 1149.95 & 1113.81 & 477.66 &  888.08 \\ 
\hline
Degraded images & (a) & (b) & (c) & (d) \\
\hline
PSNR   & 15.21   & 15.81  & 17.43 & 17.76 \\
SSIM   & 0.6781  & 0.5884 & 0.5618 & 0.5750 \\
Loss   & 1569.06 & 1374.13 & 746.22 & 1068.69  \\
\hline
\multirow{2}{*}{\makecell{Recovered images by \\ \textit{RED-SEGA} (ZO version)}}
       & \multirow{2}{*}{(a)} & \multirow{2}{*}{(b)} & \multirow{2}{*}{(c)} & \multirow{2}{*}{(d)}   \\
       &  &  &  &    \\
\hline
PSNR   & 25.68  & 25.78 & 26.14  & 26.99 \\
SSIM   & 0.9441  & 0.8966 & 0.9139 & 0.9028  \\
Loss   & 1071.34  & 1032.18 & 447.24 & 822.39  \\
Consensus error  & 8.97e-04 & 5.32e-03 & 7.16e-03 & 8.53e-04 \\
\hline
\multirow{2}{*}{\makecell{Recovered images by \\ \textit{RED-SEGA} (FO version)}}
       & \multirow{2}{*}{(a)} & \multirow{2}{*}{(b)} & \multirow{2}{*}{(c)} & \multirow{2}{*}{(d)}  \\
       &  &  &  &    \\
\hline
PSNR             & 27.44 & 27.56  & 27.88 & 29.12 \\
SSIM             & 0.9795 & 0.9491 & 0.9808 & 0.9809  \\
Loss             & 1029.21 & 994.73 & 407.58 & 781.17 \\
Consensus error  & 3.49e-05 & 1.51e-04 & 6.87e-04 & 3.11e-05 \\
\hline
\bottomrule
\end{tabular}}
\end{table}
We assume a total number of $5$ agents, including 4 reliable agents and 1 Byzantine agent over a time-varying network with dynamic connectivity. We randomly select 4 images from the BSDS300 dataset \cite{Martin2001}, which are resized into $ 64 \times 64 \times 3$ RGB images. We use four metrics: peak signal-to-noise ratio (PSNR), structural similarity index (SSIM), loss, and evaluate the performance of two tested algorithms. As shown in Fig. \ref{Fig9}, both the FO and ZO versions of \textit{RED-SEGA} of \textit{RED-SEGA} demonstrate good deblurring performance under various Byzantine attacks, with the FO version exhibiting superior performance to the ZO version. This is verified by the results in Table \ref{Tab2}.

\section{Conclusions}\label{sec7}
This paper aims to resolve a class of SRM problems with non-decomposable objectives and potentially unavailable gradients over time-varying networks, which are prevalent in data-scarce yet high-dimensional tasks. To this end, we developed a localized gradient-sketching technique, which attains variance reduction via random dimension sampling. Based on this, a decentralized preconditioned stochastic proximal gradient algorithm with variance reduction was developed based on gradient sketching, dubbed \textit{Gossip-SEGA}. To enhance the resilience of \textit{Gossip-SEGA} against Byzantine agents, \textit{RED-SEGA} is further developed via substituting the naive aggregation via weighted average with a resilient aggregation based on a norm-penalized approximation. Theoretical analysis proved that \textit{RED-SEGA} converges to a neighborhood of the global optimal solution at a linear rate, where the radius of the convergence error is explicitly characterized. The effectiveness and practicality of \textit{Gossip-SEGA} are validated via resolving a constrained least-squared problem and an image deblurring problem under various Byzantine attacks.

\bibliographystyle{IEEEtran}
\bibliography{RED-SEGA}

\clearpage
\setcounter{page}{1}
\onecolumn
\section*{Supplementary Material}\label{sec8}
\subsection{Proof of Theorem \ref{T1}}\label{sec8-1}
The optimal solution to \text{DRC} (\ref{E3-3-2}) satisfies the optimality condition
\begin{equation}\label{E8-1-1}
{0_n} \in \nabla {f_i}\left( {x_i^*} \right) + \partial {r_i}\left( {x_i^*} \right) + \phi \sum\limits_{j \in {{\bar {\mathcal{R}}}_i},i < j} {{\omega _{\left( {i,j} \right)}}\partial } {\left\| {x_i^* - x_j^*} \right\|_p} - \phi \sum\limits_{j \in {{\bar {\mathcal{R}}}_i},i > j} {{\omega _{\left( {i,j} \right)}}\partial } {\left\| {x_j^* - x_i^*} \right\|_p}.
\end{equation}
According to the definition of the subdifferential for a general $\ell_p$-norm, we have ${\partial}{\left\| {x_i^* - x_j^*} \right\|_p} = \left\{ {{y_{ij}} \in {\mathbb{R}^n} \mid \left\langle {{y_{ij}}, x_i^* - x_j^*} \right\rangle  = {{\left\| {x_i^* - x_j^*} \right\|}_p}, {{\left\| {{y_{ij}}} \right\|}_d} \le 1} \right\}$, where ${\left\| \cdot \right\|_d}$ is the dual norm of ${\left\| \cdot \right\|_p}$ satisfying $\left( {1/p} \right) + \left( {1/d} \right) = 1$. Thus, there exist $\partial r_i^* \in \partial r_i\left( {x_i^*} \right)$ and ${\tilde y_{ij}} \in {\partial}{\left\| {x_i^* - x_j^*} \right\|_p}$. For each edge $e = \left( {i,j} \right) \in {\bar {\mathcal{E}}_{\mathcal{R}}}$ with $i < j$, we denote the edge subgradient as ${{\bar y}_e}: = {\tilde y}_{ij}$. We further define the combined consensus subgradient for agent $i$ as ${{{\bar{\mathbf{y}}}_i}}: = \sum\nolimits_{j \in {{\bar {\mathcal{R}}}_i},i < j} {{\omega _{\left( {i,j} \right)}}{{\bar y}_e}}  - \sum\nolimits_{j \in {{\bar {\mathcal{R}}}_i},i > j} {{\omega _{\left( {i,j} \right)}}{{\bar y}_e}} $ such that we have
\begin{equation}\label{E8-1-2}
0_n = \nabla {f_i}\left( {x_i^*} \right) + \partial r_i^* + \phi {{{\bar{\mathbf{y}}}_i}}.
\end{equation}
Under Assumption \ref{A1}, we know that ${x^*} = {\mathbf{col}}{\left\{ {x_i^*} \right\}_{i \in \mathcal{R}}}$ is a unique global optimal solution to \text{DRC} (\ref{E3-3-2}). It is sufficient to prove that the optimal solution ${\tilde x^*}$ satisfies (\ref{E8-1-2}), i.e.,
\begin{equation}\label{E8-1-3}
0_n = \nabla {f_i}\left( {{{\tilde x}^*}} \right) + \partial {r_i}\left( {{{\tilde x}^*}} \right) + \phi {{{\bar{\mathbf{y}}}_i}}.
\end{equation}
To rigorously analyze the general $n$-dimensional case without element-wise decomposition, we denote $\Psi : = {\mathbf{col}}{\left\{ {{\psi _i}} \right\}_{i \in \mathcal{R}}} \in \mathbb{R}^{n|\mathcal{R}|}$ with ${\psi _i}: = {\nabla {f_i}\left( {{{\tilde x}^*}} \right) + \partial {r_i} \left( {{{\tilde x}^*}} \right)}$, and let $\mathbf{y} = \mathbf{col}\{ \bar{y}_e \}_{e \in \bar{\mathcal{E}}_{\mathcal{R}}} \in \mathbb{R}^{n|\bar{\mathcal{E}}_{\mathcal{R}}|}$ be the vector constructed by stacking all consensus subgradients. Since the $i$-th block of the matrix-vector product $(\Omega \otimes I_n) \mathbf{y}$ exactly equals ${{{\bar{\mathbf{y}}}_i}}$, the node-wise optimality condition (\ref{E8-1-3}) over the \textit{average} network can be compactly rewritten using the Kronecker product into a global linear system as
\begin{equation}\label{E9-5-1}
\phi (\Omega \otimes I_n) {\mathbf{y}} + \Psi = {0_{n|\mathcal{R}|}},
\end{equation}
which necessitates at least one solution $\mathbf{y}$ satisfying the dual norm constraint ${\left\| {{{\bar y}_e}} \right\|_d} \le  1$ for all $e \in \bar{\mathcal{E}}_{\mathcal{R}}$. The task is split into two parts.
\\
{\bf{Part I:}} The linear system (\ref{E9-5-1}) admits a solution if and only if $\Psi$ lies in the range space of $\Omega \otimes I_n$. Since the left null space of the incidence matrix $\Omega$ is spanned by $1_{|\mathcal{R}|}$, this implies that $(1_{|\mathcal{R}|}^\top \otimes I_n) \Psi = 0_n$ must hold true, which expands to $\sum_{i \in \mathcal{R}} ( \nabla f_i(\tilde{x}^*) + \partial r_i(\tilde{x}^*) ) = 0_n$. This holds as it is the global optimality condition of the original \text{SRM} problem. Therefore, a valid set of subgradients $\partial r_i^*$ always exists such that $\Psi \in \text{Range}(\Omega \otimes I_n)$.
\\
{\bf{Part II:}} To bound the magnitude of the subgradients, we consider the minimum $\ell_2$-norm solution to (\ref{E9-5-1}) given by the Moore-Penrose pseudo-inverse: $\mathbf{y} = - (\Omega \otimes I_n)^\dagger \Psi / \phi$. According to the equivalence of $\ell_p$-norms in $\mathbb{R}^n$, we can explicitly bound the dual norm $\|\bar{y}_e\|_d$ using its $\ell_2$-norm via a dimension-dependent constant $C_{d,2} = n^{\max\left(0, \frac{1}{d} - \frac{1}{2}\right)}$. For instance, $C_{d,2} = 1$ for $p \in [1, 2]$, and $C_{d,2} = \sqrt{n}$ for $p = \infty$. For any specific edge $e \in \bar{\mathcal{E}}_{\mathcal{R}}$, we have the following sequence of inequalities:
\begin{equation}
\begin{aligned}
\left\| \bar{y}_e \right\|_d \le& C_{d,2} \left\| \bar{y}_e \right\|_2 \\
\le & C_{d,2} \left\| \mathbf{y} \right\|_2 \\
\le & \frac{C_{d,2}}{\phi} \left\| (\Omega \otimes I_n)^\dagger \right\|_2 \left\| \Psi \right\|_2 \\
\le & \frac{C_{d,2} \sqrt{|\mathcal{R}|}}{\phi {\sigma _{\min }}\left( \Omega  \right)} \max_{i \in \mathcal{R}} \left\| \psi_i \right\|_2,
\end{aligned}
\end{equation}
where ${\sigma _{\min }}\left( \Omega  \right)$ denotes the minimum non-zero singular value of the rectangular incidence matrix $\Omega$. To guarantee the subgradient validity condition $\left\| \bar{y}_e \right\|_d \le 1$ for all edges simultaneously, it is sufficient to enforce
\begin{equation}\label{E9-5-2}
\phi \ge \frac{C_{d,2} \sqrt{|\mathcal{R}|}}{{\sigma _{\min }}\left( \Omega  \right)} \max_{i \in \mathcal{R}} \left\| \nabla {f_i}\left( {{{\tilde x}^*}} \right) + \partial r_i^* \right\|_2.
\end{equation}
When the penalty parameter $\phi$ is chosen sufficiently large to satisfy (\ref{E9-5-2}), the constructed vector $\mathbf{y}$ yields a valid set of subgradients. This ensures that the global optimal solution to \text{SRM} (\ref{E3-3-1}) exactly meets the optimality condition (\ref{E8-1-1}) of the soft-approximation counterpart \text{DRC} (\ref{E3-3-2}), which completes the proof.

\subsection{Proof of Lemma \ref{L1}}\label{sec8-2}
Given the positive definiteness of $W_i$, for all $i \in \mathcal{R}$, it follows that
\begin{equation}\label{E8-2-2}
\begin{aligned}
&{\mathbb{E}_k}\left[{\left\| {\Delta _{i, k+1}^\psi} \right\|_{{W_i}}^2}\right]\\
\overset{(\ref{E4-2-2})}{=} &{\mathbb{E}_k} \left[{\left\| {{\Delta _{i,k}^\psi} + W_i^{ - 1}{M_{i,k}}\left( {\nabla {f_i}\left( {{x_{i,k}}} \right) - {h _{i,k}}} \right)} \right\|_{{W_i}}^2}\right] \\
= & {\mathbb{E}_k}\left[ {\left\| {\left( {{\mathbf{I}} - W_i^{ - 1}{M_{i,k}}} \right)\Delta _{i,k}^\psi } \right\|_{{W_i}}^2 + \left\| {W_i^{ - 1}{M_{i,k}}\Delta _{i,k}^f} \right\|_{{W_i}}^2} \right] + 2{\mathbb{E}_k}\left[ {{{\left\langle {\left( {{\mathbf{I}} - W_i^{ - 1}{M_{i,k}}} \right)\Delta _{i,k}^\psi ,W_i^{ - 1}{M_{i,k}}\Delta _{i,k}^f} \right\rangle }_{{W_i}}}} \right]\\
\overset{(\ref{E4-2-2+})}{=} & {{{\left( {\Delta _{i,k}^f} \right)}^ \top }{\mathbb{E}_k}\left[ {{{\left( {W_i^{ - 1}{M_{i,k}}} \right)}^ \top }{W_i}\left( {W_i^{ - 1}{M_{i,k}}} \right)} \right]\Delta _{i,k}^f} + {{{\left( {\Delta _{i,k}^\psi } \right)}^ \top}{\mathbb{E}_k}\left[ {{\left( {{\mathbf{I}} - W_i^{ - 1}{M_{i,k}}} \right) ^ \top }{W_i}\left( {{\mathbf{I}} - W_i^{ - 1}{M_{i,k}}} \right) } \right]\Delta _{i,k}^\psi }.\\
\end{aligned}
\end{equation}
Since $M_{i,k}^ \top  = {M_{i,k}}$ and ${\left( {W_i^{ - 1}} \right)^ \top } = W_i^{ - 1}$, we know that
\begin{equation}\label{E8-2-3}
\begin{aligned}
{\mathbb{E}}_k\left[ \left( {{\mathbf{I}} - W_i^{ - 1}{M_{i,k}}} \right)^ \top {W_i} \left( {{\mathbf{I}} - W_i^{ - 1}{M_{i,k}}} \right) \right] = & {\mathbb{E}}_k\left[ {\left( {{W_i} - {M_{i,k}}} \right)\left( {{\mathbf{I}} - W_i^{ - 1}{M_{i,k}}} \right)} \right]\\
\overset{(\ref{E4-2-2+})}{=} & {W_i} - {M_{i}}.
\end{aligned}
\end{equation}
Plugging (\ref{E8-2-3}) back into (\ref{E8-2-2}) completes the proof.

\subsection{Proof of Lemma \ref{L2}}\label{sec8-3}
For all $i \in \mathcal{R}$ and $k \ge 0$, the following inequality holds
\begin{equation}\label{E8-3-1}
{\mathbb{E}_k}\left[ {\left\| {{g_{i,k}} - \nabla {f_i}\left( {{{\tilde x}^*}} \right)} \right\|_{{W_i}}^2} \right] \overset{(\ref{E4-2-5})}{\le} 2{\mathbb{E}_k}\left[ {\left\| {{{\left( {{\mathbf{I}} - {\theta _{i,k}}W_i^{ - 1}{M_{i,k}}} \right)}\Delta _{i,k}^\psi}} \right\|_{{W_i}}^2} \right] + 2{\mathbb{E}_k}\left[ {\left\| {{\theta _{i,k}}W_i^{ - 1}{M_{i,k}}\Delta _{i,k}^f} \right\|_{{W_i}}^2} \right],
\end{equation}
which is owing to the basic inequality. We next handle ${\mathbb{E}_k}\left[ {\left\| {{{\left( {{\mathbf{I}} - {\theta _{i,k}}W_i^{ - 1}{M_{i,k}}} \right)}\Delta _{i,k}^\psi}} \right\|_{{W_i}}^2} \right]$ as follows:
\begin{equation}\label{E8-3-2}
\begin{aligned}
{\mathbb{E}_k}\left[ {\left\| {{{\left( {{\mathbf{I}} - {\theta _{i,k}}W_i^{ - 1}{M_{i,k}}} \right)}\Delta _{i,k}^\psi}} \right\|_{{W_i}}^2} \right] = &  {{{\left( {\Delta _{i,k}^\psi } \right)}^ \top }{\mathbb{E}_k}\left[ {\left\| {{\mathbf{I}} - {\theta _{i,k}}W_i^{ - 1}{M_{i,k}}} \right\|_{{W_i}}^2} \right]\Delta _{i,k}^\psi } \\
\overset{(\ref{E4-2-2+})}{=} & {{{\left( {\Delta _{i,k}^\psi } \right)}^ \top }{\mathbb{E}_k}\left[ {{W_i} - 2{\theta _{i,k}}{M_{i,k}} + \theta _{i,k}^2{M_{i,k}}} \right]\Delta _{i,k}^\psi }\\
\overset{(\ref{E4-2-4})}{=} &{\left( {\Delta _{i,k}^\psi } \right)^ \top } {\mathbb{E}_k}\left[ {{W_i} - 2{W_i} + \theta _{i,k}^2{M_{i,k}}} \right] \Delta _{i,k}^\psi\\
= & \left\| \Delta _{i,k}^\psi \right\|_{{\mathbb{E}_k}\left[ {\theta _{i,k}^2{M_{i,k}}} \right] - {W_i}}^2 \\
= & \left\| \Delta _{i,k}^\psi \right\|_{{\Phi _i}   - {W_i}}^2.
\end{aligned}
\end{equation}
To prove $\Phi_i - W_i \succcurlyeq 0$, we recall that for all $i \in \mathcal{R}$, $W_i^{-1} \succ 0$, $(W_i^{-1})^\top = W_i^{-1}$, and $M_{i,k}^\top = M_{i,k}$. Consequently, we have
\begin{equation}\label{E8-4-14}
\begin{aligned}
\mathbb{E}_k\left[ {\left\| {{\theta _{i,k}}{M_{i,k}}\tilde x - {W_i}\tilde x} \right\|_{W_i^{ - 1}}^2} \right] = & \mathbb{E}_k\left[ {\theta _{i,k}^2{{\tilde x}^ \top }{{\left( {{M_{i,k}}} \right)}^ \top }W_i^{ - 1}{M_{i,k}}\tilde x} \right] - 2{{\tilde x}^ \top }\mathbb{E}_k\left[ {{\theta _{i,k}}{M_{i,k}}} \right]\tilde x + {{\tilde x}^ \top }{W_i}\tilde x\\
\overset{(\ref{E4-2-2+})}{=} & {{\tilde x}^ \top }\mathbb{E}_k\left[ {\theta _{i,k}^2{M_{i,k}}} \right]\tilde x - 2{{\tilde x}^ \top }\mathbb{E}_k\left[ {{\theta _{i,k}}{M_{i,k}}} \right]\tilde x + {{\tilde x}^ \top }{W_i}\tilde x\\
\overset{(\ref{E4-2-4})}{=} & {{\tilde x}^ \top }\left( {\mathbb{E}_k\left[ {\theta _{i,k}^2{M_{i,k}}} \right] - {W_i}} \right)\tilde x\\
=& {{\tilde x}^ \top }\left( {{\Phi _i} - {W_i}} \right)\tilde x\\
 \ge & 0,
\end{aligned}
\end{equation}
which implies $\Phi_i - W_i \succcurlyeq 0$. We proceed to tackle ${\mathbb{E}_k}\left[ {\left\| {{\theta _{i,k}}W_i^{ - 1}{M_{i,k}}\Delta _{i,k}^f} \right\|_{{W_i}}^2} \right]$
\begin{equation}\label{E8-3-3}
\begin{aligned}
{\mathbb{E}_k}\left[ {\left\| {{\theta _{i,k}}W_i^{ - 1}{M_{i,k}}\Delta _{i,k}^f} \right\|_{{W_i}}^2} \right] = & {{{\left( {\Delta _{i,k}^f} \right)}^ \top }{\mathbb{E}_k}\left[ {\theta _{i,k}^2{{\left( {{M_{i,k}}} \right)}^ \top }W_i^{ - 1}{M_{i,k}}} \right]\Delta _{i,k}^f} \\
\overset{(\ref{E4-2-2+})}{=} & {\left( \Delta _{i,k}^f  \right)^ \top }{\Phi _i}\left( \Delta _{i,k}^f  \right)\\
= & \left\| {\Delta _{i,k}^f } \right\|_{\Phi_i}^2.
\end{aligned}
\end{equation}
Plugging (\ref{E8-3-2}) and (\ref{E8-3-3}) back into (\ref{E8-3-1}) completes the proof.

\subsection{Proof of Theorem \ref{T2}}\label{sec8-4}
Under the conditions given in Theorem \ref{T1}, we know that $x_i^* = {{\tilde x}^*}$ for all $i \in \mathcal{R}$ and the optimality condition of \text{DRC} (\ref{E3-3-2}) is given by
\begin{equation}\label{E8-4-1}
x_i^* = {\mathbf{prox}}_{\alpha ,{r_i}}^{{W_i}}\left\{ {x_i^* - \alpha \left( {\nabla {f_i}\left( {x_i^*} \right) + \partial \chi _i^*} \right)} \right\},\forall i \in \mathcal{R}.
\end{equation}
By invoking (\ref{E8-4-1}), we have that for all $i \in \mathcal{R}$ and $k \ge 0$,
\begin{equation}\label{E8-4-2}
\begin{aligned}
&\mathbb{E}_k\left[ {\left\| {{x_{i,k + 1}} - {x_i^*}} \right\|_{W_i}^2} \right]\\
\overset{(\ref{E4-3-3})}{=} &\mathbb{E}_k\left[ {\left\| {{{\mathbf{prox}}_{\alpha ,{r_i}}^{{W_i}}\left\{ {{{\bar x}_{i,k}}} \right\}}- {x_i^*}} \right\|_{W_i}^2} \right]\\
\overset{(\ref{E4-3-1})}{=} & \mathbb{E}_k\left[ \left\| {{\mathbf{prox}}_{\alpha ,{r_i}}^{{W_i}}\left\{ {{x_{i,k}} - \alpha {W_i^{ - 1}} \left( {{g_{i,k}} + \partial {\chi _{i,k}} + \partial {\delta _{i,k}}} \right)} \right\}} - {\mathbf{prox}}_{\alpha ,{r_i}}^{{W_i}}\left\{ {x_i^* - \alpha {W_i^{ - 1}} \left( {\nabla {f_i}\left( {{{\tilde x}^*}} \right) + \partial \chi _i^*} \right)} \right\} \right\|_{{W_i}}^2 \right]\\
\le & \mathbb{E}_k\left[ \left\| {x_{i,k}} - x_i^* -  \alpha W_i^{ - 1}\left( {{g_{i,k}} - \nabla {f_i}\left( {x_i^*} \right)} \right) - \alpha W_i^{ - 1}\left( {\partial {\chi _{i,k}} - \partial \chi _i^* + \partial {\delta _{i,k}}} \right) \right\|_{{W_i}}^2 \right]\\
= & \left\| {{x_{i,k}} - x_i^*} \right\|_{{W_i}}^2 - 2\alpha \mathbb{E}_k\left[ {{{\left\langle {{x_{i,k}} - x_i^*,{g_{i,k}} - \nabla {f_i}\left( {x_i^*} \right)} \right\rangle }}} \right] - 2\alpha \mathbb{E}_k\left[ {{{\left\langle {{x_{i,k}} - x_i^*,\partial {\chi _{i,k}} - \partial \chi _i^* + \partial {\delta _{i,k}}} \right\rangle }}} \right]\\
&+ {\alpha ^2}\mathbb{E}_k\left[ {\left\| {{g_{i,k}} - \nabla {f_i}\left( {x_i^*} \right) + \partial {\chi _{i,k}} - \partial \chi _i^* + \partial {\delta _{i,k}}} \right\|_{{W_i^{ - 1}}}^2} \right],
\end{aligned}
\end{equation}
where the inequality follows from the non-expansiveness of the weighted proximal operator \cite{Li2021}, which can be verified via the Cholesky decomposition of the positive-definite matrix $W_i$ \cite{Golub2013}. According to \cite[Proposition 1]{Hu2025a}, we have that for any vector $v \in {\mathbb{R}^n}$ and all $k \ge 0$,
\begin{equation}\label{E8-4-3}
\left| {{{\left[ {\partial {{\left\| {{x_{i,k}} - v} \right\|}_p}} \right]}_l}} \right| \le 1,\forall l = 1,2, \ldots ,n.
\end{equation}
Invoking (\ref{E8-4-3}), one can derive that for all $i \in \mathcal{R}$ and $k \ge 0$,
\begin{subequations}\label{E8-4-4}
\begin{align}
\label{E8-4-4-1}&\mathbb{E}_k\left[\left\| {\partial {\delta _{i,k}}} \right\|_{{W_i^{-1}}}^2\right] \le \frac{{n{\phi ^2}}}{{{\lambda _{\min }}\left( {{W_i}} \right)}}{\left| {{\mathcal{B}_i}} \right|^2},\\
\label{E8-4-4-2}&\mathbb{E}_k\left[\left\| {\partial {\chi _{i,k}} - \partial \chi _i^*} \right\|_{{W_i^{-1}}}^2\right] \le \frac{{4n{\phi ^2}}}{{{\lambda _{\min }}\left( {{W_i}} \right)}}{\left| {{\mathcal{R}_i}} \right|^2}.
\end{align}
\end{subequations}
We proceed to handle the last term in the right-hand-side (RHS) of (\ref{E8-4-2}) as follows:
\begin{equation}\label{E8-4-5}
\begin{aligned}
&{\alpha ^2}\mathbb{E}_k\left[ {\left\| {{g_{i,k}} - \nabla {f_i}\left( {x_i^*} \right) + \partial {\chi _{i,k}} - \partial \chi _i^* + \partial {\delta _{i,k}}} \right\|_{W_i^{-1}}^2} \right]\\
\le & 2{\alpha ^2}\mathbb{E}_k\left[ {\left\| {{g_{i,k}} - \nabla {f_i}\left( {x_i^*} \right)} \right\|_{{W_i^{-1}}}^2} \right]  + 4{\alpha ^2}\mathbb{E}_k\left[\left\| {\partial {\delta _{i,k}}} \right\|_{{W_i^{-1}}}^2\right] + 4{\alpha ^2}\mathbb{E}_k\left[\left\| {\partial {\chi _{i,k}} - \partial \chi _i^*} \right\|_{{W_i^{-1}}}^2\right] \\
\overset{(\ref{E8-4-4})}{\le} & \frac{{2{\alpha ^2}}}{{{\lambda _{\min }}{{\left( {{W_i}} \right)}^2}}}\mathbb{E}_k\left[ {\left\| {{g_{i,k}} - \nabla {f_i}\left( {x_i^*} \right)} \right\|_{{W_i}}^2} \right] + \frac{{16n{\phi ^2}{\alpha ^2}}}{{{\lambda _{\min }}\left( {{W_i}} \right)}}{\left| {{\mathcal{R}_i}} \right|^2}+ \frac{{4n{\phi ^2}{\alpha ^2}}}{{{\lambda _{\min }}\left( {{W_i}} \right)}}{\left| {{\mathcal{B}_i}} \right|^2} \\
\overset{(\ref{E5-2-2})}{\le} & \frac{{4{\alpha ^2}}}{{{\lambda _{\min }}{{\left( {{W_i}} \right)}^2}}} {\left\| {\Delta _{i,k}^\psi } \right\|_{{\mathbb{E}_k}\left[ {\theta _{i,k}^2{M_{i,k}}} \right] - {W_i}}^2 + \left\| {\Delta _{i,k}^f} \right\|_{{\Phi _i}}^2} + \frac{{4n{\phi ^2}{\alpha ^2}}}{{{\lambda _{\min }}\left( {{W_i}} \right)}}\left( {4 {{{\left| {{\mathcal{R}_i}} \right|}^2}} + {{{\left| {{\mathcal{B}_i}} \right|}^2}} } \right),
\end{aligned}
\end{equation}
where the first inequality applies the basic inequality twice. We continue to handle the second term in the RHS of (\ref{E8-4-2}). For arbitrary two variables $\tilde x,\tilde y \in {\mathbb{R}^n}$, we define the Bregman divergence with respect to the local objective function $f_i$ by ${D_{{f_i}}}\left( {\tilde x,\tilde y} \right) := {f_i}\left( {\tilde x} \right) - {f_i}\left( {\tilde y} \right) - {\left\langle {\nabla {f_i}\left( y \right),\tilde x - \tilde y} \right\rangle}$, which is nonnegative as $f_i$ is convex. Since the local objective is $\mu_i$-strongly convex in $W_i$-norm, we have
\begin{equation}\label{E8-4-6}
\begin{aligned}
- 2\alpha\mathbb{E}_k\left[\left\langle {{x_{i,k}} - x_i^*,{g_{i,k}} - \nabla {f_i}\left( {x_i^*} \right)} \right\rangle \right] = & - 2\alpha\left\langle {{x_{i,k}} - x_i^*,\Delta _{i,k}^f} \right\rangle\\
\overset{(\ref{E3-1-2})}{\le} & - {\mu _i}\alpha \left\| {{x_{i,k}} - x_i^*} \right\|_{{W_i}}^2 - 2 \alpha {D_{{f_i}}}\left( {{x_{i,k}},x_i^*} \right)\\
\overset{(\ref{E3-1-3})}{\le} & - {\mu _i}\alpha\left\| {{x_{i,k}} - x_i^*} \right\|_{{W_i}}^2 - \alpha\left\| {\Delta _{i,k}^f} \right\|_{{Q_i}}^2.
\end{aligned}
\end{equation}
We next handle the third term in the RHS of (\ref{E8-4-2}). Recalling the definition of ${\partial {\delta _{i,k}}}$, we have that for any $\varphi  > 0 $
\begin{equation}\label{E8-4-7}
\begin{aligned}
-2\alpha\mathbb{E}_k\left[{\left\langle {{x_{i,k}} - x_i^*,\partial {\delta _{i,k}}} \right\rangle } \right] \le & \varphi  \alpha \left\| {{x_{i,k}} - x_i^*} \right\|_{{W_i}}^2 + \frac{\alpha}{\varphi  }\mathbb{E}_k\left[\left\| {\partial {\delta _{i,k}}} \right\|_{{W_i^{-1}}}^2 \right]\\
\overset{(\ref{E8-4-4-1})}{\le} & \varphi \alpha \left\| {{x_{i,k}} - x_i^*} \right\|_{{W_i}}^2 + \frac{{n\alpha {\phi ^2}}}{{\varphi {\lambda _{\min }}\left( {{W_i}} \right)}}{\left| {{\mathcal{B}_i}} \right|^2},
\end{aligned}
\end{equation}
where the first inequality applies the Young's inequality. Similarly, recalling the definition of ${\chi _{i,k}}$, we have that for any $\eta > 0 $,
\begin{equation}\label{E8-4-8}
\begin{aligned}
-2\alpha\mathbb{E}_k\left[\left\langle {{x_{i,k}} - x_i^*,\partial {\chi _{i,k}} - \partial \chi _i^*} \right\rangle \right] \le &\eta \alpha \left\| {{x_{i,k}} - x_i^*} \right\|_{{W_i}}^2 + \frac{{\alpha}}{\eta } \mathbb{E}_k\left[\left\| {\partial {\chi _{i,k}} - \partial \chi _i^*} \right\|_{{W_i ^ {-1}}}^2\right]\\
\overset{(\ref{E8-4-4-2})}{\le} & \eta \alpha \left\| {{x_{i,k}} - x_i^*} \right\|_{{W_i}}^2 + \frac{{4n\alpha {\phi ^2}}}{{\eta {\lambda _{\min }}\left( {{W_i}} \right)}}{\left| {{\mathcal{R}_i}} \right|^2}.
\end{aligned}
\end{equation}
Plugging (\ref{E8-4-5})-(\ref{E8-4-8}) back into (\ref{E8-4-2}) obtains
\begin{equation}\label{E8-4-9}
\begin{aligned}
\mathbb{E}_k\left[ {\left\| {{x_{i,k + 1}} - {x_i^*}} \right\|_{W_i}^2} \right]= & \left( {1 - \left( {{\mu _i} - \varphi  - \eta } \right)\alpha } \right) {\left\| {{x_{i,k}} - x_i^*} \right\|_{{W_i}}^2} - \alpha \left\| {\Delta _{i,k}^f} \right\|_{{Q_i}}^2 + \frac{{4{\alpha ^2}}}{{{\lambda _{\min }}{{({W_i})}^2}}} \left\| {\Delta _{i,k}^\psi } \right\|_{{\Phi _i} - {W_i}}^2 + \left\| {\Delta _{i,k}^f} \right\|_{{\Phi _i}}^2 \\
& + \frac{{4n{\phi ^2}}}{{{\lambda _{\min }}({W_i})}}\left( {4{{\left| {{\mathcal{R}_i}} \right|}^2} + {{\left| {{\mathcal{B}_i}} \right|}^2}} \right){\alpha ^2} + \frac{{n{\phi ^2}}}{{{\lambda _{\min }}({W_i})}}\left( {\frac{4}{\eta }{{\left| {{\mathcal{R}_i}} \right|}^2} + \frac{1}{\varphi }{{\left| {{\mathcal{B}_i}} \right|}^2}} \right)\alpha
\end{aligned}
\end{equation}
We recall the definitions: $\mu : = {\min _{i \in \mathcal{R}}}{\mu _i}$, $\tilde{\Phi} = \text{diag}\{\lambda_{\min}(W_i)^{-2}\Phi_i\}_{i\in\mathcal{R}}$ and $\tilde{W} := \text{diag}\{\lambda_{\min}(W_i)^{-2}W_i\}_{i\in\mathcal{R}}$. By summing both sides of (\ref{E8-4-9}) over $i \in \mathcal{R}$ and rearranging the result, we obtain
\begin{equation}\label{E8-4-10}
\begin{aligned}
\mathbb{E}_k\left[ {\left\| {{x_{k + 1}} - {x^*}} \right\|_W^2} \right] = & \left( {1 - \left( {\mu  - \varphi   - \eta } \right)\alpha } \right)\left\| {{x_k} - {x^*}} \right\|_W^2 + \left\| {\Delta _k^F} \right\|_{4{\alpha ^2}{\tilde \Phi} - \alpha Q}^2 + 4{\alpha ^2} \left\| {\Delta _k^\psi } \right\|_{{\tilde \Phi} - \tilde W}^2  \\
&+ 4n{\phi ^2}{\alpha ^2}\sum\limits_{i \in \mathcal{R}} {\frac{{4{{\left| {{\mathcal{R}_i}} \right|}^2} + {{\left| {{\mathcal{B}_i}} \right|}^2}}}{{{\lambda _{\min }}\left( {{W_i}} \right)}}} + n\alpha {\phi ^2}\sum\limits_{i \in \mathcal{R}} {\frac{1}{{{\lambda _{\min }}\left( {{W_i}} \right)}}\left( {\frac{4}{\eta }{{\left| {{\mathcal{R}_i}} \right|}^2} + \frac{1}{\varphi }{{\left| {{\mathcal{B}_i}} \right|}^2}} \right)}.
\end{aligned}
\end{equation}
For simplicity, we then set $\eta  = \varphi   = \gamma /2$, such that (\ref{E8-4-10}) becomes
\begin{equation}\label{E8-4-11}
\begin{aligned}
\mathbb{E}_k\left[ {\left\| {{x_{k + 1}} - {x^*}} \right\|_W^2} \right] \le & \left( {1 - \gamma  \alpha } \right)\left\| {{x_k} - {x^*}} \right\|_W^2 + \left\| {\Delta _k^F} \right\|_{4{\alpha ^2}{\tilde \Phi} - \alpha Q}^2 + 4{\alpha ^2} \left\| {\Delta _k^\psi } \right\|_{{\tilde \Phi} - \tilde W}^2\\
& + 4n{\phi ^2}{\alpha ^2}\sum\limits_{i \in \mathcal{R}} {\frac{1}{{{\lambda _{\min }}\left( {{W_i}} \right)}}\left( {4{{\left| {{\mathcal{R}_i}} \right|}^2} + {{\left| {{\mathcal{B}_i}} \right|}^2}} \right)} \\
&  + n\alpha {\phi ^2}\sum\limits_{i \in \mathcal{R}} {\frac{1}{{{\lambda _{\min }}\left( {{W_i}} \right)}}\left( {\frac{4}{\eta }{{\left| {{\mathcal{R}_i}} \right|}^2} + \frac{1}{\varphi }{{\left| {{\mathcal{B}_i}} \right|}^2}} \right)}.
\end{aligned}
\end{equation}
According to the gradient-learning sequence (\ref{E5-2-1}), we have for any constant $c > 0$,
\begin{equation}\label{E8-4-12}
c\alpha \mathbb{E}_k\left[ {\left\| {\Delta _{k + 1}^\psi } \right\|_W^2} \right] = c\alpha \left\| {\Delta _k^\psi } \right\|_{W - M}^2 + c\alpha \left\| {\Delta _k^F} \right\|_{M}^2.
\end{equation}
Recalling the definitions of $\gamma$, $\Gamma _1$, and $\Gamma _2$, we combine (\ref{E8-4-11})-(\ref{E8-4-12}) to obtain
\begin{equation}\label{E8-4-13}
\begin{aligned}
\mathbb{E}_k\left[ {U_{k + 1}}\right] = & \mathbb{E}_k\left[ {\left\| {{x_{k + 1}} - {x^*}} \right\|_W^2} \right] + c\alpha \mathbb{E}_k\left[ {\left\| {\Delta _{k + 1}^\psi } \right\|_W^2} \right]\\
\le & \left( {1 - \gamma  \alpha } \right)U_k^W + \left\| {\Delta _k^F} \right\|_{{\text{4}}{\alpha ^2}{\Phi} - \alpha \left( {Q - cM} \right)}^2 + \left\| {\Delta _k^\psi } \right\|_{\left( {\left( {\gamma  c - 4} \right)W + 4{\Phi}} \right){\alpha ^2} - c\alpha M}^2 + {\Gamma _1}{\alpha ^2} + {\Gamma _2}\alpha.
\end{aligned}
\end{equation}
For $i \in \mathcal{R}$, since $S_i$ is not a zero matrix or zero vector and $M_i$ and $\Phi_i$ are positive-definite matrices. For any $k \ge 0$, if we choose $0 < c < {{\lambda _{\min }}\left( Q \right)}/\left( {2{\lambda _{\max }}\left( {M} \right)} \right)$ and $0 < \alpha  \le {{\lambda _{\min }}\left( {Q - cM} \right)}/\left( 4{\lambda _{\max }}\left( \tilde \Phi  \right) \right)$, then we can drop the second term in the RHS of (\ref{E8-4-13}). From Lemma \ref{L2}, we know that ${\Phi} - {W} \succcurlyeq 0$ such that ${\tilde \Phi} - {\tilde W} \succcurlyeq 0$. For any $k \ge 0$, if we choose $0 < \alpha  \le {c \lambda _{\min }}\left( {M} \right) / {\lambda _{\max }}\left( {\gamma  cW + 4\left( {\tilde \Phi  - \tilde W} \right)} \right)$, then the third term in the RHS of (\ref{E8-4-13}) can also be dropped. By summarizing all the upper bounds on the constant step-size, i.e., $0 < \alpha  \le \min \left\{ {\frac{{c{\lambda _{\min }}\left( {M} \right)}}{{{\lambda _{\max }}\left( {\gamma  cW + 4\left( {\tilde \Phi  - \tilde W} \right)} \right)}},\frac{{{\lambda _{\min }}\left( {Q - cM} \right)}}{{4{\lambda _{\max }}\left(\tilde \Phi  \right)}}} \right\}$, then (\ref{E8-4-13}) reduces to
\begin{equation}\label{E8-4-16}
\mathbb{E}_k\left[{U_{k + 1}}\right] \le \left( {1 - \gamma  \alpha } \right)U_k^W + {\Gamma _1}{\alpha ^2} + {\Gamma _2}\alpha.
\end{equation}
Applying the telescopic cancellation to (\ref{E8-4-16}) over $k$ and taking the total expectation on the both sides of the result completes the proof.

\subsection{Proof of Corollary \ref{C1}}\label{sec8-5}
We define $\nabla F\left( {{x^*}} \right): = {\mathbf{col}}{\left\{ {\nabla {f_i}\left( {{{\tilde x}^*}} \right)} \right\}_{i \in \mathcal{R}}}$, $\nabla F\left( {{x_k}} \right): = {\mathbf{col}}{\left\{ {\nabla {f_i}\left( {{x_{i,k}}} \right)} \right\}_{i \in \mathcal{R}}}$, $\partial  {\chi _k}: = {\mathbf{col}}{\left\{ \partial {\chi _{i,k}} \right\}_{i \in \mathcal{R}}}$, and $\partial {\delta _k}: = {\mathbf{col}}{\left\{ {\partial {\delta _{i,k}}} \right\}_{i \in \mathcal{R}}}$. Since the weight matrices are assumed to be $W_i = \mathbf{I}$ for all $i \in \mathcal{R}$, the proof can be facilitated by the \textit{RED-SEGA} update in the following compact form:
\begin{subequations}\label{E8-5-1}
\begin{align}
\label{E8-5-1-1}&{{{\bar x}_{k}}} = {x_k} - {\alpha _k}\left( {{g_k} + \partial {\chi _k} + \partial {\delta _k}} \right),\\
\label{E8-5-1-2}&{x_{k + 1}} = {\mathbf{prox}}_{{\alpha _k},R}\left\{ {{{\bar x}_k}} \right\},
\end{align}
\end{subequations}
where ${{\bar x}_{k}}: = {\mathbf{col}}{\left\{ {{\bar x}_{i,k}} \right\}_{i \in \mathcal{R}}}$ and $R: = \sum\nolimits_{i \in \mathcal{R}} {{r_i}} $. Building upon Lemmas \ref{L1} and \ref{L2}, summing both sides of (\ref{E5-2-1}) and (\ref{E5-2-2}) over $i \in \mathcal{R}$ under the condition $W_i = \mathbf{I}$ yields the following two inequalities:
\begin{subequations}\label{E8-5-2}
\begin{align}
\label{E8-5-2-1} & \mathbb{E}_k\left[ {{{\left\| {\Delta _{k + 1}^\psi } \right\|}^2}} \right] =  \left\| {\Delta _k^\psi } \right\|_{{\mathbf{I}} - M}^2 + \left\| {\Delta _k^F} \right\|_M^2,\\
\label{E8-5-2-2} & \mathbb{E}_k\left[ {\left\| {{g_k} - \nabla F\left( {{x^*}} \right)} \right\|_W^2} \right] \le 2\left\| {\Delta _k^\psi } \right\|_{\Phi  - \mathbf{I}}^2 + 2\left\| {\Delta _k^F} \right\|_\Phi ^2.
\end{align}
\end{subequations}
According to the optimality condition: ${x^*} = {\mathbf{pro}}{{\mathbf{x}}_{{\alpha _k},R}}\left\{ {{x^*} - {\alpha _k}\left( {\nabla F\left( {{x^*}} \right) + \partial {\chi ^*}} \right)} \right\}$, we define ${\Delta _k}: = {g_k} - \nabla F\left( {{x^*}} \right) + \partial {\chi _k} - \partial {\chi ^*} + \partial {\delta _k}$ such that
\begin{equation}\label{E8-5-3}
\begin{aligned}
&\mathbb{E}_k\left[ {{{\left\| {{x_{k + 1}} - {x^*}} \right\|}^2}} \right]\\
= & {\mathbb{E}_k}\left[ {{{\left\| {{\mathbf{pro}}{{\mathbf{x}}_{{\alpha _k},R}}\left\{ {{x_k} - {\alpha _k}\left( {{g_k} + \partial {\chi _k} + \partial {\delta _k}} \right)} \right\} - {\mathbf{pro}}{{\mathbf{x}}_{{\alpha _k},R}}\left\{ {{x^*} - {\alpha _k}\left( {\nabla F\left( {{x^*}} \right) + \partial {\chi ^*}} \right)} \right\}} \right\|}^2}} \right]
\le &\mathbb{E}_k\left[ {{\left\| {{x_k} - {x^*} - {\alpha _k}\Delta _k} \right\|^2}} \right]\\
= & {\left\| {{x_k} - {x^*}} \right\|^2} + \alpha _k^2\mathbb{E}_k\left[ \left\|\Delta _k \right\|^2 \right] - 2{\alpha _k}\mathbb{E}_k\left[ {\left\langle {{x_k} - {x^*}, \Delta _k} \right\rangle } \right],
\end{aligned}
\end{equation}
where the inequality applies the non-expansiveness of the proximal operator ${\mathbf{prox}}_{{\alpha _k},R}\left\{ \cdot \right\}$. According to (\ref{E8-4-3}), one can derive
\begin{subequations}\label{E8-5-4}
\begin{align}
\label{E8-5-4-1} & \mathbb{E}_k\left[{\left\| {\partial {\delta _k}} \right\|^2}\right] \le n{\phi ^2}\sum\limits_{i \in \mathcal{R}} {\left| {{\mathcal{B}_i}} \right|}^2,\\
\label{E8-5-4-2} & \mathbb{E}_k\left[{\left\| {\partial {\chi _k} - \partial {\chi ^*}} \right\|^2} \right] \le 4n{\phi ^2}\sum\limits_{i \in \mathcal{R}} {\left| {{\mathcal{R}_i}} \right|}^2.
\end{align}
\end{subequations}
By utilizing (\ref{E8-5-4-1})-(\ref{E8-5-4-2}), we proceed to handle the second term in the RHS of (\ref{E8-5-3}) as follows:
\begin{equation}\label{E8-5-5}
\begin{aligned}
\mathbb{E}_k\left[ {{\left\| \Delta _k \right\|^2}} \right] \le & 2\mathbb{E}_k\left[ {\left\| {{g_k} - \nabla F\left( {{x^*}} \right)} \right\|^2 + 2\left\| {\partial {\chi _k} - \partial {\chi ^*}} \right\|^2 + 2\left\| {\partial {\delta _k}} \right\|^2} \right]\\
\le & 2\mathbb{E}_k\left[ {{{\left\| {{g_k} - \nabla F\left( {{x^*}} \right)} \right\|}^2}} \right] + 16n{\phi ^2}\sum\limits_{i \in \mathcal{R}} {\left( {{{\left| {{\mathcal{R}_i}} \right|}^2} + {{\left| {{\mathcal{B}_i}} \right|}^2}} \right)} \\
\le & 4{\left\| {\Delta _k^\psi } \right\|_{\Phi  - {\mathbf{I}}}^2 + \left\| {\Delta _k^F} \right\|_\Phi ^2} + 16n{\phi ^2}\sum\limits_{i \in \mathcal{R}} {\left( {{{\left| {{\mathcal{R}_i}} \right|}^2} + {{\left| {{\mathcal{B}_i}} \right|}^2}} \right)},
\end{aligned}
\end{equation}
where the first inequality uses the basic inequality twice and the last inequality is owing to (\ref{E8-5-2-2}). To proceed, we proceed to handle the last term in the RHS of (\ref{E8-5-3}) as follows:
\begin{equation}\label{E8-5-6}
\begin{aligned}
& \mathbb{E}_k\left[ {\left\langle {{x_k} - {x^*},{g_k} - \nabla F\left( {{x^*}} \right) + \partial {\chi _k} - \partial {\chi ^*} + \partial {\delta _k}} \right\rangle } \right] \\
= & \mathbb{E}_k\left[ {\left\langle {{x_k} - {x^*},{g_k} - \nabla F\left( {{x^*}} \right)} \right\rangle } \right] + \mathbb{E}_k\left[\left\langle {{x_k} - {x^*},\partial {\chi _k} - \partial {\chi ^*}} \right\rangle\right] + \mathbb{E}_k\left[\left\langle {{x_k} - {x^*},\partial {\delta _k}} \right\rangle\right] \\
= &\mathbb{E}_k\left[ {\left\langle {{x_k} - {x^*},{g_k} - \nabla F\left( {{x_k}} \right) + \nabla F\left( {{x_k}} \right) - \nabla F\left( {{x^*}} \right)} \right\rangle } \right] + \mathbb{E}_k\left[ \left\langle {{x_k} - {x^*},\partial {\chi _k} - \partial {\chi ^*}} \right\rangle\right] + \mathbb{E}_k\left[ \left\langle {{x_k} - {x^*},\partial {\delta _k}} \right\rangle  \right]\\
= & \mathbb{E}_k\left[ {\left\langle {{x_k} - {x^*},\nabla F\left( {{x_k}} \right) - \nabla F\left( {{x^*}} \right)} \right\rangle } \right] + \mathbb{E}_k\left[\left\langle {{x_k} - {x^*},\partial {\delta _k}} \right\rangle \right] + \mathbb{E}_k\left[\left\langle {{x_k} - {x^*},\partial {\chi _k} - \partial {\chi ^*}} \right\rangle \right]\\
\end{aligned}
\end{equation}
where the last equality follows from the fact that $\mathbb{E}_k\left[ {\left\langle {{x_k} - {x^*},{g_k} - \nabla F\left( {{x_k}} \right)} \right\rangle } \right] = 0$. Different to \cite{Hu2025a}, $\partial {\chi _k}$ is based on the current graph whereas $\partial\chi^*$ comes from the average graph, which implies that we cannot directly follows the subgradient monotonicity to handle the last term in the RHS of (\ref{E8-5-6}). However, recalling the definition of $\bar {P}(x)$, we know that $\partial\chi^*$ is inherently a valid subgradient of $\bar{P}(x)$ evaluated at the global optimal solution $x^*$, i.e., $\partial {\chi ^*} \in \bar P\left( {{x_k}} \right)$. Let $\mathcal{F}_k$ denote the filtration generated by the historical randomness up to iteration $k$. Under Assumption \ref{A3}, taking the conditional expectation of the instantaneous subgradient $\partial\chi_k$ with respect to the network randomness $\zeta_k$ yields the expected average graph subgradient $\bar \partial {\chi _k}: = {\mathbb{E}_{{\zeta _k}}}\left[ {\partial {\chi _k}|{\mathcal{F}_k}} \right] \in \partial \bar P\left( {{x_k}} \right)$. Since $\bar{\partial}\chi_k$ and $\partial\chi^*$ belong to the subdifferential of the identical convex function $\bar{P}(\cdot)$, invoking the law of total expectation and the standard subgradient monotonicity allows us to bound the cross-term
\begin{equation}\label{E8-5-7}
\begin{aligned}
- \mathbb{E}_k\left[\left\langle {{x_k} - {x^*},\partial {\chi _k} - \partial {\chi ^*}} \right\rangle \right] \le & - {\mathbb{E}_k}\left[ {{\mathbb{E}_{{\zeta _k}}}\left[ {\left\langle {{x_k} - {x^*},\partial {\chi _k} - \partial {\chi ^*}} \right\rangle }|{\mathcal{F}_k} \right]} \right]\\
= &  - \mathbb{E}_k\left[ {\left\langle {{x_k} - {x^*},\bar \partial {\chi _k} - \partial {\chi ^*}} \right\rangle } \right]\\
\le & 0,
\end{aligned}
\end{equation}
which implies that the last term in the RHS of (\ref{E8-5-6}) can dropped. Given any constant $\gamma  > 0$, we next handle the second term in the RHS of (\ref{E8-5-6}) as follows:
\begin{equation}\label{E8-5-8}
- 2\mathbb{E}_k\left[\left\langle {{x_k} - {x^*},\partial \delta_k} \right\rangle\right] \le \gamma \left\| {{x_k} - {x^*}} \right\|^2 + \frac{n}{\gamma }{\phi ^2}\sum\limits_{i \in \mathcal{R}} {{{\left| {{\mathcal{B}_i}} \right|}^2}},
\end{equation}
where we apply the Young's inequality. Since $F$ is a strongly-convex and smooth function, we have
\begin{equation}\label{E8-5-9}
\begin{aligned}
- 2\left\langle {x_k} - {x^*},\nabla F\left( {{x_k}} \right) - \nabla F\left( {{x^*}} \right) \right\rangle \le & - \mu \left\| {{x_k} - {x^*}} \right\|^2 - \left\| {\nabla F\left( {{x_k}} \right) - \nabla F\left( {{x^*}} \right)} \right\|_Q^2.
\end{aligned}
\end{equation}
Plugging (\ref{E8-5-9})-(\ref{E8-5-7}) back into (\ref{E8-5-6}) reduces to
\begin{equation}\label{E8-5-10}
\begin{aligned}
- {\text{2}}\mathbb{E}_k\left[ {\left\langle {{x_k} - {x^*},{g_k} - \nabla F\left( {{x^*}} \right) + \partial {\chi _k} - \partial {\chi ^*} + \partial {\delta _k}} \right\rangle } \right] \le & \left( {\gamma  - \mu } \right){\left\| {{x_k} - {x^*}} \right\|^2} - \left\| {\nabla F\left( {{x_k}} \right) - \nabla F\left( {{x^*}} \right)} \right\|_Q^2 + \frac{n}{\gamma }{\phi ^2}\sum\limits_{i \in \mathcal{R}} {{{\left| {{\mathcal{B}_i}} \right|}^2}}.
\end{aligned}
\end{equation}
We recall the definitions of ${\Delta _k^\psi }$ and $\Delta _k^F$, such that substituting ((\ref{E8-5-5})) and ((\ref{E8-5-10})) back into (\ref{E8-5-3}) gives
\begin{equation}\label{E8-5-11}
\begin{aligned}
\mathbb{E}_k\left[ {{{\left\| {{x_{k + 1}} - {x^*}} \right\|}^2}} \right] \le & \left( {1 - \gamma {\alpha _k}} \right){\left\| {{x_k} - {x^*}} \right\|^2} + 4\alpha _k^2\left\| {\Delta _k^\psi } \right\|_{\Phi  - {\mathbf{I}}}^2  - 4\left\| {\Delta _k^F} \right\|_{{\alpha _k}Q - 4\Phi \alpha _k^2}^2 + \frac{n}{\gamma }{\phi ^2}\sum\limits_{i \in \mathcal{R}} {{{\left| {{\mathcal{B}_i}} \right|}^2}}\\
& + 16n{\phi ^2}\alpha _k^2\sum\limits_{i \in \mathcal{R}} {\left( {{{\left| {{\mathcal{R}_i}} \right|}^2} + {{\left| {{\mathcal{B}_i}} \right|}^2}} \right)}.
\end{aligned}
\end{equation}
According to (\ref{E8-5-2-1}), we have for any constant $\tilde c > 0$
\begin{equation}\label{E8-5-12}
\tilde c {\alpha _k}\mathbb{E}_k\left[ {{{\left\| {\Delta _{k + 1}^\psi } \right\|}^2}} \right] = \tilde c {\alpha _k}\left\| {\Delta _k^\psi } \right\|_{{\mathbf{I}} - M}^2 + \tilde c {\alpha _k}\left\| {\Delta _k^F} \right\|_M^2.
\end{equation}
Adding (\ref{E8-5-12}) to the both sides of (\ref{E8-5-11}) yields
\begin{equation}\label{E8-5-13}
\begin{aligned}
\mathbb{E}_k\left[ {{U_{k + 1}}} \right]  = & \mathbb{E}_k\left[ {{\left\| {{x_{k + 1}} - {x^*}} \right\|}^2} \right]+ \mathbb{E}_k\left[ \tilde c{\alpha _k}{{\left\| {\Delta _{k + 1}^\psi } \right\|}^2} \right] \\
\le & \left( {1 - \gamma {\alpha _k}} \right){\left\| {{x_k} - {x^*}} \right\|^2} + \left\| {\Delta _k^F} \right\|_{{\alpha _k}Q - 4\alpha _k^2\Phi }^2 + 4\alpha _k^2\left\| {\Delta _k^\psi } \right\|_{\Phi  - {\mathbf{I}}}^2 + \alpha _k^2 {\Gamma _3}  + \alpha _k {\Gamma _4} + \tilde c{\alpha _k}\left\| {\Delta _k^\psi } \right\|_{{\mathbf{I}} - M}^2 + \tilde c{\alpha _k}\left\| {\Delta _k^F} \right\|_M^2\\
= & \left( {1 - \gamma {\alpha _k}} \right){{U_k}} + \left\| {\Delta _k^\psi } \right\|_{\left( {\left( {\gamma \tilde c - 4} \right) + 4\Phi } \right)\alpha _k^2 - \tilde c{\alpha _k}M}^2 +  \left\| {\Delta _k^F} \right\|_{4\alpha _k^2\Phi  - {\alpha _k}\left( {Q - \tilde c M} \right)}^2 + \alpha _k^2 {\Gamma _3}  + \alpha _k {\Gamma _4}
\end{aligned}
\end{equation}
If the decaying step-size is chosen as $0 < {\alpha _k} \le {\lambda _{\min }}\left( {Q - \tilde cM} \right)/\left( {4{\lambda _{\max }}\left( \Phi  \right)} \right)$ when $0 < \tilde c < {\lambda _{\min }}\left( Q \right)/\left( {2{\lambda _{\max }}\left( M \right)} \right)$, we can drop the third term in the RHS of (\ref{E8-5-13}). Recalling the definitions of $\Gamma _3$ and $\Gamma _4$, (\ref{E8-5-13}) reduces to
\begin{equation}\label{E8-5-14}
\begin{aligned}
\mathbb{E}_k\left[ {{U_{k + 1}}} \right] \le & \left( {1 - \gamma {\alpha _k}} \right){U_k} + \left\| {\Delta _k^\psi } \right\|_{\left( {\left( {\gamma \tilde c - 4} \right) + 4\Phi } \right)\alpha _k^2 - \tilde c{\alpha _k}M}^2 + \alpha _k^2 {\Gamma _3}  + \alpha _k {\Gamma _4}.
\end{aligned}
\end{equation}
Following (\ref{E8-4-14}), one can verify that $\Phi  - {\mathbf{I}} \succcurlyeq 0$. If we further set the decaying step-size by $0 < {\alpha _k} \le \tilde c {\lambda _{\min }}\left( M \right)/{\lambda _{\max }}\left( {\gamma \tilde c + {\text{4}}\left( {\Phi  - {\mathbf{I}}} \right)} \right)$, then the second term in the RHS of (\ref{E8-5-13}) can also be dropped such that taking the total expectation on the both sides of (\ref{E8-5-14}) reduces to
\begin{equation}\label{E8-5-15}
\mathbb{E}\left[ {{U_{k + 1}}} \right] \le \left( {1 - \gamma {\alpha _k}} \right)\mathbb{E}\left[{U_k} \right] + \alpha _k^2 {\Gamma _3}  + \alpha _k {\Gamma _4}.
\end{equation}
Summarizing two upper bounds on the decaying step-size gives
\begin{equation}\label{E8-5-16}
0 < {\alpha _k} \le \min \left\{ {\frac{{{\lambda _{\min }}\left( {Q - \tilde c M} \right)}}{{4{\lambda _{\max }}\left( \Phi  \right)}},\frac{{\tilde c{\lambda _{\min }}\left( M \right)}}{{{\lambda _{\max }}\left( {\gamma \tilde c + 4\left( {\Phi  - {\mathbf{I}}} \right)} \right)}}} \right\}.
\end{equation}
For a sufficiently large positive constant $\Theta$, we need to prove
\begin{equation}\label{E8-5-16+}
 \mathbb{E}\left[ {{U_k}} \right] \le \Xi /\left( {\xi  + k} \right) + \Theta, \forall k \ge 0,
\end{equation}
 by induction. For $k=0$, we know from (\ref{E8-5-14}) that
\begin{equation}\label{E8-5-17}
{U_1} \le \left( {1 - \frac{\mu }{2}{\alpha _0}} \right){U_0} + \alpha _0^2{\Gamma _3} + {\alpha _0}{\Gamma _4}.
\end{equation}
We set ${\alpha _k} = \beta  /\left( {k + \xi } \right)$ with ${\alpha _0} = \beta  /\xi $. If $\Xi  \ge \left( {\xi  - \gamma \beta  } \right){U_0} + \left( {{\beta  ^2}{\Gamma _3}/\xi } \right) + \beta  {\Gamma _4} - \xi \Theta $ with $0 < \beta  \gamma  \le \xi $, the following inequality holds true
\begin{equation}\label{E8-5-18}
\left( {1 - \gamma {\alpha _0}} \right){U_0} + \alpha _0^2{\Gamma _3} + {\alpha _0}{\Gamma _4} \le \frac{\Xi }{\xi } + \Theta.
\end{equation}
We assume that for $k = k'$, $\forall k' \ge 1$, the following inequality holds
\begin{equation}\label{E8-5-19}
\begin{aligned}
\mathbb{E}\left[ {{U_{k' + 1}}} \right] \le & \left( {1 - \gamma {\alpha _{k'}}} \right) {\mathbb{E}\left[U_{k'}\right]}  + \alpha _{k'}^2{\Gamma _3} + {\alpha _{k'}}{\Gamma _4}\\
\le & \frac{\Xi }{{\xi  + k'}} + \Theta.
\end{aligned}
\end{equation}
Then, we will prove that for $k = k' + 1$, the following inequality still holds true
\begin{equation}\label{E8-5-20}
\begin{aligned}
\mathbb{E}\left[ {{U_{k' + 2}}} \right] \le & \left( {1 - \gamma {\alpha _{k' + 1}}} \right)\mathbb{E}\left[ {{U_{k' + 1}}} \right] + \alpha _{k' + 1}^2{\Gamma _3} + {\alpha _{k' + 1}}{\Gamma _4}\\
 \le & \frac{\Xi }{{\xi  + k' + 1}} + \Theta.
\end{aligned}
\end{equation}
By choosing $\Theta  \ge{\Gamma _4}/\gamma $ and $\xi > \beta  \max \{ \gamma, \varpi \}$ with $\beta  > 1/\gamma$, we obtain
\begin{equation}\label{E8-5-21}
\begin{aligned}
\mathbb{E}\left[ {{U_{k' + 2}}} \right] \le &\left( {1 - \gamma {\alpha _{k' + 1}}} \right)\mathbb{E}\left[ {{U_{k' + 1}}} \right] + \alpha _{k' + 1}^2{\Gamma _3} + {\alpha _{k' + 1}}{\Gamma _4}\\
= & \left( {1 - \frac{{\gamma \beta  }}{{\xi  + k' + 1}}} \right)\mathbb{E}\left[ {{U_{k' + 1}}} \right] + \frac{{{\beta  ^2}{\Gamma _3}}}{{{{\left( {\xi  + k' + 1} \right)}^2}}} + \frac{{\beta  {\Gamma _4}}}{{\xi  + k' + 1}}\\
\le & \left( {1 - \frac{{\gamma \beta  }}{{\xi  + k' + 1}}} \right)\left( {\frac{\Xi }{{\xi  + k'}} + \Theta } \right)  + \frac{{{\beta  ^2}{\Gamma _3}}}{{{{\left( {\xi  + k' + 1} \right)}^2}}} + \frac{{\beta  {\Gamma _4}}}{{\xi  + k' + 1}}\\
= & \left( {1 - \frac{{\gamma \beta  }}{{\xi  + k' + 1}}} \right)\frac{\Xi }{{\xi  + k'}} + \Theta  + \frac{{{\beta  ^2}{\Gamma _3}}}{{{{\left( {\xi  + k' + 1} \right)}^2}}} + \frac{\beta  }{{\xi  + k' + 1}}\left( {{\Gamma _4} - \gamma \Theta } \right)\\
\le & \left( {1 - \frac{{\gamma \beta  }}{{\xi  + k' + 1}}} \right)\frac{\Xi }{{\xi  + k'}} + \Theta  + \frac{{\Xi \left( {\gamma \beta   - 1} \right)}}{{{{\left( {\xi  + k' + 1} \right)}^2}}}\\
\le & \frac{\Xi }{{\xi  + k'}} - \frac{\Xi }{{\left( {\xi  + k'} \right)\left( {\xi  + k' + 1} \right)}} + \Theta \\
= & \frac{\Xi }{{\xi  + k' + 1}} + \Theta,
\end{aligned}
\end{equation}
where the second inequality applies (\ref{E8-5-19}) and the third inequality is owing to the facts that both $\Theta  \ge{\Gamma _4}/\gamma $ and $\beta  > 1/\gamma$ are selected. This confirms that (\ref{E8-5-16+}) holds. The proof is completed by fixing $\Theta  = {\Gamma _4}/\gamma$ and then taking the total expectation on the both sides of the result.

\vfill

\end{document}